\documentclass[11pt, reqno]{amsart}
\usepackage{lmodern}
\usepackage{amsmath, amsthm, amssymb, amsfonts}
\usepackage[normalem]{ulem}
\usepackage{hyperref}
\usepackage[all,cmtip]{xy}
\usepackage{verbatim}
\usepackage{nccmath}
\usepackage{stmaryrd}
\usepackage{caption}
\usepackage{mathrsfs}
\usepackage{mathtools}
\usepackage{esvect}
\usepackage{cite}
\usepackage{bbm}
\usepackage{eucal}

\usepackage{mathbbol}
\usepackage{tabularx}
\usepackage[toc,page]{appendix}

\usepackage{tikz-cd}

\theoremstyle{plain}
\newtheorem{thm}{Theorem}[section]
\newtheorem{cor}[thm]{Corollary}

\newtheorem{lemma}[thm]{Lemma}
\newtheorem{prop}[thm]{Proposition}
\newtheorem{conj}[thm]{Conjecture}
\newtheorem*{conjno}{Conjecture}

\newtheorem{conjecture}{Conjecture}  

\renewcommand{\arraystretch}{2}

\theoremstyle{definition}
\newtheorem{defn}[thm]{Definition}

\makeatletter
\newcommand\ackname{Acknowledgements}
\if@titlepage
\newenvironment{acknowledgements}{%
	\titlepage
	\null\vfil
	\@beginparpenalty\@lowpenalty
	\begin{center}%
		\bfseries \ackname
		\@endparpenalty\@M
\end{center}}%
{\par\vfil\null\endtitlepage}
\else

\fi
\makeatother

\theoremstyle{remark}
\newtheorem{rmk}[thm]{Remark}

\newcommand{\BA}{{\mathbb{A}}}

\newcommand{\BC}{{\mathbb{C}}}

\newcommand{\BE}{{\mathbb{E}}}

\newcommand{\BG}{{\mathbb{G}}}

\newcommand{\BL}{{\mathbb{L}}}

\newcommand{\BN}{{\mathbb{N}}}

\newcommand{\BP}{{\mathbb{P}}}
\newcommand{\BQ}{{\mathbb{Q}}}
\newcommand{\BR}{{\mathbb{R}}}

\newcommand{\BT}{{\mathbb{T}}}

\newcommand{\BZ}{{\mathbb{Z}}}

\newcommand{\CA}{{\mathcal A}}

\newcommand{\CC}{{\mathcal C}}
\newcommand{\CD}{{\mathcal D}}
\newcommand{\CE}{{\mathcal E}}
\newcommand{\CF}{{\mathcal F}}
\newcommand{\CG}{{\mathcal G}}
\newcommand{\CH}{{\mathcal H}}
\newcommand{\CI}{{\mathcal I}}

\newcommand{\CK}{{\mathcal K}}
\newcommand{\CL}{{\mathcal L}}
\newcommand{\CM}{{\mathcal M}}
\newcommand{\CN}{{\mathcal N}}
\newcommand{\CO}{{\mathcal O}}
\newcommand{\CP}{{\mathcal P}}
\newcommand{\CQ}{{\mathcal Q}}

\newcommand{\CS}{{\mathcal S}}

\newcommand{\CU}{{\mathcal U}}

\newcommand{\CX}{{\mathcal X}}
\newcommand{\CY}{{\mathcal Y}}

\newcommand{\FM}{{\mathfrak{M}}}

\newcommand{\rT}{{\mathrm{T}}}
\newcommand{\vsf}{{\mathsf{v}}}
\newcommand{\wsf}{{\mathsf{w}}}

\newcommand{\Qsf}{{\mathsf{Q}}}
\newcommand{\rsf}{{\mathsf{r}}}
\newcommand{\dsf}{{\mathsf{d}}}
\newcommand{\usf}{{\mathsf{u}}}
\newcommand{\asf}{{\mathsf{a}}}
\newcommand{\msf}{{\mathsf{m}}}
\newcommand{\Asf}{{\mathsf{A}}}

\newcommand{\Fsf}{{\mathsf{F}}}
\newcommand{\Usf}{{\mathsf{U}}}
\newcommand{\Lsf}{{\mathsf{L}}}
\newcommand{\Ws}{{\check{W}}}
\newcommand{\Wp}{{\hat{W}}}

\newcommand{\ch}{{\mathrm{ch}}}

\newcommand{\rk}{{\mathrm{rk}}}

\newcommand{\Ms}{\check{M}}
\newcommand{\Mp}{\hat{M}}
\newcommand{\Msp}{\hat{\check{M}}}
\newcommand{\VWs}{\check{\mathsf{VW}}}
\newcommand{\VWp}{\hat{\mathsf{VW}}}

\newcommand{\QMs}{\check{\mathsf{QM}}}
\newcommand{\QMp}{\hat{\mathsf{QM}}}
\newcommand{\GWs}{\check{\mathsf{GW}}}
\newcommand{\GWp}{\hat{\mathsf{GW}}}
\newcommand{\QMsp}{\hat{\check{\mathsf{QM}}}}
\newcommand{\SW}{\mathsf{SW}}
\newcommand{\Ns}{\check{N}}

\newcommand{\Nsp}{\hat{\check{N}}}
\newcommand{\Thetas}{\check{\Theta}}
\newcommand{\Thetap}{\hat{\Theta}}

\newcommand{\alphas}{\check{\alpha}}
\newcommand{\alphap}{\hat{\alpha}}

\newcommand{\tr}{{\mathrm{tr}}}
\newcommand{\CHom}{{\mathcal{H} om}}

\DeclareFontFamily{OT1}{rsfs}{}
\DeclareFontShape{OT1}{rsfs}{n}{it}{<-> rsfs10}{}
\DeclareMathAlphabet{\curly}{OT1}{rsfs}{n}{it}

\newcommand\Ext{\operatorname{Ext}}
\newcommand\Hom{\operatorname{Hom}}
\newcommand{\p}{\mathbb{P}}

\newcommand\Spec{\operatorname{Spec}}
\newcommand{\Jac}{{\mathrm{Jac}}}

\newcommand{\Mbar}{{\overline M}}

\newcommand{\td}{\mathrm{td}}

\newcommand{\Pic}{\mathop{\rm Pic}\nolimits}

\newcommand{\DT}{\mathsf{DT}}

\newcommand{\ev}{{\mathrm{ev}}}
\newcommand{\evsf}{{\mathsf{ev}}}

\newcommand{\thickslash}{\mathbin{\!\!\pmb{\fatslash}}}

\newcommand{\FMs}{\check{\mathfrak{M}}}
\newcommand{\FMsr}{\check{\mathfrak{M}}_{\mathrm{rg}}}
\newcommand{\FMp}{\hat{\mathfrak{M}}}
\newcommand{\FMsp}{\hat{\check{\mathfrak{M}}}}
\newcommand{\HAz}{\hat{\mathfrak{M}}}

\newcommand{\Gsf}{\mathsf{G}}

\begin{document}
	
	\title[Enumerative mirror symmetry  and S-duality]
	{Enumerative mirror symmetry for moduli spaces of Higgs bundles and S-duality}
	
	\author{Denis Nesterov}
	\address{University of Vienna, Faculty of Mathematics}
	\email{denis.nesterov@univie.ac.at}
	\maketitle
	\begin{abstract} We derive conjectures, called genus 1 Enumerative mirror symmetry for moduli spaces of Higgs bundles, which relate curve-counting invariants of moduli spaces of Higgs $\mathrm{SL}_\rsf$-bundles to curve-counting invariants of moduli spaces of Higgs $\mathrm{PGL}_\rsf$-bundles. This contrasts with Enumerative mirror symmetry for Calabi-Yau 3-folds which relates curve-counting invariants to periods. We also provide extensive mathematical evidence for these conjectures. 
		
		The conjectures are obtained with the help of the theory of quasimaps to moduli spaces of sheaves, Tanaka--Thomas's construction of Vafa--Witten theory, Jiang--Kool's enumerative S-duality of Vafa--Witten invariants and Manschot--Moore's calculations.  We use the latter together with some basic computations to give a complete list of conjectural expressions for genus 1 quasimap invariants for all prime ranks. They have many interesting properties, among which is quantum $\chi$-independence. The wall-crossing to Gromov--Witten theory is also thoroughly discussed.

	\end{abstract}
	\setcounter{tocdepth}{1}
	\tableofcontents
	\section{Introduction}
	\subsection{Overview} 
\begin{figure} [h!]	\vspace{2.3cm}
	\scriptsize
	\[
	\begin{picture}(200,75)(-30,-50)
		\thicklines
		\put(25,-25){\line(1,0){25}}
		\put(100,-25){\line(1,0){25}}
		\put(25,-25){\line(0,1){40}}
		\put(25,26){\makebox(0,0){\textsf{S-duality}}}
		\put(25,35){\line(0,1){40}}
		\put(125,-25){\line(0,1){40}}
		\put(125,26){\makebox(0,0){\textsf{T-duality}}}
		\put(25,75){\line(1,0){25}}
		\put(100,75){\line(1,0){25}}
		\put(125,35){\line(0,1){40}}
		\put(140,85){\makebox(0,0){$\sigma$-model of $M_G$}}
		\put(15,85){\makebox(0,0){Yang-Mills for $G$ on $X\times C$}}
		\put(75,75){\makebox(0,0){\textsf{dim. red.}}}
		\put(15,-35){\makebox(0,0){Yang-Mills for $G^L$ on $X\times C$}}
		\put(145,-35){\makebox(0,0){$\sigma$-model of $M_{G^L}$}}
		\put(75,-25){\makebox(0,0){\textsf{dim. red.}}}
	\end{picture}
	\]
	\caption{Dimensional reduction}
	\label{square}
	\vspace{-0.3cm}
\end{figure}
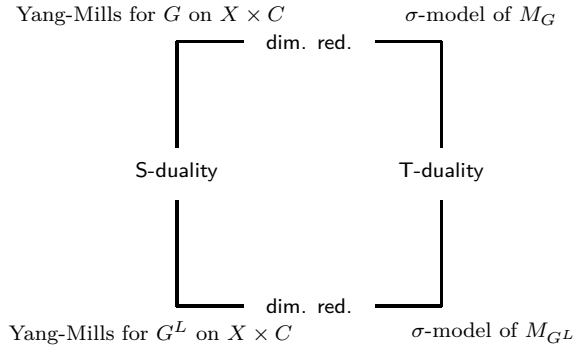

In \cite{KW}, Kapustin and Witten used \textit{dimensional reduction} to relate a twisted 4-dimensional $N$=4 super  Yang--Mills theory on a product of two curves $X\times C$ to a twisted 2-dimensional $\sigma$-model of a moduli space of Higgs bundles on the curve $X$. With respect to this dimensional reduction, S-duality between 4-dimensional theories associated to a group $G$ and its Langlands dual $G^L$ reduces to T-duality between $\sigma$-models of moduli spaces of Higgs bundles associated to $G$ and $G^L$, as it is shown in Figure \ref{square}. 

In this article, we produce a mathematical realisation of the dimensional reduction on the level of enumerative invariants (intersection theory on moduli spaces of maps and sheaves) for moduli spaces of Higgs $\mathrm{SL}_\rsf$- and $\mathrm{PGL}_\rsf$-bundles on $X$, denoted by $\Ms$ and $\Mp$, respectively. 

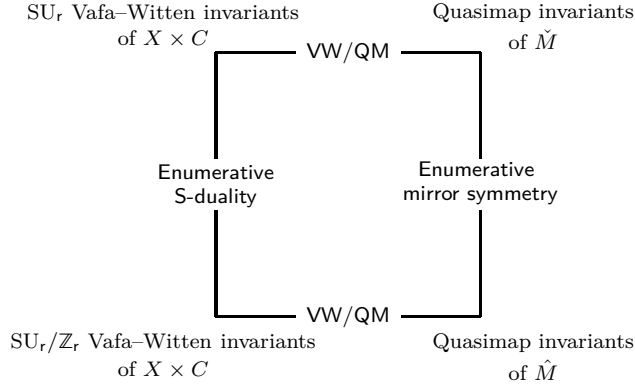
\begin{figure} [h!]	\vspace{2.5cm}
	\scriptsize
	\[
	\begin{picture}(200,75)(-30,-50)
		\thicklines
		\put(25,-25){\line(1,0){30}}
		\put(95,-25){\line(1,0){30}}
		\put(25,-25){\line(0,1){40}}
		\put(25,30){\makebox(0,0){\textsf{Enumerative}}}
		\put(25,20){\makebox(0,0){\textsf{S-duality}}}
		\put(25,35){\line(0,1){40}}
		\put(125,-25){\line(0,1){40}}
		\put(125, 31){\makebox(0,0){\textsf{Enumerative }}}
		\put(125,21){\makebox(0,0){\textsf{mirror symmetry}}}
		\put(25,75){\line(1,0){30}}
		\put(95,75){\line(1,0){30}}
		\put(125,35){\line(0,1){40}}
		\put(145,90){\makebox(0,0){Quasimap invariants}}
		\put(145,80){\makebox(0,0) {of $\Ms$}}
		\put(5,90){\makebox(0,0){$\mathrm{SU}_\rsf$ Vafa--Witten invariants}}
		\put(5,80){\makebox(0,0){of $X\times C$}}
		\put(75,75){\makebox(0,0){\textsf{VW/QM}}}
		\put(5,-35){\makebox(0,0){$\mathrm{SU}_\rsf/\BZ_\rsf$ Vafa--Witten invariants}}
		\put(5,-45){\makebox(0,0){of $X\times C$}}
		\put(145,-35){\makebox(0,0){Quasimap invariants}}
		\put(145,-45){\makebox(0,0){of $\Mp$}}
		\put(75,-25){\makebox(0,0){\textsf{VW/QM}}}
	\end{picture}
	\]
	\caption{Vafa--Witten/Quasimaps correspondence}
	\label{square2}
\end{figure}

For us, the twisted 4-dimensional $N$=4 super Yang--Mills theory turns into the Vafa--Witten theory in the sense of Tanaka--Thomas \cite{TT,TT2}; the twisted 2-dimensional $\sigma$-model turns into the theory of quasimaps to moduli spaces of sheaves \cite{N,NK3}. The dimensional reduction then becomes a correspondence between Vafa--Witten invariants and quasimap invariants, summarised in Theorem \ref{vwqm} for genus 1 invariants. With respect to this correspondence, Jiang--Kool's enumerative S-duality for Vafa--Witten invariants \cite{JK} translates into a statement about quasimap invariants, which we call \textit{Enumerative mirror symmetry for moduli spaces of Higgs bundles} - Conjecture \ref{MS1} and \ref{Enumerativemirror}. Overall,  Figure \ref{square} becomes Figure \ref{square2}.

Using physical calculations of \cite{MM}, we also give explicit formulas for quasimap invariants - Conjecture \ref{wnonzero} and \ref{wzero}.  
Here, we focus on genus 1 quasimap invariants, i.e.\ $g(C)=1$. 
 However, our correspondence between quasimaps and Higgs sheaves holds for arbitrary genera and comes in two flavours - quasimaps of \cite{N} or quasisections of \cite[Section 2]{Nqm}. If $g(C)=1$, there is no need to distinguish between these two approaches, as they coincide. 

The basic calculations are done in Section \ref{calc}, where we prove the conjectures for $g(X)=1$, which amounts to determining Vafa--Witten invariants for abelian surfaces, whose physical derivation is given in \cite[(5.25)]{LL}. We also determine certain degree 0 invariants for all genera from topological invariants, Proposition \ref{ggreat2}. Some parts of Conjecture \ref{wnonzero} and \ref{wzero} will be proven for all genera in the sequel paper \cite{Nqm}, where we also calculate some higher rank invariants. These calculations allow us to guess natural extensions of conjectures to an arbitrary prime rank in Section \ref{higherrank}. The resulting statements satisfy all basic constraints and expectations. Our conjectures completely determine genus 1 invariants for $\Ms$ and $\Mp$ in all prime ranks and have many interesting properties. Among them is independence from degree of Higgs bundles, which can be called \textit{quantum $\chi$-independence}, see Conjecture \ref{independence}.

By adopting results of \cite{N,NK3} and \cite{Ob}, we also show that quasimap invariants admit a wall-crossing to Gromov--Witten invariants, such that the wall-crossing invariants are given by integrals on Quot schemes on the curve $X$, Theorem \ref{wallcrossing} and \ref{wallcrossing2}. In many cases, Gromov--Witten and quasimap invariants are equal, Corollary \ref{gweqqm} and \ref{gweqqm2}. The wall-crossing will also be central for calculations in \cite{Nqm}. 

	\begin{figure} [h!]	\vspace{2cm}
	\scriptsize
	\[
	\begin{picture}(200,0)(-30,-30)
		\thicklines
		\put(25,0){\line(1,0){30}}
		\put(95,0){\line(1,0){30}}
		\put(75,5){\makebox(0,0){\textsf{quasimap}}}
			\put(75,-5){\makebox(0,0){\textsf{wall-crossing}}}
		\put(150,10){\makebox(0,0){Quasimap invariants}}
		\put(150,0){\makebox(0,0) {of $\Ms$ or $\Mp$}}
		\put(0,10){\makebox(0,0){Gromov--Witten invariants}}
		\put(0,0){\makebox(0,0){of $\Ms$ or $\Mp$}}
	\end{picture}
	\]
	\caption{Quasimap wall-crossing}
	\label{square3}
\end{figure}
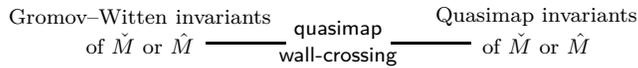
	\subsection{Topological and Categorical mirror symmetries}Hyper-K\"ahler nature of moduli spaces of Higgs bundles makes their mirror symmetry immensely rich. In its most elevated form, mirror symmetry of moduli spaces of Higgs bundles incorporates  Geometric Langlands conjectures, as explained by Kapustin--Witten \cite{KW}. But even in its more classical mathematical incarnations, it exhibits very special features (in comparison with the Calabi--Yau 3-fold situation). More precisely, if the complex structure is fixed to be the one for which Hitchin fibration is holomorphic, both topological and categorical invariants of dual moduli spaces of Higgs bundles, $\Ms$ and $\Mp$, agree without being exchanged, i.e.\ there is no conventional A-model/B-model swap\footnote{We refer to \cite[Table 2]{KW} for more details on how A-models and B-models are exchanged for various complex structures of the hyper-K\"ahler triplet.}. For example, this can be seen in Hausel--Thaddeus's Topological mirror symmetry \cite{HT}, which asserts equality of Hodge numbers without reflections along the diagonal line of Hodge diamonds as in the Calabi--Yau 3-fold situation.
	\[\textbf{Topological mirror symmetry:} \quad h^{p,q}(\Ms)=h_{\mathrm{orb},\alphap}^{p,q}(\Mp).\]
	 The same phenomenon is also present in Categorical mirror symmetry of Donagi--Pantev \cite{DP}  - generically over the base of the Hitchin map, there is an equivalence of bounded derived categories of coherent sheaves with no Fukaya categories involved.
	\[\hspace{-0.6cm}\textbf{Categorical mirror symmetry:} \quad D^{b}(\Ms)\approx D^{b}(\Mp).\]
	This leaves us wander what happens with curve-counting invariants (like Gromov--Witten invariants) of $\Ms$ and $\Mp$. In other words, what is the statement of Enumerative mirror symmetry for  $\Ms$ and $\Mp$? 
		\[\hspace{-0.8cm}\textbf{Enumerative mirror symmetry:} \quad \quad \quad \quad ?\quad \quad  \quad  .\]
	Despite an extensive  progress in both Topological and Categorical mirror symmetries (e.g.\ the former was proven in \cite{GWZ}, while latter was shown to hold generically in \cite{DP}), nothing has been said about Enumerative mirror symmetry.
	
	 In a nutshell, our statement of Enumerative mirror symmetry also relates curve-counting invariants to curve-counting invariants without a mention of B-model invariants. However, it has an interesting feature of being a relation that involves moduli spaces of Higgs bundles of all degrees, even of degree 0 (which is an Artin stack).
	

\subsection{Enumerative mirror symmetry} Genus 1 curve-counting invariants provide a good testing ground for Enumerative mirror symmetry, because other genera require non-trivial insertions. Moreover, genus 1 invariants determine the trace of an operator given by 2-point genus 0 invariants. 

 Let us now present a statement of genus 1 Enumerative mirror symmetry conjecture, for that we firstly need a few definitions. 
\subsubsection{Moduli spaces} Let $\rsf$ and $\dsf$ be a prime number\footnote{For the most part of the present work, this assumption is unnecessary, but it will significantly simply the exposition. The main part which requires this assumption is Section \ref{azal}} and an integer, such that $\rsf\geq 2$ and $0 \leq \dsf<\rsf$. We define $\FMsr(\dsf)$ to be a rigidified\footnote{We quotient out the $\BC^*$-automorphisms of sheaves given by multiplication of a scalar, see Section \ref{SLbundles} for more details.} moduli stack of rank $\rsf$ and degree $\dsf$ Higgs sheaves with a fixed determinant and a traceless Higgs fields on a smooth projective curve $X$. Let 
\[\Ms(\dsf) \subset \FMsr(\dsf)\]
 be the locus of semistable Higgs sheaves. Let 
 \[\Gamma_X:= \Jac(X)[\rsf]\]
  be the group of $\rsf$-torsion line bundles on $X$, it acts on $\FMs_{\mathrm{rg}}(\dsf)$ by tensoring sheaves with line bundles. We define the dual $\mathrm{PGL}_\rsf$-spaces
  \[\FMp(\dsf):= [\FMsr(\dsf)/\Gamma_X] \quad \text{and} \quad \Mp(\dsf):=[\Ms(\dsf)/\Gamma_X].\] 
  
  We will need two classes, the class of the theta line bundle (the ample generator of the Picard group) on $\FMsr(\dsf)$ and the class of the  $\BZ_\rsf$-gerbe of  SL-liftings of the universal family (the class used by \cite{HT}) on $\FMsr(\dsf)$, denoted as 
  \[ \Thetas \in H^2(\FMsr(\dsf),\BQ) \quad \text{and} \quad \alphas \in H^2(\FMsr(\dsf), \BZ_\rsf), \]
  respectively. Both classes descend to classes on $\FMp(\dsf)$, 
  \[ \Thetap \in H^2(\FMp(\dsf),\BQ) \quad \text{and} \quad \alphap \in H^2(\FMp(\dsf), \BZ_\rsf). \]
  \subsubsection{Quasimaps, $\dsf\neq0$} Assume for the moment that $\dsf\neq 0$.  Let $E$ be an elliptic curve.  A quasimap of degree $(\wsf,\asf) \in \BZ \oplus \BZ_\rsf$  from $E$ to $\Ms(\dsf)$ is a map
  \[f\colon E \rightarrow  \FMsr(\dsf),\]
  such that 
  \begin{itemize} \label{defconditions}
  	\item $f$ generically maps to the stable locus $\Ms(\dsf) \subset \FMsr(\dsf)$,
  	\item $f^*\Thetas=\wsf \in H^2(E,\BZ)\cong \BZ$,
  	\item $f^*\alphas=\asf \in H^2(E, \BZ_\rsf)\cong \BZ_\rsf$. 
  \end{itemize}
  We define quasimaps of degree $(\wsf, \asf)$ to $\Mp(\dsf)$ in exactly the same way, using the classes $\Thetap$ and $\alphap$. For $\wsf\neq0$,  let  
  \[ Q_{E}(\Ms(\dsf), \asf,\wsf)^{\bullet} \quad \text{and} \quad Q_{E}(\Mp(\dsf), \asf, \wsf)^{\bullet}\] 
  be the moduli spaces of quasimaps from $E$ to $\Ms(\dsf)$ and $\Mp(\dsf)$ up to translations\footnote{We identify quasimaps, if they equal up to a translation of $E$.} of $E$. 
  The (reduced) expected dimension of these spaces is 0, so we can define a virtual number of quasimaps from $E$, 
  
  \begin{align*}
  	& \QMs^{\asf, \bullet}_{\dsf,\wsf}:= \int_{[Q_{E}(\Ms(\dsf),\asf,\wsf)^{\bullet}]^{\mathrm{vir}}} 1\in \BQ,\\
  	& \QMp^{\asf, \bullet}_{\dsf,\wsf}:= \int_{[Q_{E}(\Mp(\dsf),\asf, \wsf)^{\bullet}]^{\mathrm{vir}}} 1\in \BQ.
  	 \end{align*}
We secretly use virtual localisation to define the integrals above, we refer to Definitions \ref{qmsell} and \ref{qmsell2} for the precise definitions. 
 
 \subsubsection{Quasimaps, $\dsf=0$}If $\dsf=0$, then $\Ms(\dsf)$ is an Artin stack. Counting curves in Artin stacks is problematic and there exists no general approach to it so far.  In this work, we take the following route, which is in accordance with \cite[Section 6]{Wit}. Using the relation between quasimaps and sheaves, which was explored in \cite{N} and adjusted to this setting in Section \ref{sectionvw}, we can treat quasimaps as semistable sheaves.  Hence we say that a quasimap $f\colon E \rightarrow \FMs(0)$ is \textit{semistable}, if 
 \begin{itemize}
 	\item the associated Higgs sheaf on $X\times E$ is semistable with respect to a polarisation of the form $\CO_X(1)\boxtimes \CO_E(k)$ for $k\gg 1$. 
 	\end{itemize}
By Proposition \ref{stablefiber}, the condition above is automatic for $\dsf\neq 0$ and is independent of $k$, as long as $k\gg 1$. A similar definition can be given in the $\mathrm{PGL}_\rsf$  case.  We then use the theory of Joyce--Song \cite{JS} and Joyce \cite{J} to define the quasimap invariants for $\Ms(0)$ and $\Mp(0)$, 
 \[\QMs^{\asf, \bullet}_{0,\wsf},  \QMp^{\asf, \bullet}_{0,\wsf} \in \BQ,\] 
 we refer to Section \ref{noncoprime} for more details. 
 
 \subsubsection{Conjecture} \label{conjsec}
 Consider the formal generating series in $\BQ[\![q]\!]$, 
 \begin{align*}
 	\QMs^{\asf}_{\dsf}(q)&:= \sum_{\wsf>0} \QMs^{\asf, \bullet}_{\dsf, \mathsf w}q^{\mathsf w}\\
 	\QMp^{\asf}_{\dsf}(q)&:= \sum_{\wsf>0} \QMp^{\asf, \bullet}_{\dsf, \mathsf w}q^{\mathsf w}.
 \end{align*}
Let $\tilde{\eta}(q):=\prod^{\infty}_{k=1}(1-q^k)$, we then define
\begin{align*}
\Usf_1(q)&:=\log \tilde{\eta}(q^4)\\
\Usf_2(q)&:=  \log \tilde{\eta}(q)\\
\Usf_3(q)&:=  \log\tilde{\eta}(-q).
\end{align*}
Let the symmetric group $S_3$ act on the linear span $ \BQ\langle \Usf_i(q) \rangle $ by permutation, 
\[ \sigma \cdot \Usf_i(q):=\Usf_{\sigma(i)}(q), \quad \sigma \in S_3.\]
\begin{conjno}[Conjecture \ref{wnonzero}, \ref{MS1}] If $g \geq2$ and $\rsf=2$, then
	\begin{align*}
		\QMs^0_{\dsf}(q)&=(-1)^\dsf(2-2g)2^{4g-1}  \Usf_1(q) \\	
		\QMs^1_{\dsf}(q)&=(2-2g) 2^{2g-1 }(\Usf_2(q)+(-1)^\dsf\Usf_3(q)) \\
	\\
		\QMp^0_{\dsf}(q)&=(-1)^{\dsf}(2-2g)2^{2g-1}\Usf_1(q) \\
		\QMp^1_{\dsf}(q)&=(2-2g) (2^{4g-1}  \Usf_2(q) +(-1)^{\dsf} 2^{2g-1}  \Usf_3(q)),
	\end{align*}
such that the series above satisfy the following relation, called Enumerative mirror symmetry,
\[(12)\cdot 2\QMs^{\asf}_{\dsf}(q)=\sum_{\dsf'} \sum_{\mathsf a'} (-1)^{  \mathsf d \cdot \mathsf a'+\dsf'\cdot \asf' } \hat{\mathsf{QM}}^{\asf'}_{\dsf'}(q).\]

\end{conjno}

We refer to Conjecture \ref{wnonzeror} for a statement in the case of an arbitrary prime rank $\rsf$. The conjecture is derived using the correspondence between quasimaps and Vafa--Witten invariants discussed in Section \ref{sectionvw} and \ref{noncoprime}, which allows us to translate physical calculations done in \cite[Section 7]{MM} into the language of quasimaps. Additional features of the conjecture come from the structural results on the resulting generating series, see Section \ref{sduality} for more details. 

The corresponding statements of enumerative S-duality are given in  Section \ref{sduality}. Enumerative mirror symmetry for $\wsf=0$ and any prime rank $\rsf$ is given in Conjecture \ref{Enumerativemirror}.  Since degree 0 genus 1 invariants recover Euler characteristics, we obtain a peculiar statement involving Euler characteristics of moduli spaces of Higgs bundles, whose geometric meaning is obscured by the presence of invariants of $\hat M(0)$. This statement differs from \cite{HT}.
\subsection{The series $\Usf_i$}Let us now say a few words about the meaning of $\Usf_i$. We will write $\Msp(\dsf)$ to mean either $\Mp(\dsf)$ or $\Mp(\dsf)$.    Consider  the $\BC_t^*$-action on $\Msp(\dsf)$ which scales Higgs fields. We have a decomposition of the $\BC_t^*$-fixed loci, 
\begin{align*}
	\Msp(\dsf)^{\BC_t^*}&=\Nsp(\dsf) \cup \Msp(\dsf)^{\mathrm{nil}}, 
\end{align*}
where $\Nsp(\dsf)$ are moduli spaces of semistable bundles;  $\Msp(\dsf)^{\mathrm{nil}}$ are loci of semistable Higgs bundles with nilpotent Higgs fields. Then $\Usf_1(q)$ corresponds to invariants of $\Msp(\dsf)^{\mathrm{nil}}$, while $\Usf_2(q)$ and $\Usf_3(q)$ correspond to invariants of $\Nsp(\dsf) $, 

\begin{align*}
	\BQ \langle \Usf_1(q) \rangle & \iff \Msp(\dsf)^{\mathrm{nil}}  \\
	\BQ \langle \Usf_2(q), \Usf_3(q) \rangle  & \iff  \Nsp(\dsf).
\end{align*}
More precisely, projections of  $\QMsp^{\asf}_{\dsf}(q)$ to $\BQ\langle \Usf_1(q) \rangle$ or to $\BQ\langle \Usf_2(q),\Usf_3(q)\rangle$  are exactly the invariants of the corresponding components of the fixed loci. 
Hence Enumerative mirror symmetry partially exchanges the invariants of $\Nsp(\dsf)$ and of $\Msp(\dsf)^{\mathrm{nil}}$.  This is not surprising in the light of S-duality, because these invariants correspond to \textit{monopole} and \textit{instanton} branches of the Vafa--Witten theory. 

\subsection{Other formulation of Enumerative mirror symmetry} The conjecture presented above has a nice property of being a statement about formal generating series without analytic continuations.  Let us  now present a different formulation of Enumerative mirror symmetry that involves a more traditional change-of-variables relation. We define 
\begin{align*}
	\QMs^{\asf}_{\dsf}(q)'&:=q\frac{d\QMs^{\asf}_{\dsf}(q)}{dq} +\check{\mathrm{c}}\\
	\QMp^{\asf}_{\dsf}(q)'&:= q\frac{d\QMp^{\asf}_{\dsf}(q)}{dq} + \hat{\mathrm{c}},
\end{align*}
where constants\footnote{The constants appear in the formulas of \cite{MM}, but it is not clear how to define them mathematically, e.g.\ they are not present in the formulas of \cite[Appendix C]{GK};   they should correspond to $\wsf=0$ invariants, but the construction of these invariants differ for $\wsf=0$ and in any case they do not give the right constants.} $\check{\mathrm{c}}$ and $\hat{\mathrm{c}}$ are chosen in the way that the resulting generating series transform like quasimodular forms of weight 2 for $\Gamma(2)$. They can be obtained from the conjectural formulas in Section \ref{formulas}. For example, in the case of $\QMs^{0}_{1}(q)'$, it is equal to $(2g-2)2^{4g-2}/3$.

\begin{conjno}[Enumerative mirror symmetry in $\tau$] Assume $\rsf=2$, let $q=e^{\pi i\tau}$, then
	\[	\frac{2}{ \tau^2}\QMs^{\asf}_{\dsf}(-1/\tau)'+\frac{\mathrm{const}}{\pi i \tau }=\sum_{\dsf'} \sum_{\mathsf a'} (-1)^{  \mathsf d \cdot \mathsf a'+\mathsf d' \cdot \mathsf a} \hat{\mathsf{QM}}^{\asf'}_{\dsf'}(\tau)'.\]
\end{conjno}
The statement of  Enumerative mirror symmetry in terms of $\Usf)i$ is  a consequence of modularity properties of the generating series and Enumerative mirror symmetry for the variable $\tau$, as it is explained in Section \ref{sduality}. The permutation of $\Usf_i$ is induced by the variable transformation $\tau \rightarrow -1/\tau$.




\subsection{Plan of the paper} The paper is roughly divided into three parts.  The first part, comprised of Sections \ref{intromoduli} - \ref{noncoprime}, is an adaption of \cite{N,NK3} to the set-up of Higgs sheaves and Higgs-Azumaya algebras (considered in \cite{Jiang}). The final product of these sections is the statement of Vafa--Witten/Quasimaps correspondence on the level of invariants, Theorem \ref{vwqm}, and the definition of invariants for the $\dsf=0$ case, Definition \ref{defnnoncoprime}.

The second part is comprised of Sections \ref{sectionMS} - \ref{calc}. In Section \ref{sectionMS}, we state the conjectures of Enumerative mirror symmetry for genus 1 quasimaps  invariants, Conjectures \ref{wnonzero} - \ref{wzero2}. In Section \ref{MM}, we show how to derive these conjectures from the calculations of \cite{MM} and Vafa--Witten/Quasimaps correspondence, if $\rsf=2$. In Section \ref{calc}, we do the basic calculations for $\wsf=0$, proving the conjectures in the case $g(X)=1$.  Even though $\Ms(\dsf)$ is a point for $\dsf\neq 0$,   the case of $g(X)=1$ is already non-trivial, because $\Ms(0)$ is far from being a point. 

The third part is Section \ref{sectionws}. This section is an adaption of the wall-crossing results of \cite{N,NK3} (whose one of the main ingredients is the master space of \cite{YZ}) and \cite{Ob}. In this section, we establish wall-crossing formulas between quasimap invariants and Gromov--Witten invariants, such that the wall-crossing invariants are given by Euler characteristics of Quot schemes on the curve $X$, Theorem \ref{wsE} and \ref{wsE2}.
\subsection{Relation to other work} 
A study similar to ours was undertaken in \cite{OT}. The approach of \cite{OT} is more physical in its nature and, as far as we understand, for genus 1 their statement is restricted to $\wsf=0$ case, i.e.\ Conjecture \ref{Enumerativemirror}. Moreover, we have doubts about their correspondence between Vafa--Witten invariants and curve-counting invariants for arbitrary genera, because such correspondence requires curve-counting invariants to be twisted, e.g.\ see \cite[Section 2.3]{Nqm}. 

The story described here resembles in many ways the one of \cite{OPa}, where genus 1 Gromov--Witten invariants of Hilbert scheme of points on a K3 surface were considered. In fact, \cite{N,NK3} and consequently this work were inspired by  \cite{OPa} and \cite{OP1}. 
\subsection{Further directions} \label{sectionfurther}
\subsubsection{Genus 0} By degenerating the elliptic curve $E$ to a nodal genus 1 curve,  one can show that the genus 1 invariants are equal to the trace of the operator given by 2-point genus 0 invariants, see \cite[(18)]{OPa}. We believe that arbitrary genus 0 invariants should also satisfy some form of Enumerative mirror symmetry, which therefore would vastly extend the genus 1 case. 

 For Vafa--Witten invariants, $S$-duality holds for arbitrary $\mu$-insertions, as is shown in \cite{MM}. Expressing genus 0 invariants in terms of $\mu$-insertions would be feasible approach to the problem, exchanging descendent insertions with relative insertions is a well-known practice \cite{AO}. In some sense, genus-0 Enumerative symmetry holds in some form almost tautologically, because one can view  $\mu$-insertions as quasimap insertions coming from cohomology of stacks $\FMsp(\dsf)$. 

 Moreover, using \cite{Ru}, one can define a quantum multiplication on the twisted orbifold cohomology. Hence one may speculate that \cite{HT}  can be extended to the level of cohomology rings or even quantum cohomology rings. 
\subsubsection{Geometric representation theory} In principle, it is reasonable to expect a story parallel to \cite{MO}. Especially in the light of the recent work \cite{HMMS}, where an action of an algebra on the cohomology of moduli spaces of Higgs bundles is constructed. Moreover, the form of expressions in Conjectures \ref{wnonzero} and \ref{wzero} suggests that the quantum multiplication by a divisor should indeed have an expression similar to \cite[Theorem 1.3.2]{MO}. 

 We also believe that Yun's Global Springer theory \cite{Yun} should play an important role in the story relating enumerative invariants with geometric representation theory for moduli spaces of Higgs bundles.

\subsubsection{Compact hyper-K\"ahler varieties} Recently, compact dual hyper-K\"ahler Lagrangian fibrations were constructed for all known compact hyper-K\"ahler deformation types in \cite{Kim}. The dual fibration is given by a finite-group quotient of a fibration. We think that our Enumerative mirror symmetry should hold in some form for this situation too, i.e.\ curve counting invariants of dual hyper-K\"ahler  fibrations should be related.

\subsubsection{Enumerative geometry of Artin stacks} The statement of Enumerative mirror symmetry in genus 1 unavoidably involves quasimap invariants  of Artin stacks $\Ms(0)$ and $\Mp(0)$, defined via the correspondence with sheaves. 

Some steps towards defining Gromov--Witten invariants of Artin stacks were recently made in \cite{DI} and in \cite{ArP}. Our construction of quasimap invariants for $\Ms(0)$ and $\Mp(0)$ suggests that it is indeed possible to define some version of Gromov--Witten invariants for Artin stacks. However, one needs to take extra care of stability. 

\subsection{Notation and convention}  We work over $\BC$. We denote the torus that scales Higgs fields by $\BC_t^*$, while the torus that scales $\p^1$ (with weight 1 at $0 \in \p^1$)  by $\BC_z^*$. Other $\BC^*$-actions will not have any subscripts and it will be clear from the context what they are. 

We also denote 

\[ \mathbf t := \text{weight 1 representation of } \BC_t^* \text{ on } \BC,\]
\[  \mathbf z := \text{weight 1 representation of } \BC_z^* \text{ on } \BC,\]
such that $-t:=e_{\BC_t}(\mathbf t )$ and  $z:=e_{\BC_z}(\mathbf z)$ are the associated classes in the equivariant cohomology of a point. Sign in the expression for $e_{\BC_t}(\mathbf t )$ is there to match the conventions of \cite[Appendix C]{GK2}, because their torus $\BC_t^*$ scales tangent directions with weight 1, while our torus $\BC_t^*$ scales cotangent directions with weight 1. 

We define the discriminant of a sheaf $\Delta(F):=\mathrm c_1^2(F)-2\rsf \mathrm{ch}_2(F)$. We will denote $\BZ_\rsf:= \BZ/\rsf\BZ$. We will also identify $\rsf$-roots of unity $\mu_\rsf$ of $\BC$ with $\BZ_\rsf$. 

\subsection{Acknowledgments} I thank Georg Oberdieck, Martijn Kool,  Richard Thomas, Henry Liu, Jan Manschot  for useful discussions on related topics. Special thanks goes to Martijn Kool for drawing my attention to \cite{MM}, which made explicit conjectural expressions for quasimap invariants possible.  

	\section{Moduli spaces of Higgs bundles} \label{intromoduli}
	\subsection{Higgs $\mathrm{SL} _\rsf-$bundles} \label{SLbundles}
		Let $X$ be a smooth proper connected curve over $\BC$ with $g:=g(X)\geq 1$. Throughout the article we fix a prime number $\rsf \geq 2$. 
		\begin{defn} A Higgs sheaf is pair $(G, \phi )$, a coherent sheaf $G$ on $X$ and $\phi \in \Hom(G,G\otimes \omega_X)$. 
			
			A Higgs sheaf is \textit{semistable}, if it is torsion-free and for all non-zero $\phi$-invariant subsheaves $G' \subset G$ we have $\mu(G')\leq \mu(G)$. It is \textit{stable}, if the inequality is strict. 
			\end{defn}
		
 By the spectral curve construction \cite{BNR}, the category of pairs $(G, \phi)$ is equivalent to the category of compactly supported 1-dimensional coherent sheaves $\CG$ on the total space of $\omega_X$, which we denote by $\rT^*X$. The equivalence is given by the pushforward  along  the natural map $\rT^*X \rightarrow X$.  We will refer to $\CG$  as a \textit{1-dimensional sheaf}. 
\\

	A family of Higgs sheaves on $X$ over a base scheme $B$ is pair $(G_B, \phi_B)$, where $G_B$ is a sheaf on $X\times B$ flat over $B$ and $\phi_B\in \Hom(G_B, G_B\otimes p^*_X\omega_X)$.

	\begin{defn} Gievn a class $\vsf:=(\rsf,\dsf) \in H^{\mathrm{even}}(X,\BZ) \cong \BZ\oplus \BZ$, let 
		\begin{align*}
			\FMs(\dsf)\colon (Sch/ \BC)^{\circ} &\rightarrow (Grpd) \\
			B &\mapsto \left\{(G_B,\phi_B) \ \Bigl\rvert \ \arraycolsep=0.1pt\def\arraystretch{1} \begin{array}{c}\det(G_b)\cong L,  \tr(\phi_b)=0, \\[.001cm] 	\ch(G_b)=\vsf\end{array} \right\}			
			\end{align*}
		be a moduli stack of Higgs sheaves on $X$.

\end{defn}
By tensoring sheaves with line bundles, the moduli stacks become isomorphic for different choices of $L$, hence we drop $L$ from the notation. Moreover, the moduli stacks $\FMs (\dsf)$ and $\FMs( \dsf')$ are isomorphic, if 
\[\dsf=\dsf' \ \mathrm{mod} \ \rsf.\]
To eliminate this redundancy, throughout the article we require 
\[0 \leq \dsf<\rsf.\]
We will thereby often identify $\dsf$ with its class in $\BZ_\rsf:=\BZ/\rsf\BZ$.
\\

\begin{defn} Let 
\[\check{\CM}(\dsf) \subset \FMs(\dsf),\]
be the open locus of semistable Higgs sheaves. Let
\[\Ms(\dsf):= \check{\CM}(\dsf)\thickslash \BC^{*}\]
be its rigidification by $\BC^*$-scaling action. 

\end{defn}
The space $\Ms(\dsf)$ in turn embeds into the rigidified stack of  Higgs sheaves, 
 \[\Ms(\dsf)\subset \FMsr(\dsf) :=\FMs(\dsf) \thickslash \BC^{*}.\]
 Overall, this can be represented by the following square, 

	\[
\xymatrix{
	\check{\CM}(\dsf)\ar@{^{(}->}[r] \ar[d]^{\BC^{*}-\text{gerbe}} &
 \FMs(\dsf) \ar[d]^{\BC^{*}-\text{gerbe}}
	\\
	\Ms(\dsf)\ar@{^{(}->}[r] & \FMsr(\dsf)
}
\]
\noindent Let
\[ \check{\CN}(\dsf) \subset  \check{\CM}(\dsf)\]
be the moduli stack of semistable sheaves on $X$, embedded as the locus of Higgs sheaves with a trivial Higgs field. By $\check{N}(\dsf)$ we will denote the corresponding rigidification by $\BC^*$.

\subsection{Determinant line bundles.} \label{sectiondet}

Let $(\Gsf, \phi)$ be the universal Higgs sheaf on 
\[X \times \FMs(\dsf). \]
We define the determinant line bundle map as the composition
\begin{multline*}\lambda \colon K_{0}(X)\xrightarrow{p_{S}^{!}} K_{0}(X\times \FMs(\dsf)) \xrightarrow{\cdot [\Gsf]}  K_{0}(X\times \FMs(\dsf))\\ \xrightarrow{p_{\FMs(\dsf)!}} K_{0}(\FMs(\dsf))\xrightarrow{\det} \Pic(\FMs(\dsf)).
	\end{multline*}
 \begin{rmk} The definition of $\lambda$ requires some care, because  $\FMs(\dsf)$ is not of finite type, hence the determinant map cannot be defined naively. The correct approach is provided by Waldhausen K-theory \cite{Wa}. See also \cite[Section 3]{STV}. 
 	\end{rmk}
 	

In general, the $\BC^{*}$-weight of the line bundle $\lambda(u)$ is equal to the Euler characteristics,
\[w_{\BC^{*}}(\lambda(u))=\chi(\vsf \cdot u):=\int_X \vsf \cdot u \cdot \td_X.\]
There are two types of classes that we will be of interest to us. A class $u \in K_{0}(X)$, such that $\chi(\vsf \cdot u)=1$, gives a trivilisation of the  $\BC^{*}$-gerbe 
\[ \FMs(\dsf)  \rightarrow \FMsr(\dsf) \]
 or, in other words, a universal family on $\FMsr(\dsf)$, which is given as a pullback of $(\Gsf, \Phi)$ by the section defined by $\lambda(u)$. 
\\

While for a class $u \in K_{0}(X)$, such that $\chi(\vsf \cdot u )=0$, the line bundle $\lambda(u)$ descends to $\FMsr(\dsf)$. Let 
\[K_0(X)^{\perp}=\{u \in K_{0}(X) \mid  \chi(\vsf \cdot u)=0 \}  \subset K_{0}(X),\]
then  $\lambda$ restricted to $K_0(X)^{\perp}$ descends to a map to $\Pic(\FMsr(\dsf))$, 
\[\lambda \colon K_0(X)^{\perp} \rightarrow \Pic(\FMsr(\dsf)).\]
The map $\lambda$ factors through the numerical $K$-group, 
\[\lambda \colon K_0(X)_{\mathrm{num}}^{\perp} \rightarrow \Pic(\FMsr(\dsf)).\]
By associating rank and degree to a class  in  $K_0(X)_{\mathrm{num}}$, we obtain an isomorphism 
\[K_0(X)_{\mathrm{num}} \cong \BZ \oplus \BZ, \quad u \mapsto (\rk(u), \deg(u)).\]
By Riemann--Roch theorem,  a class $(r',d') \in K_0(X)_{\mathrm{num}}$ is perpendicular to $(\rsf,\dsf)$ with respect to $\chi$, if and only if 
\begin{equation} \label{perp}
(r',d')=m \cdot (\rsf, \rsf(g-1)-\dsf), \quad m \in \BZ. 
\end{equation}
Therefore, 
\[K_0(X)_{\mathrm{num}}^{\perp}\cong \BZ.\]

\begin{thm}	 \label{thet}
	The restriction of $\lambda$ to $\check{M}(\dsf)$ is an isomorphism, 
	\[\lambda_{\check{M}(\dsf)} \colon \BZ \xrightarrow{\sim} \Pic(\Ms(\dsf)). \]
	\end{thm}
\textit{Proof.} By composing $\lambda_{\check{M}(\dsf)}$ with the restriction $\Pic(\check{M}(\dsf)) \rightarrow \Pic(N(\dsf))$, we obtain 
\[\lambda_{N(\dsf)} \colon K_0(X)_{\mathrm{num}}^{\perp} \rightarrow \Pic(\Ns(\dsf)),\]
which is just the determinant-line-bundle map associated to the moduli space $\Ns(\dsf)$. 
By \cite{DN} (see also \cite{KN} and \cite{BLS}), the map $\lambda_{\Ns(\dsf)}$ is an isomorphism.  By \cite{SP}, restriction of line bundles from $\check{M}(\dsf)$ to $\Ns(\dsf)$ is an isomorphism.
\qed 
\\

The ample generator of $\Pic(\Ns(\dsf))$ is  the \textit{theta line bundle} $\Theta_{\Ns(\dsf)}$, 
\[\Pic(\Ns(\dsf))=\BZ  \Theta_{\Ns(\dsf)}.\]
By Theorem \ref{thet}, we obtain that 
\[\Pic(\Ms(\dsf))= \BZ \check{\Theta},\]
where $\Thetas$ is the line bundle which restricts to $\Theta_{\Ns(\dsf)}$ on $\Ns(\dsf)$. 
Let us denote the corresponding class in  $K_0(X)_{\mathrm{num}}$ by $\theta$,
\[\theta=(-\rsf,\dsf-\rsf(g-1)) \in K_0(X)_{\mathrm{num}}, \quad \lambda(\theta)=\Thetas,\]
it is associated to $m=-1$ in $(\ref{perp})$. We will also denote corresponding line bundle on $\FMsr(\dsf)$ by $\Thetas$. 
\\

If $\dsf\neq 0$, there exists an element  $\usf \in K_{0}(X)$, such that 
\[\chi(\vsf \cdot \usf)=1,\]
the line bundle $\lambda(\usf)$ is of $\BC^*$-weight 1, hence it trivialises the $\BC^*$-gerbe 
\[ \FMs(\dsf)  \rightarrow \FMsr(\dsf).\]
We will denote the trivialising section by 
\[s_\usf \colon \FMsr(\dsf) \rightarrow  \FMs(\dsf).\]
\vspace{0.1cm}

\noindent \textbf{Important} We fix $\usf \in  K_0(X)_{\mathrm{num}}$, such that $\chi(\vsf \cdot \usf)=1$, throughout the article for all $\dsf\neq 0$. This fixes the universal family of Higgs sheaves on $\FMsr(\dsf)$, which we also denote by $(\Gsf, \phi)$. 

\vspace{0.3cm}
\subsection{Higgs $\mathrm{PGL}_r$-bundles } \label{PGLbundles}
We will now introduce the notion of generalised Azumaya algebras defined by Lieblich \cite{Lie}. Lieblich distinguishes between \textit{pre-generalised Azumaya algebras} and \textit{generalised Azumaya algebras}, such that the latter is given by stackification of the moduli functor associated to the former. By \cite[Section 6.4]{Lie}, it is not necessary to distinguish between the two in the case of surfaces and curves. Hence for us,   \textit{generalised Azumaya algebras} will mean Lieblich's  \textit{pre-generalised Azumaya algebras}.
\begin{defn} A \textit{generalised Azumaya algebra} of rank $\rsf$ on a scheme $Y$  is a perfect object $\CA\in D(Y)$ endowed with the structure of a weak algebra, such that there exists an \'etale cover $\{U_i\}$ of $Y$ and an isomorphism of weak algebras
	\[\CA_{|U_i} \cong R\CE nd_U(G_i,G_i) \]
 for a rank $\rsf$ sheaf\footnote{The original definition of Lieblich, \cite[Definition 5.2.1]{Lie} assumed that $G_i$ is torsion-free, we, however, consider arbitrary sheaves. If we do it on a curve,  all the results of Lieblich apply to this set-up verbatim.} $G_i$ on $U_i$.  We say that $\CA$ is \textit{torsion-free}, if $G_i$ are torsion free for all $i$.
	
	A \textit{generalised Higgs--Azumaya algebra} is a pair $(\CA,\phi)$, where $\CA$ is a generalised Azumaya algebras and $\phi \in H^0(\CA \otimes \omega_Y,Y)$. Slightly abusing the terminology,  we will refer to generalised Azumaya algebra as \textit{azalgebra} and to generalised Higgs--Azumaya algebra as \textit{Higgs azalgebra}.
\end{defn}

Given an azalgebra $\CA$ on a scheme $Y$, we can associate a class 
\[\delta(\CA) \in H^2_{\text{\'et}}(Y,\BZ_\rsf),\]
constructed as follows. We firstly define a stack $\CY_{\CA}$, whose $B$-valued points are trivialisations of $\CA$ with a trivial determinant,
\begin{equation} \label{trivilisations}
\CY_{\CA}(B)=\left\{(G_B,\psi_B, \gamma_B)  \ \Bigl\rvert \ \arraycolsep=0.1pt\def\arraystretch{1} \begin{array}{c} \psi_B \colon \CA_{|B} \xrightarrow{\sim} R \CE nd_B(G_B,G_B), \\[.001cm] 	\gamma_B \colon \det(G_B)\xrightarrow{\sim} \CO_B \end{array} \right\} 
\end{equation}
where an isomorphism between triples $(G_B,\psi_B, \gamma_B)$ and  $(G'_B,\psi_B, \gamma_B)$ is an isomorphism between $G_B$ and $G_B'$ preserving the maps. The stack $\CY_\CA$ is an \'etale $\BZ_\rsf$-gerbe and therefore has an associated class $[ \CY_{\CA}]\in H^2_{\text{\'et}}(Y,\BZ_\rsf)$. We define 
\[ \delta(\CA):= [ \CY_{\CA}] \in H^2_{\text{\'et}}(Y,\BZ_\rsf),\]
we will refer to $\delta(\CA)$ as the \textit{degree} of $\CA$.
\begin{defn} We say that an azalgebra $\CA$ on $X$ is in the class 
	\[\vsf=(\rsf,\dsf) \in 
H^0(X,\BZ) \oplus H^2(X,\BZ_\rsf)\cong \BZ\oplus \BZ_\rsf,\]
 if it is of rank $\rsf$ and degree $\dsf$. 
\end{defn}


A family of Higgs azalgberas on $X$ over a base scheme $B$ is pair $(\CA_B, \phi_B)$, where $\CA_B$ is an azalgebra\footnote{Since $X$ is a smooth curve, a family of perfect objects on $X\times B$ is perfect and vice versa. } on $X\times B$ and $\phi_B\in H^0(\CA_B\otimes p^*_X\omega_X)$. 
\begin{defn}
Let
\begin{align*} 
	\HAz(\dsf) \colon (Sch/ \BC)^{\circ} &\rightarrow (Grpd)\\ 
	B &\mapsto \left\{(\CA_B,\phi_B) \ \Bigl\rvert \ \arraycolsep=0.1pt\def\arraystretch{1} \begin{array}{c}\tr(\phi_b)=0, \\[.001cm] 	(\rk(\CA_b),\delta(\CA_b))=\vsf \end{array} \right\}			
\end{align*}
be the stack of Higgs azalgebras on $X$.
\end{defn}
The space $\HAz(\dsf)$ is indeed an Artin stack by the arguments of \cite[Section 6.4]{Lie} (also follows from \cite{To} and \cite{Lie2}). \\

On a curve $X$, all azalgebras are globally of the form $R \CE nd_X(G,G)$ for some sheaf $G$ on $X$ which is in the class $\vsf=(\rsf,\dsf)$. This can be seen as  follows. Let $Y$ be an arbitrary scheme. Consider the natural map 
\begin{equation} \label{map}
H_{\text{\'et}}^2(Y,\BZ_\rsf) \rightarrow H_{\text{\'et}}^2(Y,\BC^*)
\end{equation}
induced by the inclusion of the sheaf $\BZ_\rsf \hookrightarrow \BC^*$ as $\rsf$-roots of unity.
\begin{defn} \label{brauergerbe}
	Let $\mathrm{Br}(\CA)  \in H_{\text{\'et}}^2(Y,\BC^*)$ be the image of $\delta(\CA)$ under map (\ref{map}), called the \textit{Brauer class} of $\CA$. Let $\CY$ be a $\BC^*$-gerbe that represents the class $\mathrm{Br}(\CA)$.
	\end{defn}
 In other words, to obtain  $\CY$, we forget the data of $\gamma$ in the definition of $\CY_{\CA}$.  By construction of $\CY_{\CA}$ and $\CY$, there exists a (twisted) sheaf $G_{\CY}$ on $\CY$ , such that 
\[R\CE nd_\CY(G_{\CY},G_{\CY})\cong \CA.\]
Since $X$ is a smooth curve, 
\[ H^2_{\text{\'et}}(X,\BC^*)=0,\]
 the associated gerbe $\CX$ is therefore a trivial $\BC^*$-gerbe and  $G_{\CX}$ descends to a sheaf $G$ on $X$, such that  $R\CE nd_X(G,G)\cong \CA$, as claimed.

We thereby  can define a map between stacks $\FMsr(\dsf)$ and $\HAz(\dsf)$,
\begin{equation} \label{etale2}
	\FMsr(\dsf) \rightarrow \HAz(\dsf), \quad  (G,\phi)\mapsto (R\CE nd_X(G,G),\phi).   
\end{equation}
The results of Lieblich imply that this map is a torsor. 

\begin{lemma}  \label{torsor}
	The map (\ref{etale2}) is a $\Gamma_X$-torsor under the natural action of $\Gamma_X$ on $\FMsr(\dsf)$  given by tensoring a sheaf with a line bundle.  
\end{lemma}
\textit{Proof.} The claim  follows from the (derived) Skolem--Noether theorem of Lieblich, \cite[Corollary 5.1.6]{Lie}. \qed 

\begin{defn}
Let \[\Mp(\dsf) \subset \HAz(\dsf)\] 
be the open locus consisting of Higgs azalgebras $\CE nd_X(G,G)$ on $X$ such that $G$ is a semistable bundle. We refer to elements of $\Mp(\dsf)$ as (semistable) Higgs $\mathrm{PGL}_r$-bundles. 
\end{defn}


By the results of \cite{BLS}, the theta line bundle\footnote{In \cite{BLS}, a gerbe of $\FMsr(\dsf)$, defined in \ref{gerbeclasses} is considered,  the power of their determinant line bundle $\CD^\rsf$ descends to $\Thetas$ on $\FMsr(\dsf)$. The same power descends from the gerbe to $\FMp(\dsf)$. } $\Thetas$ descends from the locus of Higgs locally free sheaves ($\mathrm{SL} _\rsf$-bundles)
\[\FMsr^{\mathrm{lf}}(\dsf) \subset \FMsr(\dsf)\] to the locus of  \textit{classical} Higgs-Azumaya algebras ($\mathrm{PGL}_r$-bundles)
\[\FMp^{\mathrm{lf}}(\dsf)  \subset \HAz(\dsf).\]
Extending this result to the entire $\HAz(\dsf)$ seems a bit too involved for our purposes, as we  need $\Thetas$ mainly for the definition of the degree of a quasimap. However, if we view
\[[\Thetas] \in H^2(\FMsr(\dsf), \BQ )\]
as a class in some cohomology theory $H^2(\FMsr(\dsf), \BQ )$, then $[\Thetas]$ descends to a class 
\[[\Thetap] \in H^2(\HAz(\dsf), \BQ ),\]
because it is $\Gamma_X$-invariant. This can be made more precise by  the construction in Definition \ref{torsorcurve}. 

\begin{rmk}Let us say a few words about $H^2(\FMsr(\dsf), \BQ )$. Perhaps, the best way to define it is via limit of singular cohomologies of quasicompact open substacks. In the end, we use $[\Thetap]$ only to define the degree of quasimaps and any of the possible ways will work. 
\end{rmk}
\subsection{Gerbe classes} \label{gerbeclasses}
 Let $\FMs_{\mathrm{SL} }(\dsf)$ be the moduli stack, whose $B$-valued points are defined as follows
\begin{equation*}
	\FMs_{\mathrm{SL} }(\dsf)(B)=\left\{(G_B, \phi_B, \gamma_B) \ \Bigl\rvert \ \arraycolsep=0.1pt\def\arraystretch{1} \begin{array}{c}(G_B, \phi_B,) \text{ is a family over $B$,} \\
	[.001cm] \gamma_B\colon \det(G_B) \xrightarrow{\sim} \pi^*_XL \end{array}  \right\}.
\end{equation*}
In other words, we add the data of identification of the determinant $\det(G_B)$ with $\pi^*_XL$, which has an effect of killing almost all $\BC^*$-scaling automorphisms of sheaves, only $\rsf$-roots of unity remain. 

Equivalently, $\FMs_{\mathrm{SL} }(\dsf)$ can be defined as a fiber product
\[
\begin{tikzcd}[row sep=scriptsize, column sep = scriptsize]
	&  \FMs_{\mathrm{SL} }(\dsf) \arrow[d]  \arrow[r,hook] &   \FM(\dsf) \arrow[d,"(\mathrm{det}\text{,}\tr)"]\\
	& \Spec(\BC)  \arrow[r, "\textit{(L,}0)", hook] &  \mathrm{T}^*\CP ic(X)
\end{tikzcd}
\]
where $\FM(\dsf)$ is the stack of Higgs sheaves without a fixed determinant, $\rT^*\CP ic(X)$ is the total space of the cotangent bundle over the Picard stack (rank 1 Higgs sheaves) and $[(L,0)]$ is a point corresponding to a line bundle $L$ with a zero cotangent vector. 

The stack $\FMs_{\mathrm{SL} }(\dsf)$ is a $\BZ_\rsf$-gerbe over $\FMsr(\dsf)$. We then define 
\[\check \alpha:= [\FMs_{\mathrm{SL} }(\dsf)] \in H_{\text{lis-\'et}}^2(\FMsr(\dsf),\BZ_\rsf)\]
to be its class. Here $H_{\text{lis-\'et}}^2(\FMsr(\dsf),\BZ_r)$ is defined to be the group of isomorphism classes of lisse-\'etale $\BZ_\rsf$-gerbes on $\FMsr(\dsf)$. 

\begin{rmk}The class $\check \alpha$ can be given a simpler interpretation. Let $\Gsf$ be the universal sheaf on $X\times \FMsr(\dsf)$.  Then $\check \alpha$ should be the pullback of the class $[\mathrm{c}_1(\Gsf)]\in H^2(X\times \FMsr(\dsf) ,\BZ_\rsf) $ by an embedding   $ \iota_x\colon x\times \FMsr(\dsf) \hookrightarrow X\times \FMsr(\dsf)$ for some point $x \in X$.  However, since $\FMsr(\dsf)$ is a non-quasicompact Artin stack, we choose to refrain from making any precise statements due to possible technical obstacles. In particular, one has to define the cohomology group $H^2(X\times \FMsr(\dsf) ,\BZ_\rsf)$ carefully and then prove a correspondence between $\BZ_\rsf$-gerbes and classes in $H^2(X\times \FMsr(\dsf) ,\BZ_\rsf)$. 
	\end{rmk}

One can verify that $\FMs_{\mathrm{SL} }(\dsf)$ descends to a $\BZ_\rsf$-gerbe over $\HAz(\dsf)$, whose class we denote by 
\[\hat{\alpha} \in H_{\text{lis-\'et}}^2(\HAz(\dsf),\BZ_\rsf).\]
In fact, it can be given as a pullback of the degree  $\delta(\Asf)$ by a fiber embedding as above, where $\Asf$ is the universal azalgebra over $\HAz(\dsf)$ (notice that the degree of azalgebras was defined in exactly the same way as the gerbe $\FMs_{\mathrm{SL} }(\dsf)$).

\begin{rmk}Over the stable locus $\Ms(\dsf)$, one can view $\check{\alpha}$ as a class in the ordinary singular cohomology with $\BZ_\rsf$-coefficients, using the natural comparison isomorphisms 
\[ H_{\text{lis-\'et}}^2(\Ms(\dsf),\BZ_\rsf) \cong H_{\text{\'et}}^2(\Ms(\dsf),\BZ_\rsf)\cong H_{\text{sing}}^2(\Ms(\dsf),\BZ_\rsf)\]
and the fact that $H_{\text{\'et}}^2(\Ms(\dsf),\BZ_\rsf)$ classifies $\BZ_\rsf$-gerbes. Over $\Mp(\dsf)$, the situation is more complicated, because $\Mp(\dsf)$ is an orbifold. The author was not able to find good references for the comparison and the classification results in this set-up.  
\end{rmk}


\section{$\mathrm{SL} _\rsf$-quasimaps, $\dsf\neq 0$} \label{slpre}

Through this section we assume $\dsf\neq 0$.  Let $(M, \FM)$ be a pair $(\Ms(\dsf), \FMsr(\dsf))$ or a pair $(\check{\CM}(\dsf),  \FMs(\dsf))$.
\begin{defn} \label{quasimapsdef} 
	A map $f\colon (C,\mathbf x) \rightarrow \FM $ is a \textit{quasimap} to $(M, \FM)$ of genus $g$ and of degree $\wsf \in \BZ$,  if 
	\begin{itemize}
		\item $(C,\mathbf x)$ is a  connected marked nodal curve of genus $g$,
		\item $\deg(f):=\deg(f^*\Thetas)=\wsf$;
		\item $| \{ p\in C \mid f(p) \in \FM \setminus M \} | < \infty$.
	\end{itemize}
	We will refer to the set $\{ p\in C\mid f(p) \in \FM \setminus M \}$ as \textit{base} points. A quasimap $f$ is \textit{prestable}, if 
	\begin{itemize}
		\item $\{ t\in C \mid f(t) \in \FM\setminus M \} \cap \{\mathrm{nodes}, \mathbf x \}=\emptyset$. 
	\end{itemize} 
\end{defn}
For a smooth curve $C$, we define 
\begin{align*}
\Lambda&:=H^0(X,\BZ)\oplus H^2(X,\BZ)\cong \BZ \oplus \BZ\\
\Lambda_{C}&:=H^1(X,\BZ) \otimes H^1(C,\BZ). 
\end{align*}

By K\"unneth's decomposition theorem, we have 
\begin{equation} \label{decomposition}
	\begin{split}
H^{\mathrm{even}}(X\times C,\BZ)
&\cong \Lambda \oplus \Lambda_C \oplus \Lambda . 
\end{split}
\end{equation}

\noindent Let $f \colon C \rightarrow \FMs(\dsf)$ be a quasimap of degree $\wsf$. By the definition of $\FMs(\dsf)$, the quasimap $f$ is given by a family of Higgs sheaves $(F,\phi)$ on $X \times C$, i.e.\ $F$ is a sheaf flat over $C$ and $\phi \in \Hom(F,F\otimes p^*_X\omega_X)$. The Chern character of $F$ has three components with respect to the decomposition in (\ref{decomposition}),
\[\ch(F)=(\ch(F)_{\mathrm{f}}, \ch(F)_{\mathrm{m}}, \ch(F)_{\mathrm{d}}) \in \Lambda \oplus \Lambda_C \oplus \Lambda.\]
Note that $\ch(F)$ is an integral class, because $H^2(X\times C,\BZ)$ is an even lattice. 
\\

Let us now analyse the components of $\ch(F)$. Firstly, since $F$ is a family of sheaves with a fixed determinant line bundle, $\det(F)$ is of the form $L\boxtimes L' $, hence the middle component is zero,
\[ \ch(F)_{\mathrm{m}}=0.\] 
  The subscripts "f" and "d" stand for $\textit{fiber}$ and $\textit{degree}$, respectively.  So as the notation suggests, \[\ch(F)_{\mathrm{f}}=\vsf,\]
which can be seen by pulling back $\ch(F)$ to a fiber over $C$ and using the flatness of $F$.  On the other hand, the degree component of $\ch(F)$ and $\deg(f)$ are related in the following way. 

\begin{lemma} \label{degreechern}
	\begin{align*}\deg(f)&= \int_{X} \ch(\theta) \cdot \ch(F)_{\mathrm{d}} \cdot \td_X \\
		&=\dsf\ch_1(F)_{\mathrm d}-\rsf \ch_2(F)_{\mathrm d} \\
		&=\mathrm c_1^2(F)/2-\rsf \mathrm{ch}_2(F) \\
		&=\Delta(F)/2.
		\end{align*}
	\end{lemma}

\textit{Proof.} By the functoriality of the determinant line bundle construction,
\[\deg(f)=\deg(f^*\Thetas)=\deg(\lambda_{F}(\theta)),\]
where $\lambda_{F}(\theta)$ is the determinant line bundle associated to the family $F$ and a class $\theta \in K_{0}(X)$. Using the Grothendieck--Riemann--Roch theorem and the projection formula, we obtain
\begin{align*}
	\deg(\lambda_F(\theta))&=\int_{C} \ch(p_{C!}(p^{!}_{S}\theta\cdot [F]))\\
	&=\int_{X\times C}\ch(p^!_{S}\theta\cdot [F]) \cdot p^*_{X}\td_X\\
	&=\int_X \ch(\theta) \cdot p_{S*}\ch(F) \cdot \td_X\\
	&=\int_{X} \ch(\theta) \cdot \ch(F)_{\mathrm{d}} \cdot \td_X.
\end{align*}
By (\ref{perp}), the last quantity can be made more explicit,
\begin{equation} 
	\begin{split}
		\int_{X} \ch(\theta) \cdot \ch(F)_{\mathrm{d}} \cdot \td_X =\dsf\ch_1(F)_{\mathrm d}-\rsf \ch_2(F)_{\mathrm d}. 
	\end{split}
\end{equation}
Finally,
\begin{align*}
	\dsf\ch_1(F)_{\mathrm d}-\rsf \ch_2(F)_{\mathrm d}&=\mathrm c_1^2(F)/2-\rsf \mathrm{ch}_2(F)\\
	&=\Delta(F)/2.
\end{align*}
\qed 
\\

Composing a quasimap $f \colon C \rightarrow \FMsr(\dsf)$  with a section $s_\usf \colon \FMsr(\dsf) \rightarrow \FMs(\dsf)$, we obtain a quasimap 
\[f_u:=s_\usf \circ f \colon C \rightarrow \FMs(\dsf),\]
which has an associated  family $(F,\phi)$.  By construction, the section $s_{\usf}$ is given by the descend of the sheaf $\Gsf\otimes\lambda(\usf)^{-1}$ from  $\FMs(\dsf)$ to $\FMsr(\dsf)$.  A family $F$ is a pullback of $\Gsf\otimes\lambda(\usf)^{-1}$ by a quasimap $f_\usf\colon C  \rightarrow \FMs(\dsf)$. By the definition of $\lambda(\usf)$, the family $F$ therefore satisfies 
\begin{equation} \label{detcond}
\det(p_{C*}(p_{S}^{*}\usf\otimes F))\cong \CO_{C},
\end{equation}
which, in particular, implies that 
\begin{equation} \label{u}
\int_{X} \ch(\usf) \cdot \ch(F)_{\mathrm{d}} \cdot \td_X=0.
\end{equation}

Combining $(\ref{u})$ with Lemma \ref{degreechern}, we obtain the following lemma. 

\begin{lemma} \label{lemmaeqchern}If $f \colon C \rightarrow \FMsr(\dsf)$ is of degree $\wsf$, the degree component of the associated family $F$, 
\[\ch(F_\usf)_{\mathrm d}=(\ch_1(F)_{\mathrm d},\ch_2(F)_{\mathrm d}) \in \Lambda \cong \BZ \oplus \BZ,\]
is the unique solution to the following system of equations
	\begin{equation} \label{eqchern}
		\begin{split}
&\int_{X} \ch(\usf) \cdot \ch(F)_{\mathrm{d}} \cdot \td_X=0, \\
& \int_{X} \ch(\theta) \cdot \ch(F)_{\mathrm{d}} \cdot \td_X=\wsf.
\end{split}
\end{equation}
We will denote the solution by $\check \wsf$, 
\[\ch(F)_{\mathrm d}=\check \wsf.\]
\end{lemma}

 \subsection{Moduli spaces of quasimaps}

Consider a quasimap $f\colon C \rightarrow  \FMsr(\dsf)$ and a point $p \in C$ in the regular locus of $C$. By Langton's semistable reduction \cite{Lang}, we can modify the quasimap $f$ at the point $p$, to obtain a quasimap  
\[f_p \colon C \rightarrow  \FMsr(\dsf)\]
which maps to the stable locus $\Ms(\dsf)\subset \FMsr(\dsf)$ at $p$ (if $p$ is not a base point, then $f_p=f$).   We refer to $f_p$ as \textit{stabilisation} of $f$ at $p$.

\begin{defn}
	We define the \textit{length} of a point $p \in C$ to be 
	\[\ell(p):=\deg(f)-\deg(f_p).\]
	By the proof of \cite[Proposition 3.7]{N}, $\ell(p)\geq 0$ and $\ell(p)=0$, if and only if $p$ is not a base point. 
\end{defn}
In what follows, by $0^+$ we will denote a number $A \in \BR_{>0}$, such that $A \ll 1$. 
\begin{defn} \label{stabilityqm} Given $\epsilon \in \BR_{>0}\cup\{0^+, \infty\}$, a prestable quasimap $f\colon(C, \mathbf x) \rightarrow  \FMsr(\dsf)$ is $\epsilon$-stable, if 
	\begin{itemize}
		\item[$\mathbf{(i)}$] $\omega_C(\mathbf x) \otimes f^*\Thetas^\epsilon$ is positive;
		\item[$\mathbf{(ii)}$] $\epsilon \ell(p)\leq 1$ for all $p \in C$.
	\end{itemize}
	We will refer to $0^+$-stable and $\infty$-stable quasimaps as stable quasimaps and stable maps respectively. 

\end{defn} A \textit{family} of quasimap over a base $B$ is a family of nodal curves $\CC$ over $B$ with a map $f\colon \CC \rightarrow \FMsr(\dsf)$ such that the geometric fibers of $f$ over $B$ are quasimaps.  
\begin{defn}
Let
\begin{align*}
	Q^{\epsilon}_{g,n}(\check{M}(\dsf), \wsf) \colon (Sch/ \BC)^{\circ} &\rightarrow (Grpd) \\
	B&\mapsto \{\text{families of $\epsilon$-stable quasimaps over }B\}
\end{align*}
be the moduli space of $\epsilon$-stable quasimaps of genus $g$ and the degree $\wsf$ with $n$ marked points. 
\end{defn}
The $\BC_t^*$-action on $\Ms(\dsf)$ given by the scaling of Higgs fields induces a $\BC_t^*$-action on $Q^{\epsilon}_{g,n}(\Ms(\dsf), \wsf)$. Since  only the fixed locus $\Ms(\dsf)^{\BC_t^*}$ is proper, this action will be crucial for defining enumerative invariants.

\begin{lemma} \label{extension0}
	Let $\CC \rightarrow \Delta$ be a family of nodal curves over a DVR $\Delta$ and $\{p_i\}\subset \CC$ be finitely many closed points in the regular locus of the central fiber.  Any quasimap 
	$f^{\circ} \colon \CC^{\circ}=\CC \setminus \{p_i\} \rightarrow \FMsr(\dsf)^{\BC_t^*}$ extends to a quasimap $f\colon\CC \rightarrow \FMsr(\dsf)^{\BC_t^*}$, which is unique up to unique isomorphism. 
\end{lemma}
\textit{Proof. } Let $ \CF^\circ$ be the sheaf on $\rT^*X\times  \CC^\circ$ associated to a quisimap $f^\circ$.  We extend $\CF^\circ$ to a coherent sheaf $\CF$ on $\rT^*X
 \times \CC$, quotienting subsheaves of lower dimension, if necessary. The sheaf $\CF$ is therefore flat over $\Delta$ and $\BC_t^*$-invariant, because the generic fiber over $\Delta$ is $\BC_t^*$-invariant. The central fiber $\CF_{\Spec \BC}$ of $\CF$ defines a quasimap, if it is pure, because $\CC_{\Spec \BC}$ is regular at $p_i$. If $\CF_{\Spec \BC}$ is non-pure, we can remove the subsheaves of lower dimension inductively as follows. Let $\CF^{0}=\CF$ and  $\CF^{i}$ be defined by short exact sequences, 
\[0 \rightarrow \CF^{i} \rightarrow \CF^{i-1} \rightarrow \CQ^{i} \rightarrow 0,\]
such that $\CQ^{i}$ is the quotient of $\CF^{i-1}_{\Spec \BC}$ by the maximal subsheaf of lower dimension. It is not difficult to check, that at each step the torsion of $\CF_{\Spec \BC}^{i}$ is supported at slices $\rT^*X \times p_{i}$, therefore all $\CF^i$'s are isomorphic to $\CF^0$ over $\rT^*X \times \tilde{\CC}$. Moreover, $\BC_t^*$-invariance is preserved at each step, because the sheaf remains unchanged at the generic fiber.  By the standard argument (see e.g.\ \cite[Theorem 2.B.1]{HL}), this process terminates, i.e.\ $\CF^i=\CF^{i+1}$ and $\CF^i_{k}$ is pure for $i\gg 0$. Let us  redefine the sheaf $\CF$, we set $\CF=\CF^i$ for some $i\gg0$, then the sheaf $\CF$ induces a quasimap to $\FMs(\dsf)^{\BC_t^*}$. By composing it with the projection to $\FMsr(\dsf)^{\BC_t^*}$, we obtain an extension $f\colon \CC \rightarrow \FMsr(\dsf)^{\BC_t^*}$ of $\tilde{f}$. 

Consider now another extension $f'\colon \CC \rightarrow \FMsr(\dsf)^{\BC_t^*}$, we lift both $f$ and $f'$ to $\FMs(\dsf)^{\BC_t^*}$. Let $\CF'$ and $\CF$ be the corresponding families on $S\times \CC$ (notice, $\CF$ might differ from the previous $\CF$ by a tensor with a line bundle), by \cite[Lemma 3.17]{N} they define a family of stable sheaves relative to $\Delta$ in some relative moduli of sheaves $\CM(\rT^*X\times \CC/\Delta)$, hence they must be isomorphic up to tensoring with a line bundle by separateness of the relative moduli space of stable sheaves. The isomorphism becomes unique after projection to $\FMsr(\dsf)^{\BC_t^*}$.
\qed

\begin{prop} \label{positive}
	A prestable quasimap $f\colon C \rightarrow \FMsr(\dsf)$ is non-constant, if and only if
	\[ \deg(f)>0.\]
\end{prop}  
\textit{Proof.}
As in \cite[Lemma 3.5]{N}, we have that
\[\deg(f^*\Thetas)=\dsf\cdot \mathrm{rk}(p_{S*}F)-\rsf \cdot \mathrm{deg}(p_{S*}F).\]
The claim then follows from the same arguments as in \cite[Proposition 3.7]{N} and ampleness of $\Thetas$ on $\Ms (\dsf)$. In fact, the proof is even simpler, because the bounds  for $\ch_2$ are not required (cf. \cite[Lemma 3.8]{N}). 
\qed

\begin{thm} \label{properness}
	A moduli space $Q^{\epsilon}_{g,n}(\Ms(\dsf), \wsf)$ is a Deligne--Mumford stack of finite type. Moreover, the $\BC_t^*$-fixed locus $Q^{\epsilon}_{g,n}(\Ms(\dsf), \wsf)^{\BC_t^*}$ is proper. 
	\end{thm}
\textit{Proof.} Using Proposition \ref{positive}, the proof of the first claim follows the same arguments as in \cite[Theorem 3.21]{N}.   Let us now show that the fixed locus $Q^{\epsilon}_{g,n}(\Ms(\dsf), \wsf)^{\BC_t^*}$ is proper. Firstly, there is a natural identification, 
\[Q^{\epsilon}_{g,n}(\Ms(\dsf), \wsf)^{\BC_t^*}=Q^{\epsilon}_{g,n}(\Ms(\dsf)^{\BC_t^*}, \wsf),\]
where $\Ms(\dsf)^{\BC_t^*}$  is the fixed locus of $\Ms(\dsf)$. On the right, we have the moduli space of quasimaps  whose target is the pair  $\Ms(\dsf)^{\BC_t^*} \subset \FMs(\dsf)^{\BC_t^*}$.  We have to show that $Q^{\epsilon}_{g,n}(\Ms(\dsf)^{\BC_t^*}, \wsf)$ is proper.  In the light of properness of the fixed locus $\Ms(\dsf)^{\BC_t^*}$ and Lemma \ref{extension0}, properness of $Q^{\epsilon}_{g,n}(\Ms(\dsf)^{\BC_t^*}, \wsf)$ follows from the same arguments as in \cite[Section 4.3]{CFKM}. 

\subsection{Obstruction theory}  \label{obstructionsection} The moduli stack $\FMsr(\dsf)$ can be very singular, hence one cannot naively use its tangent complex $\BT_{\FMsr(\dsf)}$ to define the obstruction theory of quasimaps via the complex $R\pi_{1*}L\mathsf{f}^*\BT_{\FMsr(\dsf)}$.  Instead, one has to consider the derived enhancement of $\FMsr(\dsf)$  and the corresponding (virtual) tangent complex to obtain a perfect obstruction theory.  We start with the construction of such derived enhancement.

Let $\BR\rT^*\CP ic(X)$ be the (non-rigidified) derived cotangent bundle of the Picard stack of $X$.  Let  $\BR \FM(\dsf)$ be the (non-rigidified) derived stack of Higgs sheaves without a fixed determinant (both exist by \cite{TV}). We fix a (non-stacky) point  $[(L,0)] \in \BR \rT^*\CP ic(X)$ corresponding to a line bundle $L$ of degree $\dsf$ with a zero cotangent vector. We define  $\BR \FMs_{\mathrm{SL} }(\dsf)$ via the following derived fibred product
\[
\begin{tikzcd}[row sep=scriptsize, column sep = scriptsize]
	& \BR \FMs_{\mathrm{SL} }(\dsf) \arrow[d]  \arrow[r,hook] &  \BR \FM(\dsf) \arrow[d,"(\mathrm{det}\text{,}\tr)"]\\
	& \Spec(\BC)  \arrow[r, "\textit{(L,}0)", hook] & \BR \mathrm T^*\CP ic(X)
\end{tikzcd}
\]
 The classical truncation $\FMs_{\mathrm{SL} }(\dsf))=\mathrm{t}_0(\BR \FMs_{\mathrm{SL} }(\dsf))$ is a $\BZ_{\rsf}$-gerbe over 
 $\FMsr(\dsf)$ (see Section \ref{gerbeclasses}), this is because on the level of $\BC^*$-stabilisers, the map $(\mathrm{det}\text{,}\tr)$ is given by 
\[ \BC^*\xrightarrow{z\rightarrow z^\rsf} \BC^*.\]  
We define  
 \[\BR \FMsr := \BR \FMs_{\mathrm{SL} }(\dsf) \thickslash \BZ_\rsf,\] 
 it is a derived enhancement of  $\FMsr(\dsf)=\FMs_{\mathrm{SL} }(\dsf) \thickslash \BZ_\rsf$. 
 
 On  $\BR \FMs_{\mathrm{SL} }(\dsf)$, we have an exact triangle of tangent complexes associated to the fiber diagram, 
\begin{equation} \label{sequence}
	\BT_{ \BR \FMs_{\mathrm{SL} }(\dsf) } \rightarrow \BT_{ \BR \FM(\dsf)}  \rightarrow \BT_{\BR \mathrm T^*\CP ic(X)} \rightarrow. 
\end{equation}
Once pull-backed via the natural inclusion $\iota\colon \FMs_{\mathrm{SL} }(\dsf) \hookrightarrow \BR \FMs_{\mathrm{SL} }(\dsf)$,  the triangle (\ref{sequence}) becomes a triangle on $\FMs_{\mathrm{SL} }(\dsf)$,
\begin{equation} \label{sequence2}
	\BT_0^{\mathrm{vir}} \rightarrow \BT^{\mathrm{vir}} \rightarrow  H^*(\omega_X\oplus\CO_X[1]) \rightarrow,
\end{equation}
such that $\BT_0^{\mathrm{vir}} $ defines a perfect obstruction theory on $\FMs_{\mathrm{SL} }(\dsf)$ (see \cite[Section 1.1]{STV} for the explanation of how derived structure induces an obstruction theory).  In fact, the whole triangle descends to $\FMsr(\dsf)$, such that the descend of $\BT_0^{\mathrm{vir}}$ is the perfect obstruction theory associated to the enhancement $\BR \FMsr(\dsf) $. Abusing the notation, we will denote the descends of the terms in the triangle (\ref{sequence2}) by the same symbols.

 On $\FMsr(\dsf)$, the triangle (\ref{sequence2}) can be made more explicit. Let 
\[( \Gsf,\phi) - X\times \FMsr(\dsf), \quad \pi \colon X\times \FMsr(\dsf) \rightarrow \FMsr(\dsf) \]
be the universal Higgs sheaf with the Higgs field on $X\times \FMsr(\dsf)$ and the canonical projection, respectively.   Then (\ref{sequence2})  fits into the following diagrams of triangles
\begin{equation} \label{diagram1}
 \begin{tikzcd}[row sep=scriptsize, column sep = small]
	&\BT_0^{\mathrm{vir}}[-1]  \arrow[d] \arrow[r] &  \BT^{\mathrm{vir}}[-1] \arrow[r,] \arrow[d] & H^*(\omega_X[-1]\oplus\CO_X)  \arrow[d] \\
	&R\CHom_{\pi}( \Gsf, \Gsf)_0   \arrow[d," \circ \phi - \phi \circ"] \arrow[r] & R\CH om_{ \pi}( \Gsf,\Gsf) \arrow[r,"\tr"] \arrow[d," \circ \phi - \phi \circ"] & H^*(\CO_X)  \arrow[d,"0"]\\
	& R \CHom_{ \pi}( \Gsf,\Gsf \otimes \omega_X)_0 \arrow[r] & R \CHom_{\pi}( \Gsf,\Gsf \otimes \omega_X) \arrow[r,"\tr "]& H^*(\omega_X)
\end{tikzcd} 
\end{equation}
where both vertical sequences and horizontal sequences are distinguished triangles\footnote{This diagram also shows that the triangle (\ref{sequence2}) indeed descends, because the last two rows of the diagram are pullbacked to the same triangles defined for the universal family on $\FMs_{\mathrm{SL}}(\dsf) $.}. We refer to \cite{TT} for more details (in particular, \cite[Corollary 2.26]{TT}) on this diagram in the setting of surfaces, which also applies to curves verbatim.
\\

We will now use the obstruction theory $\BT_0^{\mathrm{vir}}$ to a define map-theoretic obstruction theory of $Q^{\epsilon}_{g,n}(\Ms(\dsf), \wsf)$ in the style of Gromov--Witten theory.    Let 
\[ \mathsf{f}\colon \CC_{g,n}\rightarrow \FMsr(\dsf), \quad \pi_{1}\colon  \CC_{g,n} \rightarrow Q^{\epsilon}_{g,n}(\Ms(\dsf), \wsf)\]
be the universal map and the canonical projection from the universal curve. 
	\begin{prop}\label{obst}
	The complex $R\pi_{1*}L\mathsf{f}^*\BT_0^{\mathrm{vir}}$ defines a relative obstruction theory on $Q^{\epsilon}_{g,n}(\Ms(\dsf), \wsf)$,
	\[\psi\colon (R\pi_{1*}L\mathsf{f}^*\BT_0^{\mathrm{vir}})^{\vee} \rightarrow \BL_{Q^{\epsilon}_{g,n}(\Ms(\dsf), \wsf)/\FM_{g,n}}.\]
\end{prop}
\textit{Proof.} 
The obstruction theory 
\[\psi\colon (R\pi_{1*}L\mathsf{f}^*\BT_0^{\mathrm{vir}})^{\vee} \rightarrow  \BL_{Q^{\epsilon}_{g,n}(\Ms(\dsf), \wsf) / \FM_{g,n}}\]
is given by a derived mapping stack of quasimaps $\BR Q^{\epsilon}_{g,n}( \Ms(\dsf), \wsf)$ to the derived stack $\BR \FMsr(\dsf)$, which exists by Lurie's representability theorem \cite{Lur}.  \qed 
\\

  The obstruction theory $R\pi_{1*}L\mathsf{f}^*\BT_0^{\mathrm{vir}}$ can be given a sheaf-theoretic interpretation, which will be the origin of comparison results between quasimap invariants and Vafa--Witten invariants. We will also need it to prove perfectness of $R\pi_{1*}L\mathsf{f}^*\BT_0^{\mathrm{vir}}$ in Proposition \ref{perf2}.   
  		
  		Let 
\[(\Fsf, \Phi):=\mathsf{f}^*(\Gsf,\phi) - X\times \CC_{g,n},\quad \pi_2\colon X\times \CC_{g,n}  \rightarrow Q^{\epsilon}_{g,n}(\Ms(\dsf), \wsf),\]
be the universal family of Higgs sheaves associated to $Q^{\epsilon}_{g,n}(\Ms(\dsf), \wsf)$ and the canonical projection. Let $\BE_0^{\mathrm{vir}}$ be defined by 
\begin{equation} \label{diagram2}
	\begin{tikzcd}[row sep=scriptsize, column sep = small]
		&\BE_0^{\mathrm{vir}}[-1]  \arrow[d] \arrow[r] &  \BE^{\mathrm{vir}}[-1] \arrow[r,] \arrow[d] & R\pi_{2*}(\omega_{\pi_2}[-1]\oplus\CO) \arrow[d] \\
		&R\CH om_{\pi_2}(\Fsf,\Fsf)_0   \arrow[d," \circ \Phi - \Phi \circ"] \arrow[r] & R\CH om_{\pi_2}(\Fsf,\Fsf) \arrow[r,"\tr"] \arrow[d," \circ \Phi - \Phi \circ"] & R\pi_{2*}\CO\arrow[d,"0"]\\
		& R \CHom_{\pi_2}(\Fsf,\Fsf \otimes \omega_{\pi_2})_0 \arrow[r] & R \CHom_{\pi_2}(\Fsf,\Fsf \otimes \omega_{\pi_2}) \arrow[r,"\tr "]& R\pi_{2*}\omega_{\pi_2}
	\end{tikzcd} 
\end{equation}
where both vertical and horizontal sequences are distinguished triangles (see \cite[Corollary 2.26]{TT} for more details). 
	\begin{prop}\label{comparison}
		The complex $R\pi_{1*}L\mathsf{f}^*\BT_0^{\mathrm{vir}}$ is canonically isomorphic to the complex $\BE_0^{\mathrm{vir}}$.
	\end{prop}
	
	\textit{Proof.}  Consider the following diagram
	\[
	\begin{tikzcd}[row sep=scriptsize, column sep = scriptsize]
		X\times \CC_{g,n} \arrow[r,"\mathrm{id}\times \mathsf{f}"] \arrow[d] \arrow[dd, bend right=80,"\pi_2" near start]
		& X \times \FMsr(\dsf)\arrow[d,"\pi"] \\
		\CC_{g,n} \arrow[r,"\mathsf{f}"] \arrow[d, "\pi_{1}"] 
		& \FMsr(\dsf) \\
		Q^{\epsilon}_{g,n}(\Ms(\dsf), \wsf)
	\end{tikzcd}
	\]
	By the functionality of the trace map and the base change theorem applied to the diagram, all the complexes in the last two rows of (\ref{diagram1}) transform to the complexes in (\ref{diagram2}) with respect to the functor $R\pi_{1*}L\mathsf{f}^*$. The result then follows from the fact that the vertical sequences are also distinguished triangles.  For more details, we refer to \cite[Proposition 5.5]{N}. 
	\qed

	\begin{prop}\label{perf2}
	The complex $R\pi_{1*}L\mathsf{f}^*\BT_0^{\mathrm{vir}}$  is of amplitude $[0,1]$. In particular, the obstruction theory
		\[\psi\colon (R\pi_{1*}L\mathsf{f}^*\BT_0^{\mathrm{vir}})^{\vee} \rightarrow \BL_{Q^{\epsilon}_{g,n}(\Ms(\dsf), \wsf)/\FM_{g,n}}\]
		is perfect.
	\end{prop}
\textit{Proof.} It is clear that the complex   $R\pi_{1*}L\mathsf{f}^*\BT_0^{\mathrm{vir}}$  is at most of amplitude $[-1,2]$. So we need to show that it is zero in degrees $-1$ and $2$. It is enough to check it over a point $[f\colon C\rightarrow \FMsr(\dsf)]\in Q^{\epsilon}_{g,n}(\check{M}(\dsf), \wsf)$. By Proposition \ref{comparison2}, 
\[H^{-1}(f^*\BT_0^{\mathrm{vir}},C)=\Hom(F,F)_0,\]
where $(F,\phi)$ is the family of Higgs sheaves associated to the quasimap $f$. By \cite[Lemma 3.17]{N}, the family $(F,\phi)$ is stable (see also Lemma \ref{stability}), we conclude that $H^{-1}(f^*\BT_0^{\mathrm{vir}})=\Hom(F,F)_0=0$.

 Consider now the distinguished triangle 
\begin{equation}  \label{triangle}
	\tau_{\leq 0}f^*\BT_0^{\mathrm{vir}} \rightarrow f^*\BT_0^{\mathrm{vir}} \rightarrow \tau_{\geq 1} f^*\BT_0^{\mathrm{vir}}\rightarrow ,
\end{equation}
where $\tau_{\dots}$ is the truncation of complex with respect to the standard $t$-structure. Taking the long exact sequence of cohomologies associated to (\ref{triangle}), we obtain
\[\dots \rightarrow H^2(\tau_{\leq 0}f^*\BT_0^{\mathrm{vir}},C) \rightarrow H^2(f^*\BT_0^{\mathrm{vir}}, C) \rightarrow H^2(\tau_{\geq 1} f^*\BT_0^{\mathrm{vir}},C).\]
Let us now analyse the terms in the long exact sequence.  Firstly,  $f$ maps generically to the stable locus $\check{M}(\dsf)$, over which $\BT_{0}^{\mathrm{vir}}$ is a locally free sheaf concentrated in degree 0. Hence the object $f^*\BT_0^{\mathrm{vir}}$ is  generically a locally free sheaf concentrated in degree 0. In particular, this implies that $\tau_{\geq 1} f^*\BT_0^{\mathrm{vir}}$ is a $0$-dimensional sheaf in degree 1 supported on the base points of $f$, therefore 
\[ H^2(\tau_{\geq 1} f^*\BT_0^{\mathrm{vir}},C)=0.\]
Since $C$ is a curve, 
\[ H^2(\tau_{\leq 0}f^*\BT_0^{\mathrm{vir}},C) =0.\] 
By the long exact sequence, we therefore obtain that $H^2(f^*\BT_0^{\mathrm{vir}}, C)=0$. \qed

	\subsection{Invariants} We will now define the quasimap invariants. Since only the fixed locus $Q^{\epsilon}_{g,n}(\Ms(\dsf), \wsf)^{\BC_t^*}$ is proper, we define the invariants via localisation. By \cite{GP},  the fixed part of the obstruction theory from Proposition \ref{obst} defines the virtual fundamental class
	\[[Q^{\epsilon}_{g,n}(\Ms(\dsf), \wsf)^{\BC_t^*}]^{\mathrm{vir}} \in H_*(Q^{\epsilon}_{g,n}(\Ms(\dsf), \wsf)^{\BC_t^*},\BQ), \]
	while the moving part defines a virtual normal bundle $N^{\mathrm{vir}}$ of $Q^{\epsilon}_{g,n}(\Ms(\dsf), \wsf)^{\BC_t^*}$ inside $Q^{\epsilon}_{g,n}(\Ms(\dsf), \wsf)$ and therefore the corresponding equivariant  Euler class, 
	\[e_{\BC_t^*}(N^{\mathrm{vir}})\in H^*(Q^{\epsilon}_{g,n}(\Ms(\dsf), \wsf)^{\BC_t^*},\BQ)[t^{\pm}].\]
	Let us recall all the necessary structures needed for the definition of invariants:
	\begin{itemize}
		\item evaluation maps at marked points
		\[\mathrm{ev}_{i}\colon  Q^{\epsilon}_{g,n}(\Ms(\dsf), \wsf) \rightarrow \Ms(\dsf), \quad i=1, \dots, n;\]
		\item projection to a moduli space of stable curve
		\[ \pi \colon Q^{\epsilon}_{g,n}(\Ms(\dsf), \wsf) \rightarrow \Mbar_{g,n};\]
		\item cotangent line bundles 
		\[\CL_{i}: =s^{*}_{i}(\omega_{\CC_{g,n}/Q^{\epsilon}_{g,n}(\Ms(\dsf),\wsf)}), \quad i=1, \dots, n,\]
	\end{itemize}
	where $s_{i}\colon  Q^{\epsilon}_{g,n}(\Ms(\dsf),\wsf) \rightarrow \CC_{g,n}$ are universal markings. We denote 
	\[\psi_{i}:=\mathrm{c}_{1}(\CL_{i}) \in H_{\BC_t^*}^*((Q^{\epsilon}_{g,n}(\Ms(\dsf),\wsf),\BQ), \quad i=1,\dots, n.\]
	
	\begin{defn} \label{invariants} The \textit{descendent} $\epsilon$-\textit{invariants} are 
		\begin{multline*}
			\langle \alpha; \psi_1^{m_{1}}\gamma_{1}, \dots, \psi_n^{m_{n}}\gamma_{n} \rangle^{\epsilon}_{g,n,\wsf}:=\\
			\int_{[Q^{\epsilon}_{g,n}(\Ms(\dsf),\wsf)^{\BC_t^*}]^{\mathrm{vir}}}\frac{\pi^* \alpha \cdot \prod^{i=n}_{i=1}\psi_{i}^{m_{i}}\ev^{*}_{i}(\gamma_{i},)}{e_{\BC_t^*}(N^{\mathrm{vir}})}\in \BQ[t^{\pm}],
		\end{multline*} 
		where $\gamma_{1}, \dots, \gamma_{n} \in H_{\BC_t^*}^{*}(\Ms, \BQ)$ are some classes, $\alpha \in R^*(\Mbar_{g,n}, \BQ)$ is a tautological class and $m_{1}, \dots m_{n}$ are non-negative integers. 
		
	\end{defn}

\subsection{Genus 1 quasimap invariants} \label{genus1sl} 
Let us now consider genus 1  invariants of $\Ms(\dsf)$. Virtual dimension of the moduli space $Q^{\epsilon}_{1,1}(\Ms(\dsf),\wsf)$  is1, i.e.\ an insertion is needed to get an invariant. We will take a tautological insertion corresponding to a non-stacky point $(E,0) \in \Mbar_{1,1}$, where $(E,0)$ is marked smooth genus 1 curve (an elliptic curve). Equivalently, we can consider the following moduli space. 

\begin{defn} \label{fixedC}If $\wsf\neq 0$, let
	\[Q^{\epsilon}_{E}(\Ms(\dsf),\wsf) \colon (Sch/ \BC) \rightarrow (Grpd) \]
	 be the moduli space of $\epsilon$-stable quasimaps from $E$ (without identifications by automorphisms of $E$). If $\wsf=0$,  by $Q^{0^+}_{E}(\Ms(\dsf),0)$ we will denote the moduli space of degree $0$ maps from $E$ to $\Ms(\dsf)$ . 
\end{defn}
The virtual dimension of $Q^{\epsilon}_{E}(\Ms(\dsf),\wsf)$  is 0 and the intersection theory of this space is the equivalent to the one of $Q^{\epsilon}_{1,1}(\Ms(\dsf),\wsf)$ with a tautological insertion corresponding to $\alpha=[(E,0)]$. 
 We now sketch a construction of an equivariant surjective cosection for $Q^{\epsilon}_{E}(\Ms(\dsf),\wsf)$, which forces us to consider a further modification of the moduli space $Q^{\epsilon}_{E}(\Ms(\dsf),\wsf)$. In what follows, we denote
  \[ \mathbf t := \text{standard representation of } \BC_t^* \text{ on } \BC.\]
\begin{prop} \label{cosectionelliptic}
	If  $\wsf\neq 0$, then the obstruction theory of $Q^{\infty}_{g,n}(\Ms(\dsf),\wsf)$ has a $\BC^*_t$-equivariant surjective cosection
	\[h^1(R\pi_{1*}L\mathsf{f}^*\BT_0^{\mathrm{vir}}) \rightarrow \mathbf t\CO, \]
	So does the obstruction theory of $Q^{0^+}_{E}(\Ms(\dsf),\wsf)$, 
	\[h^1(R\pi_{1*}L\mathsf{f}^*\BT_0^{\mathrm{vir}}) \rightarrow \mathbf t\CO.\]
\end{prop}
\textit{Proof.} The proof of the proposition is different for different values of $\epsilon$. If $\epsilon=\infty$, i.e.\ if we are in the setting of Gromov--Witten theory of $\Ms(\dsf)$, we use the argument of \cite[Section 3.4]{OP10b}. We will show it over a point $[f\colon C\rightarrow \Ms(\dsf)] \in Q^{\infty}_{g,n}(\Ms(\dsf),\wsf)$. Consider the natural composition of maps
\begin{equation} \label{mapcosection}
	f^*\Omega_{\Ms(\dsf)} \rightarrow \Omega_C \rightarrow \omega_C,
\end{equation}
whose cokernel is 0-dimensional. Applying the equivariant isomorphism 
\[\Omega_{\Ms(\dsf)} \cong \mathbf t^{-1} T_{\Ms(\dsf)},\] 
we obtain the map
\[f^*T_{\Ms(\dsf)} \rightarrow \mathbf t \omega_C, \]
taking cohomology, we obtain a surjective map
\[H^1(f^*T_{\Ms(\dsf)}, C) \rightarrow \mathbf t \BC,\]
which can be easily globalised to give the desired cosection. This proves the claim for $\epsilon=\infty$. 

Let us now deal with $\epsilon=0^+$. It is tempting to apply the same argument to any quasimap $f\colon C\rightarrow \FMsr(\dsf)$. However, the natural map (\ref{mapcosection}) might be zero for a general quasimap. Indeed, a quasimap can be generically constant and have a base point, in this case (\ref{mapcosection}) is zero.  Hence for $\epsilon=0^+$, we have to follow a different approach - using the identification from Proposition \ref{comparison}, we can apply the argument from \cite[Section 1.3]{O}. Working again over a point $[f\colon C \rightarrow \FMsr(\dsf)]  \in Q^{0^+}_{E}(\Ms(\dsf),\wsf)$,  let $\CF$ be the associated 2-dimensional sheaf on $\rT^*X\times E$.  Translations by $E$ gives an infinitesimal deformation 
\[\BC v \subset  \Ext_{\rT^*X\times E}^1(\CF,\CF), \]
applying Serre's duality and the equivariant isomorphism 
\[\omega_{\rT^*X\times E} \cong \mathbf t^{-1} \CO_{\rT^*X\times E},\]
 we obtain a surjective map
\[ \Ext_{\rT^*X\times E}^2(\CF,\CF) \cong \mathbf t\Ext_{\rT^*X\times E}^1(\CF,\CF)^{\vee}\rightarrow \mathbf{t}(\BC v)^{\vee},\]
whose globalisation gives the desired cosection. 
\qed 
\begin{rmk} In the case of moduli spaces of sheaves on K3 surfaces, a surjective cosection is constructed for any moduli space of $\epsilon$-stable quasimaps $Q^{\epsilon}_{g,n}(M,\beta)$, using the semiregularity map, \cite{NK3}. However, we do not know how to adopt the similar approach to the case of Higgs sheaves.  
	\end{rmk}

If $\wsf \neq 0$, then equivariant cosection forces the virtual fundamental class of $Q^{\epsilon}_{E}(\Ms(\dsf),\wsf)$ for $\epsilon=0^+/\infty$ to be of the following form
\[  [Q^{\epsilon}_{E}(\Ms(\dsf),\wsf))]^{\mathrm{vir}}= t[Q^{\epsilon}_{E}(\Ms(\dsf),\wsf)]^{\mathrm{red}}, \] 
where $t$ is the equivariant parameter and $[Q^{\epsilon}_{E}(\Ms(\dsf),\wsf)]^{\mathrm{red}}$ is the reduced virtual fundamental class associated to the cone of the cosection, which is of homological degree 1.  Hence the reduced virtual dimension of $Q^{\epsilon}_{E}(\Ms(\dsf),\wsf)$ is actually 1. To get interesting enumerative invariants, we make the following definition. 
\begin{defn}\label{quotE} If $\wsf\neq 0$,  let 
\[Q^{\epsilon}_{E}(\Ms(\dsf),\wsf)^{\bullet} \colon (Sch/ \BC)^{\circ} \rightarrow (Grpd) \]
be the moduli space of $\epsilon$-stable qausimaps from $E$ up to translations by $E$. 
\end{defn}
  
   Note the following relation between $Q^{\epsilon}_{E}(\Ms(\dsf),\wsf)$ and $Q^{0^+}_{E}(\Ms(\dsf),\wsf)^{\bullet}$, 
 \begin{equation} \label{Eq}
 	[Q^{\epsilon}_{E}(\Ms(\dsf),\wsf)/E] \cong Q^{\epsilon}_{E}(\Ms(\dsf),\wsf)^{\bullet},
 	\end{equation}
 where the quotient is taken with respect to the natural action of $E$ on quasimaps given by composing a map with a translation of $E$.  The moduli space $Q^{\epsilon}_{E}(\Ms(\dsf),\wsf)^{\bullet}$ also carries a perfect obstruction theory and a $\BC_t^*$-action with a proper fixed locus. The following invariants will be the principle subject of our study. 

\begin{defn} \label{qmsell} Let $\asf \in \BZ_\rsf$. If  $\wsf \neq 0$, we define
	\[\QMs^{\asf, \bullet}_{\dsf, \wsf}= 	\begin{cases}
		\int_{[Q^{0^+}_{E}(\Ms(\dsf),\wsf)^{\bullet, \BC_t^*}]^{\mathrm{vir}}}\frac{1}{te_{\BC_t^*}(N^{\mathrm{vir}})}  & \text{if}\ \wsf=\dsf\cdot \asf \mod \rsf \\
		0 &\text{otherwise.} 
	\end{cases} \]
	
	If $\wsf=0$, we define 
	\[\QMs^{\asf}_{\dsf,0}=
	\begin{cases}
 \int_{[Q^{0^+}_{E}(\Ms(\dsf),0)^{ \BC_t^*}]^{\mathrm{vir}}}\frac{1}{e_{\BC_t^*}(N^{\mathrm{vir}})}  & \text{if}\ \asf=0 \\
	0 &\text{otherwise.} 
\end{cases}\] 
Similarly, we define genus 1 Gromov-Witten invariants for $\wsf \neq 0$,
\[\GWs^{\asf, \bullet}_{\dsf, \wsf}= 	\begin{cases}
	\int_{[Q^{\infty}_{E}(\Ms(\dsf),\wsf)^{\bullet, \BC_t^*}]^{\mathrm{vir}}}\frac{1}{te_{\BC_t^*}(N^{\mathrm{vir}})}  & \text{if}\ \wsf=\dsf\cdot \asf \mod \rsf \\
	0 &\text{otherwise.} 
\end{cases} \]
Whenever we do not need the class $\asf$, we drop it from the notation, simply writing $\QMs^{\asf, \bullet}_{\dsf, \wsf}$, $\GWs^{\bullet}_{\dsf, \wsf}$ and $\QMs_{\dsf,0}$. 

	\end{defn}

\begin{rmk} \label{asfnot}
	The notation involving $\asf$ can be confusing. However, it will be of fundamental importance for the statement of Enumerative mirror symmetry. It also has a geometric meaning, $\asf$ is the pullback of the gerbe class $\check \alpha$ by a quasimap $f$, i.e.\  $\asf=f^*\check \alpha \in H^2(E,\BZ_\rsf) \cong \BZ_\rsf$. In the case of $\Ms(\dsf)$ with $\dsf\neq 0$,  the class $\asf$ can be recovered from the data of $\wsf$ and $\dsf$, since  $[\ch_1(F)_{\mathrm d}]=f^*\check \alpha$ and $\wsf= \dsf\ch_1(F)_{\mathrm d}-\rsf \ch_2(F)_{\mathrm d}$ (see Lemma \ref{degreechern}). This might be viewed as  the consequence of the fact that $H^2(\Ms(\dsf),\BZ)\cong \BZ[\Thetas]$, which implies that $\check{\asf}$ is a multiple of $[\Thetas]$ in $H^2(\Ms(\dsf),\BZ_\rsf)$. In the case of $\Ms(0)$ or $\Mp(\dsf)$, the class $\alphas$ is an extra piece of information and it is incorporated into the definition of the degree of a quasimap, see Definition \ref{extendeddegree}. 
	\end{rmk}
These invariants can be seen as (virtual) counts of genus 1 curves of a fixed complex structure inside $\Ms(\dsf)$ up to translations of the domain curve. In the case of moduli spaces of sheaves on K3 surfaces, similar invariants are expected to be the actual (enumerative) number of curves, e.g.\ see \cite[Conjecture 1.1]{NO} and \cite[Corollary 1.5]{NO}. 
	
\section{$\mathrm{PGL}_\rsf$-quasimaps, $\dsf\neq 0$}
The moduli space $\Mp(\dsf)$ is an orbifold, we therefore have to allow the curves to be twisted, (see \cite[Section 2.1]{ACV} for the definition).

\begin{defn} \label{quasimapsdef2} 
	A map $f\colon (C,\mathbf x) \rightarrow \FMp(\dsf) $ is a \textit{quasimap} to $(\Mp(\dsf), \FMp(\dsf))$ of genus $g$ and of degree $\wsf \in \BZ$,  if 
	\begin{itemize}
		\item $(C,\mathbf x)$ is a connected twisted marked nodal curve of genus $g$ and $f$ is representable,
		\item $\deg(f):=\deg(f^*[\Thetap])=\wsf$,
		\item $| \{ p\in C \mid f(p) \in \FMp(\dsf) \setminus \Mp(\dsf) \} | < \infty$.
	\end{itemize}
	We will refer to the set $\{ p\in C\mid f(p) \in \FMp(\dsf) \setminus \Mp(\dsf) \}$ as \textit{base} points. A quasimap $f$ is \textit{prestable}, if 
	\begin{itemize}
		\item $\{ t\in C \mid f(t) \in \FMp(\dsf) \setminus \Mp(\dsf) \} \cap \{\mathrm{nodes}, \mathbf x \}=\emptyset$. 
	\end{itemize} 
\end{defn}

To avoid working with azalgebras as much as possible, we make the following definition. 

\begin{defn} \label{torsorcurve} Let $f\colon C \rightarrow \HAz(\dsf)$ be a quasimap, we define 
	\[h\colon \tilde C \rightarrow C \quad \text{and} \quad   \tilde f \colon \tilde C \rightarrow \FMsr(\dsf)\] by the following fiber product
	\[
	\begin{tikzcd}[row sep=scriptsize, column sep = scriptsize]
		& \tilde C \arrow[d, "h"]  \arrow[r,"\tilde f"] &  \FMsr(\dsf) \arrow[d]\\
		& C \arrow[r, "f"] &  \HAz(\dsf)
	\end{tikzcd}
	\]
\end{defn}

\begin{lemma} \label{pgldeg}
	Degree of a quasimap $f\colon C \rightarrow \HAz(\dsf)$ is determined by $\tilde f$,
	\[\deg(f)= \deg(\tilde f) /\rsf^{2g}.\]
\end{lemma}
\textit{Proof.}  The cardinality of $\Gamma_X$ is $\rsf^{2g}$, hence so is the degree of the map $h$ by the definition. We therefore obtain
\[\deg(f)=\deg(f^*[\Thetap])=\deg(\tilde f^*\Thetas) /\rsf^{2g}=\deg(\tilde f) /\rsf^{2g}.\]
\qed 

We will need a refined version of the degree for quasimaps to $\Mp(\dsf)$ (cf. Remark \ref{asfnot}). However, due to the limited knowledge of the author, we will define it only for smooth curves, which will be enough for the purposes of the present work. 
\begin{defn} \label{extendeddegree}
	Let $C$ be a smooth (non-twisted) curve. Recall the natural identification $H^2_{\text{\'et}}(C,\BZ_\rsf)\cong \BZ_\rsf$ and the fact that $H^2_{\text{\'et}}(C,\BZ_\rsf)$ classifies $\BZ_\rsf$-gerbes. A quasimap $f \colon C \rightarrow \FMp(\dsf)$  is said to be of extended degree $(\wsf,\asf) \in \BZ\oplus \BZ_\rsf$, if
	\begin{itemize}
	 \item $\deg(f)=\wsf$, 
	 	 \item $f^*\hat \alpha=\asf \in \BZ_\rsf,$
	 \end{itemize}
 where $\hat \alpha$ is defined in Section \ref{gerbeclasses}.
	\end{defn}
\begin{rmk}
	The above definition can be easily extended to a nodal curve, but since we are in the orbifold setting, it is not enough, if we want to define the extended degree in general. The reason we cannot make the definition above for a twisted curve is that the author does not understand \'etale cohomology of a twisted curve. In particular, the author does not know how to construct a natural degree map $\deg \colon H^2_{\text{\'et}}(C,\BZ_\rsf)\rightarrow \BZ_\rsf$, if $C$ is a twisted curve. The problem is that $H^2_{\text{\'et}}(C,\BZ_\rsf)$ might contain the group cohomology $H^2(\Gamma, \BZ_\rsf)$ of the isotropy group $\Gamma$ of some point.  
	
	Luckily, a moduli space of quasimaps from a fixed smooth curve is proper, hence the curve does not need to degenerate.  We will be mainly concerned with quasimap invariants associated to a fixed curve, so it is good enough for our purposes.  
	\end{rmk}
\begin{rmk}
	We expect that the class $\asf$ corresponds to the torsion in Borel-Moore homology $H^{BM}_2(\Mp(\dsf),\BZ)$ in the case of maps to $\Mp(\dsf)$. That would also agree with the fact that $\asf$ can be recovered from the standard degree in the case of $\Ms(\dsf)$, because $H^{BM}_2(\Ms(\dsf),\BZ)$ does not have torsion. 
	\end{rmk}
Assume $C$ is a smooth (non-twisted) curve. Let $(\CA, \phi)$ be a Higgs azalgebra associated to a quasimap $f\colon C\rightarrow \FMp(\dsf)$ of extended degree $(\wsf, \asf)$.  Since $(\CA, \phi)$ is a family of rank $\rsf$ Higgs azalgebras, it is of rank $\rsf$ itself. Let us now consider its degree. As in Section \ref{slpre},  the degree $\delta(\CA)$ has the following components with respect to Kunneth's decomposition of $H^2(X\times C, \BZ_\rsf)$,
\begin{equation*}
	\delta(\CA)=(\delta(\CA)_{\mathrm f}, \delta(\CA)_{\mathrm m}, \delta(\CA)_{\mathrm d}) \in \BZ_\rsf \oplus \Lambda_{\rsf,C} \oplus \BZ_\rsf,
\end{equation*}
where 
\[ \Lambda_{\rsf,C} := H^1(X,\BZ_\rsf)\otimes H^1(C,\BZ_\rsf).\]
As before, we have 
\[\delta(\CA)_{\mathrm f}=\mathsf d.\]
However, unlike in the case of the Higgs sheaves, degree component $\delta(\CA)_{\mathrm d}$ is not accessible directly from the degree $\deg(f)$, as it was shown in Lemma \ref{lemmaeqchern}, because the \textit{middle} component  $\delta(\CA)_{\mathrm m}$ might be non-zero. For that reason, we introduced the notion of extended degree, as it recovers the degree component 
\[\delta(\CA)_{\mathrm d}=\mathsf a\]
by construction. On the other hand, the standard degree $\deg(f)$ admits a similar expression as the one in Lemma \ref{degreechern} (notice that $\mathrm{c}_2(R\CE nd(F,F))/2=\mathrm c_1^2(F)/2-\rsf \mathrm{ch}_2(F)$).
\begin{lemma}
	\[\deg(f)=\mathrm{c}_2(\CA)/2\]
\end{lemma} 
\textit{Proof.} Consider the map $\tilde{f}$ from Definition \ref{torsorcurve}. 
Let $\tilde{F}$ be the sheaf associated to $\tilde{f}$ and $\CA$ be the azalgebra associated to $f$, then  by the definition of $\tilde f$, 
\[ h^*\CA\cong R\CE nd(\tilde F, \tilde F),\]
in particular, 
\begin{align*}
	\mathrm{c}_2(A)/2=\mathrm{c}_2(R\CE nd(\tilde F, \tilde F))/2\rsf^{2g}&=(\mathrm c_1^2(F)-2\rsf \mathrm{ch}_2(F))/2\rsf^{2g}\\
	&=\deg(\tilde f)/\rsf^{2g},
\end{align*}
where the last equality follows from Lemma \ref{degreechern}. By Lemma \ref{pgldeg}, we have
\[\deg(f)= \deg(\tilde f) /\rsf^{2g},\]
we therefore obtain 
\[\deg(f)=\mathrm{c}_2(\CA)/2.\]
\qed
\\

The middle part $\delta(\CA)_{\mathrm m}$ is undetermined by topological data associated to a quasimap.  Hence we will have to add up contributions from different values of $\delta(\CA)_{\mathrm m}$. In summary, the associated azalgebra $\CA$ has the following topological invariants. 
\begin{lemma} \label{topazal}
	 Let $(\CA, \phi)$ be a Higgs azalgebra associated to a quasimap $f\colon C\rightarrow \FMp(\dsf)$ of extended degree $(\wsf, \asf)$, then
\begin{align*}
	&\rk(\CA)=\rsf \\
	&\delta(\CA)=(\dsf, \delta(\CA)_{\mathrm m}, \asf ) \\
	&\mathrm{c}_2(\CA)=2\wsf. 	
\end{align*}
\end{lemma}
 
\subsection{Moduli spaces of quasimaps}

Given a quasimap $f \colon C \rightarrow \HAz(\dsf)$ and a point $p \in C$. By Langton's semistable reduction, we can modify the quasimap $\tilde f$ at all points in the fiber $h^{-1}(p)$, to obtain a quasimap  
\[\tilde{f}_p \colon \tilde C \rightarrow  \FMsr(\dsf)\]
which maps to the stable locus $\Ms(\dsf)\subset \FMsr(\dsf)$ at $h^{-1}(p)$ (if $p$ is not a base point of $f$, then $\tilde{f}_p=\tilde f$).   

\begin{defn}
	We define the \textit{length} of a point $p \in C$ to be 
	\[\ell(p):=\deg(f)-\deg(\tilde f_p)/\rsf^{2g}.\]
	By the proof of \cite[Proposition 3.7]{N}, $\ell(p)\geq 0$ and $\ell(p)=0$, if and only if $p$ is not a base point. 
\end{defn}
\begin{defn} \label{stabilityqm2} Given $\epsilon \in \BR_{>0}\cup\{0^+, \infty\}$, a prestable quasimap $f\colon(C, \mathbf x) \rightarrow  \HAz(\dsf)$ is $\epsilon$-stable, if 
	\begin{itemize}
		\item[$\mathbf{(i)}$] $\omega_{C}(\mathbf x) \otimes f^*[\Thetap]^\epsilon$ is positive;
		\item[$\mathbf{(ii)}$] $\epsilon \ell(p)\leq 1$ for all $p \in C$.
	\end{itemize}
	
\end{defn}	
\noindent Let
\begin{align*}
	Q^{\epsilon}_{g,n}(\Mp(\dsf), \wsf) \colon (Sch/ \BC)^{\circ} &\rightarrow (Grpd) \\
	B&\mapsto \{\text{families of $\epsilon$-stable quasimaps over }B\}
\end{align*}
be the moduli space of $\epsilon$-stable quasimaps of genus $g$ and degree $\wsf$ with $n$ marked points. The extended degree is incorporated in the setting of a fixed smooth curve in Definition \ref{extendeddegree}.  

The $\BC_t^*$-action on $\Mp(\dsf)$ given by the scaling of Higgs fields induces a $\BC_t^*$-action on $Q^{\epsilon}_{g,n}(\Mp(\dsf), \wsf)$. As before $\Mp(\dsf)$ is not proper but the fixed locus $\Mp(\dsf)^{\BC_t^*}$, hence this action will be crucial for defining enumerative invariants. 


\begin{lemma} \label{extension2}
Let $\CC \rightarrow \Delta$ be a family of twisted nodal curves over a DVR $\Delta$ and $\{p_i\}\subset \CC$ be finitely many closed points in the regular locus of the central fiber.  Any quasimap 
$f^{\circ} \colon \CC^{\circ}=\CC \setminus \{p_i\} \rightarrow \HAz(\dsf)^{\BC_t^*}$ extends to a quasimap $f\colon\CC \rightarrow \HAz(\dsf)^{\BC_t^*}$, which is unique up to unique isomorphism. 
\end{lemma}
\textit{Proof.} The strategy is similar to the one used in the proof of Lemma \ref{extension2}, where we used Langton's algorithm \cite{Lang} to produce an extension. We could directly apply Langton's algorithm to an Azumaya algebra, but then we have to know that we get an Azumaya algebra at each step of the algorithm. Moreover, we also have to take care of the Higgs field. To circumvent  these technicalities, we use the following trick.

Let $\tilde{\CF}^{\circ}$ be the sheaf on $\rT^*X\times \tilde \CC^{\circ}$ associated to the quasimap $\tilde f^{\circ}$. By definition, $\tilde f^{\circ}$ is $\Gamma_X$-equivariant, hence the sheaf  $\tilde{\CF}^{\circ}$ transforms in the following way with respect to the $\Gamma_X$-action on $\rT^*X\times \tilde \CC^{\circ}$, 
\[\rho_L\colon g_L^*\tilde \CF \cong \tilde \CF\otimes p_X^*L,\]
where $g_L\colon \rT^*X\times \tilde \CC^{\circ} \cong \rT^*X\times \tilde \CC^{\circ}$ is the action associated to the line bundle $L \in \Gamma_X$. The maps $\rho_L$ are compatible in the obvious sense. Let us call the collection of such compatible maps 
\[\{\rho_L \mid L \in \Gamma_X \}\]
a \textit{twisted equivariant} structure on $\tilde{\CF}^{\circ}$. Conversely, if we have a pair $(\tilde \CC, \tilde\CF)$, where $h\colon \tilde \CC \rightarrow C$ is a $\Gamma_X$-torsor  and $\tilde\CF$ is a twisted equivariant sheaf, then $\tilde\CF$ defines a $\Gamma_X$-equivariant quasimap $\tilde f \colon \tilde \CC \rightarrow \FMsr(\dsf)$, which in turn defines a quasimap $f\colon \CC \rightarrow \HAz(\dsf)$. 
In order to establish the claim, we therefore have to show that each step of Langton's algorithm    used in Lemma \ref{extension0} is twisted equivariant.  

Firstly, we extend the $\Gamma_X$-torsor $h^{\circ}\colon \tilde \CC^{\circ} \rightarrow \CC^{\circ}$ to $\Gamma_X$-torsor $h\colon \tilde \CC \rightarrow \CC$. This is possible, because $\CC^{\circ}$ is a complement of finitely many closed points in the regular (non-stacky) locus of the surface $\CC$ and $\tilde \CC$ is a finite torsor. 


Secondly, we extend the sheaf $\tilde\CF^{\circ}$ to a point in the fiber $ h^{-1}(p_i)$. By acting on this extension, we can translate the extension to other points in the fiber $ h^{-1}(p_i)$, thereby obtaining a twisted equivariant extension $\tilde \CF$. We quotient away the torsion if necessary (the sheaf remains twisted equivariant). As in Lemma \ref{extension0}, we now want to get rid of impurity (subsheaves of lower dimension) in the central fiber $\tilde \CF_{\Spec \BC }$. 
Langton's algorithm  gives modifications of $\tilde \CF$ in the form of a finite sequences of short exact sequences
\begin{align*}
0 \rightarrow \tilde \CF^{1}\rightarrow & \tilde \CF^0 \rightarrow Q^{1} \rightarrow  0, \\
&\vdots \\
0\rightarrow \tilde \CF^{k} \rightarrow & \tilde \CF^{k-1} \rightarrow Q^{k} \rightarrow 0, 
\end{align*}
where $\tilde \CF^0=\tilde \CF$, the central fiber of $\tilde \CF^{k}$ is pure and sheaves $Q^{i}$ are quotients of $\tilde \CF^{i-1}_{\Spec \BC}$ by maximal subsheaves of lower dimension.  A map $\tilde \CF^{i-1} \rightarrow Q^{i}$  is given by the composition
\[\tilde \CF^{i-1} \rightarrow \tilde \CF^{i-1}_{\Spec \BC} \rightarrow  Q^{i}. \]
If  $\tilde \CF^{i-1}$ is twisted equivariant, then each map in the composition is also twisted equivariant,

\[
\begin{tikzcd}[row sep=scriptsize, column sep = scriptsize]
& g_L^*\tilde \CF^{i-1} \arrow[d,"\cong" ]  \arrow[r] & g_L^*\tilde \CF^{i-1}_{\Spec \BC} \arrow[r] \arrow[d, "\cong"]& g_L^* Q^{i} \arrow[d, "\cong"]\\
&\tilde \CF^{i-1}\otimes p^*_XL \arrow[r] & \tilde \CF^{i-1}_{\Spec \BC} \otimes p^*_XL \arrow[r] &Q^{i}\otimes p^*_XL
\end{tikzcd}
\]
By induction on $i$,  each step in Langton's algorithm produces a twisted equivariant sheaf $\tilde \CF^{i-1}$. Thus the last sheaf  $\tilde \CF^{k}$ is also a twisted equivariant sheaf and therefore defines an equivariant quasimap  $\tilde f \colon \tilde \CC \rightarrow \FMsr(\dsf)$, which descends to a quasimap $f \colon \CC \rightarrow \HAz(\dsf)$.  The quasimap $f$ is the desired extension of $f^\circ$. The rest is as in Lemma \ref{extension0}. 
\qed 

\begin{prop} \label{positive2}
	A prestable quasimap $f\colon C \rightarrow \HAz(\dsf)$ is non-constant, if and only if
	\[ \deg(f)>0.\]
\end{prop}  
\textit{Proof.} The quasimaps $f$ is constant, if and only if $\tilde f$ is constant. The claim therefore follows from Proposition \ref{positive} and Lemma \ref{pgldeg}. 
\qed
\begin{thm} \label{properness2}
A moduli space $Q^{\epsilon}_{g,n}(\hat{M}(\dsf), \wsf)$ is a Deligne--Mumford stack of finite type. Moreover, the $\BC_t^*$-fixed locus $Q^{\epsilon}_{g,n}(\hat{M}(\dsf), \wsf)^{\BC_t^*}$ is proper. 
\end{thm}
\textit{Proof.} Exactly the same as in Theorem \ref{properness}. 
\qed
\subsection{Obstruction theory} The obstruction theory $\BT^{\mathrm{vir}}_0$ descends from $\FMsr(\dsf)$ to $\HAz(\dsf)$. For this section, we will denote the descend also by  $\BT^{\mathrm{vir}}_0$.  Let 
\[ \mathsf{f}\colon \CC_{g,n}\rightarrow \HAz(\dsf), \quad \pi_{1}\colon  \CC_{g,n} \rightarrow Q^{\epsilon}_{g,n}(\Mp(\dsf), \wsf)\]
be the universal map and the canonical projection from the universal curve.  Let 
\[(\mathbb{A}, \Phi):=\mathsf{f}^*(\Asf,\phi) - X\times \CC_{g,n},\quad \pi_2\colon S\times \CC_{g,N}  \rightarrow Q^{\epsilon}_{g,n}(\Mp(\dsf), \wsf),\]
be the universal Higgs azalgebra associated to $Q^{\epsilon}_{g,n}(\Mp(\dsf), \wsf)$ and the canonical projection, respectively.

Let $\BE_0^{\mathrm{vir}}$ be defined by the following diagram of exact triangles
\[ \begin{tikzcd}[row sep=scriptsize, column sep = small]
&\BE_0^{\mathrm{vir}}[-1]  \arrow[d] \arrow[r] &  \BE^{\mathrm{vir}}[-1] \arrow[r,] \arrow[d] & R\pi_{2*}(\omega_{\pi_2}[-1]\oplus\CO) \arrow[d] \\
&R\pi_{2*}(\mathbb{A})_0  \arrow[d," \cdot \Phi - \Phi \cdot"] \arrow[r] & R\pi_{2*}(\mathbb{A}) \arrow[r,"\tr"] \arrow[d," \cdot \Phi - \Phi \cdot"] & R\pi_{2*}\CO  \arrow[d,"0"]\\
& R\pi_{2*}(\mathbb{A}_0  \otimes \omega_{\pi_2}) \arrow[r] & R\pi_{2*}(\mathbb{A}  \otimes \omega_{\pi_2}) \arrow[r,"\tr "]& R\pi_{2*}\omega_{\pi_2}
\end{tikzcd} \]
\begin{prop}\label{comparison2}
The complex $\pi_{1*}\mathsf{f}^*\BT_0^{\mathrm{vir}}$ is canonically isomorphic to the complex $\BE_0^{\mathrm{vir}}$. 
\end{prop}
\textit{Proof.} Similar to Proposition \ref{comparison}. \qed 
\begin{prop}\label{perf22}
The complex $R\pi_{1*}L\mathsf{f}^*\BT_0^{\mathrm{vir}}$ is of amplitude $[0,1]$. 
\end{prop}
\textit{Proof.} Similar to Proposition \ref{perf2}.  \qed 

\begin{prop}\label{obst2}
The complex $R\pi_{1*}L\mathsf{f}^*\BT_0^{\mathrm{vir}}$ defines a relative perfect obstruction theory on $Q^{\epsilon}_{g,n}(\Mp(\dsf), \wsf,\asf)$,
\[\psi\colon (R\pi_{1*}L\mathsf{f}^*\BT_0^{\mathrm{vir}})^{\vee} \rightarrow \BL_{Q^{\epsilon}_{g,n}(\Mp(\dsf), \wsf)/\FM_{g,n}}.\]
\end{prop}
\textit{Proof.}
Similar to Proposition \ref{obst}. \qed 
\subsection{Invariants} 
Since $\Mp(\dsf)$ is an orbifold, the evaluation maps associated to the moduli space $Q^{\epsilon}_{g,n}(\Mp(\dsf), \wsf)$ take values in the inertia stack of $\Mp(\dsf)$, 
\[\ev_i \colon Q^{\epsilon}_{g,n}(\Mp(\dsf), \wsf) \rightarrow \CI\Mp(\dsf), \quad i=1 ,\dots n.\]
We also take the projection to a moduli space of stable curves (non-twisted curves)
\[\pi \colon Q^{\epsilon}_{g,n}(\Mp(\dsf), \wsf) \rightarrow \Mbar_{g,n}. \]
Finally, by $\psi$-classes we mean the  \textit{coarse} $\psi$-classes (see \cite[Section 1.7]{YZ}).
\begin{defn} The \textit{descendent} $\epsilon$-\textit{invariants with tautological insertions} are 
\begin{multline*}
	\langle \alpha; \psi_1^{m_{1}}\gamma_{1}, \dots, \psi_n^{m_{n}}\gamma_{n} \rangle^{\vee, \epsilon }_{g,n,\wsf}:=\\
	\int_{[Q^{\epsilon}_{g,n}(\Mp(\dsf),\wsf)^{\BC_t^*}]^{\mathrm{vir}}}\frac{\pi^* \alpha \cdot \prod^{i=n}_{i=1}\psi_{i}^{m_{i}}\ev^{*}_{i}(\gamma_{i},)}{e_{\BC_t^*}(N^{\mathrm{vir}})}\in \BQ[t^{\pm}],
\end{multline*} 
where $\gamma_{1}, \dots, \gamma_{n} \in H_{\mathrm{orb}, \BC_t^*}^{*}(\Mp(\dsf), \BQ):=H_{\BC_t^*}^*(\CI\Mp(\dsf),\BQ)$ are some cohomological classes, $\alpha \in R^*(\Mbar_{g,n}, \BQ)$ is a tautological class and $m_{1}, \dots m_{n}$ are non-negative integers. 
\end{defn}
\subsection{Genus 1 quasimap invariants}  In the case of $\Mp(\dsf)$, we have a story completely parallel  to the one of $\Ms(\dsf)$ discussed in Section \ref{genus1sl}.

We define moduli spaces $Q^{\epsilon}_{E}(\Mp(\dsf),\wsf)$ and $Q^{\epsilon}_{E}(\Mp(\dsf),\wsf)^{\bullet}$ in the same way as in Definition \ref{fixedC} and Definition \ref{quotE}, which are related as follows 
 \begin{equation} \label{Eq2}
	[Q^{\epsilon}_{E}(\Mp(\dsf),\wsf)/E] \cong Q^{\epsilon}_{E}(\Mp(\dsf),\wsf)^{\bullet}.
\end{equation}
In this set-up, we also incorporate the extended degree defined in Definition \ref{extendeddegree}.
\begin{defn} Let $Q^{0^+}_{E}(\Mp(\dsf),\wsf,\asf)$ be the moduli space of $0^+$-stable quasimaps of extended degree $(\wsf,\asf) \in \BZ\oplus\BZ_\rsf$.
	\end{defn}
The extended degree leads to a decomposition of $Q^{0^+}_{C}(\Mp(\dsf),\wsf,\asf)$, 
\[Q^{0^+}_{E}(\Mp(\dsf),\wsf)=\coprod_{\asf} Q^{0^+}_{E}(\Mp(\dsf),\wsf,\asf), \]
which also induces a decomposition of $Q^{0^+}_{E}(\Mp(\dsf),\wsf)^{\bullet}$. 

Similar to Proposition \ref{cosectionelliptic}, we also have surjective cosections in this setting. 
\begin{prop} \label{cosectionelliptic2}
	If  $\wsf\neq 0$, then the obstruction theory of $Q^{\infty}_{g,n}(\Mp(\dsf),\wsf)$ has a $\BC^*_t$-equivariant surjective cosection
	\[h^1(\pi_{1*}\mathsf{f}^*\BT_0^{\mathrm{vir}}) \rightarrow \mathbf t\CO. \]
	So does the obstruction theory of $Q^{0^+}_{E}(\Mp(\dsf),\wsf)$, 
	\[h^1(\pi_{1*}\mathsf{f}^*\BT_0^{\mathrm{vir}}) \rightarrow \mathbf t \CO.\]
\end{prop}
\textit{Proof.}  The moduli space $\Mp(\dsf)$ is also holomorphic symplectic.  If $\epsilon=\infty$, the proof is the same as in Proposition \ref{cosectionelliptic}. 

For $\epsilon=0^+$, we have to notice that in fact the whole stack $\FMp(\dsf)$ is symplectic and therefore $\BT_0^{\mathrm{vir}}\cong (\BT_0^{\mathrm{vir}})^{\vee}$. Moreover, since $E$ is an elliptic curve, the obstruction theory of $Q^{0^+}_{E}(\Mp(\dsf),\wsf)$ is symmetric. The rest then follows from the same arguments as in Proposition \ref{cosectionelliptic}. \qed

\begin{defn} \label{qmsell2} If $\wsf\neq 0$, we define
	\[\QMp^{\asf, \bullet}_{\dsf, \wsf}=  \int_{[Q^{0^+}_{E}(\Mp(\dsf), \asf ,\wsf)^{\bullet, \BC_t^*}]^{\mathrm{vir}}}\frac{1}{te_{\BC_t^*}(N^{\mathrm{vir}})}\in \BQ.\]
	If $\wsf=0$, we define 
	\[\QMp^\asf_{\dsf,0}= \int_{[Q^{0^+}_{E}(\Mp(\dsf), \asf,0)^{ \BC_t^*}]^{\mathrm{vir}}}\frac{1}{te_{\BC_t^*}(N^{\mathrm{vir}})} \in \BQ.\]

	Similarly, we define genus 1 Gromov--Witten invariants for a standard degree $\wsf$,
	\[\GWp^{\bullet}_{\dsf, \wsf}=  \int_{[Q^{\infty}_{E}(\Mp(\dsf),\wsf)^{\bullet, \BC_t^*}]^{\mathrm{vir}}}\frac{1}{te_{\BC_t^*}(N^{\mathrm{vir}})}\in \BQ,\]
	and 
	\[\QMp^{\bullet}_{\dsf, \wsf}:=\sum_\asf \QMp^{\asf, \bullet}_{\dsf, \wsf}, \quad \QMp_{\dsf, 0}:=\sum_\asf \QMp^\asf_{\dsf,0}.\]
\end{defn}

\section{Vafa--Witten theory} \label{sectionvw}
\subsection{Moduli spaces of  Higgs sheaves on $X\times E$}
Most of definitions and results in this section apply to an arbitrary smooth curve $C$. We, however, restrict to the case of an elliptic curve $E$ to simplify the exposition.  For a detailed treatment of the Vafa--Witten theory, we refer to \cite{TT,TT2}.

\begin{defn}
	A Higgs sheaf on $X\times E$ is a pair $(F, \phi)$, where $F$ is a sheaf on $X\times E$ and $\phi \in \Hom(F, F\otimes p^*_X\omega_X)$.  
	\end{defn}
\begin{rmk}By the spectral curve construction,  Higgs sheaves on $X\times E$ are in one-to-one correspondence with compactly supported 2-dimensional sheaves on $\rT^*X\times E$, denoted by $\CF$. 
\end{rmk}
\begin{defn} Let $\CO_{X\times E}(1)$ be an ample line bundle on $X\times E$. A Higgs sheaf $(F, \phi)$  on $X\times E$ is Gieseker semistable (or just semistable for short), if $F$ is torsion-free and  for all non-zero $\phi$-invariant subsheaves $F' \subset F$ we have 
	\[p(F')(t)\leq p(F)(t),\]
where $p(\dots)(t)$ is the reduced Hilbert polynomial associated to $\CO_{X\times E}(1)$. It is $\mu$-semistable, if
\[ \mu(F')\leq \mu(F),\]
where $\mu(\dots)$ is the slope.  It is  $(\mu$-)stable, if the inequalities are strict. 
\begin{prop} \label{stablefiber} There exists $k_0 \in \BN$, such that the set of (Giesker/$\mu$) semistable Higgs sheaves is the same for all $\CO_X(1) \boxtimes \CO_E(k)$ such that $k\geq k_0$. 
	
\end{prop}

\textit{Proof.}  It is enough to prove the claim for $\mu$-stability, because if a semistable Higgs sheaf becomes unstable under a change of  stability, then it actually becomes $\mu$-unstable. This is because the last coefficient of a Hilbert polynomial does not depend on an ample line bundle. Hence if stability of a sheaf changes, then it changes on the level of slopes. 
\\

 If $\rsf=2$ and $\dsf \neq 0$, one can use the argument presented in \cite[Theorem 5.3.2]{HL} for sheaves. The same argument works for  Higgs sheaves.

If $\rsf>2$ or $\dsf=0$, one needs to use Yoshioka's arguments from \cite{Yo}. For that we need the Bogomolov--Giesker inequality for semistable Higgs sheaves, which was established in \cite[Theorem 7]{La}.
We now sketch the argument.  By \cite[Lemma 2.3]{Yo} and \cite[Remark 2.1]{Yo}, if there is a Higgs sheaf $(F,\phi)$ which is $\mu$-semistable with respect to $\CO_X(1) \boxtimes \CO_E(k)$, but not with respect to $\CO_X(1) \boxtimes \CO_E(k+1)$, then there exists a Higgs subsheaf $(F',\phi') \subset (F,\phi)$, such that 
\[ \Delta(F')\geq 0, \quad \Delta(F/F')\geq0, \quad \mu(F')-\mu(F)=0\]
where the slope is taken with respect to $\CO_X(1) \boxtimes \CO_E(k')$, such that $k< k'\leq k+1$. By \cite[Lemma 1.1]{Yo2}, the last equality implies that 
\[(\mathrm c_1(F')/\rk(F')-\mathrm c_1(F/F')/\rk(F/F'))^2\leq -2k-2.\]
By \cite[Lemma 2.1]{Yo}, we obtain that 
\[\Delta(F) \geq (k+1)\rk(F')\rk(F/F')/\rk(F)^2 \geq (k+1)/\rk(F)^2.\] 
It follows that if $k\geq \Delta(F) \rk(F)^2$, the above inequality would lead to a contradiction. In particular, if $k\geq \Delta(F) \rk(F)^2$, the set of $\mu$-semistable Higgs sheaves is the same for different  $\CO_X(1) \boxtimes \CO_E(k)$. 
\qed
\\

By Proposition \ref{stablefiber}, the following definitions are independent of a stability $\CO_X(1) \boxtimes \CO_E(k)$, as long as $k \gg 0$. Taking stabilities $\CO_X(1) \boxtimes \CO_E(k)$ for $k\gg 0$ corresponds to taking the limit $X \rightarrow 0$ in the dimensional reduction, as it is done in \cite{KW}. 
	\end{defn}
\begin{defn} \label{defnu}
	Assume $\dsf\neq 0$. We define $\Ms_\usf(X\times E,\dsf,\wsf)$ to be the moduli space of rank $\rsf$ Higgs sheaves $(F, \phi)$ on $X\times E$,  which are subject to the following  conditions:
\begin{itemize}
	\item[(a)] $\ch(F)_{\mathrm{f}}=\vsf$, $\ch(F)_{\mathrm{m}}=0$ and $0 \leq \ch_1(F)_{\mathrm{d}}<
	\rsf$; 
	\item[(b)] 	$\Delta(F)=2\wsf$; 
	\item[(c)]  $\det(p_{E*}(p_{S}^{*}\usf \otimes F))\cong L'$ for some line bundle $L' \in \Pic(E)$;
	\item[(d)]$\tr(\phi)=0$;
	\item[(e)] semistable with respect to $\CO_X(1) \boxtimes \CO_E(k)$ for some $k\gg0$.
\end{itemize}
\end{defn}
The moduli space $\Ms_\usf(X\times E,\dsf,\wsf)$ has its advantages, e.g.\ it has a direction relation to the moduli space of quasimaps, as it is shown in Proposition \ref{ident}. However, if $\dsf=0$, another space is more appropriate. 
\begin{defn} \label{defnL}  Let $0 \leq \asf < \rsf $. 
We define
$\Ms_{\mathsf L}(X\times E,\dsf, \asf, \wsf)$ to be the moduli space of rank $\rsf$ Higgs sheaves on $X\times E$, which are subject to the following conditions
\begin{itemize}
	\item[(a)] $\ch(F)_{\mathrm{f}}=\vsf$, $\ch(F)_{\mathrm{m}}=0$ and $ \ch_1(F)_{\mathrm{d}}=\asf $;
	\item[(b)]$\Delta(F)=2\wsf$; 
	\item[(c)] $\det(F)\cong \mathsf L$ for some line bundle $\Lsf \in \Pic(X \times E)$;
	\item[(d)]$\tr(\phi)=0$;
	\item[(e)] semistable with respect to $\CO_X(1) \boxtimes \CO_E(k)$ for some $k\gg0$.
\end{itemize}
	\end{defn}
\begin{rmk} The reason we have to specify $\asf$ for the moduli space $\Ms_{\mathsf L}(X\times E,\dsf, \asf, \wsf)$  is that we do not make an assumption on $\dsf$. If $\dsf\neq 0$, then by 
	\[ \Delta(F)/2= \dsf \cdot \ch_1(F)_{\mathrm d}-\rsf \cdot \ch_2(F)_{\mathrm d}, \]
	  the class $\ch_1(F)_{\mathrm{d}}$ is determined by $\wsf$, $\dsf$ and the condition $0 \leq \ch_1(F)_{\mathrm{d}}<\rsf$. See also Lemma \ref{ul}. 
	\end{rmk}
\begin{rmk} Note that we introduce the condition $0 \leq \asf< \rsf$ to get rid of redundancy associated to the operation of tensoring moduli spaces with line bundles. We did the same for $\Ms(\dsf)$ in Section \ref{SLbundles}.
	
	Also note that the moduli spaces are isomorphic for different choices of $\usf$, $L'$ and $\Lsf$. However, we include $\usf$ and $\Lsf$ in the notation to distinguish between two different determinant conditions that define these moduli spaces.
	
	The reason that we make the assumption on $\dsf$ for the definition of $\Ms_\usf(X\times E,\dsf,\wsf)$ is that $\usf \in K_{0}(X)$ was defined to be an element such that $\chi(\vsf \cdot \usf)=1$, which exists only if $\dsf \neq 0$. 
	\end{rmk}

\begin{lemma} \label{stability}
	 If $\dsf\neq 0$, then there exists $k_0 \in \BN$, such that the following holds. A Higgs sheaf $(F,\phi)$ is $\mu$-semistable with respect to $\CO_X(1) \boxtimes \CO_E(k)$ for some $k\geq  k_0$, then a general fiber of $(F,\phi)$ is semistable. Moreover, if a general fiber of $(F,\phi)$ is stable, then $(F,\phi)$ is $\mu$-stable with respect to $\CO_X(1) \boxtimes \CO_E(k)$ for all $k\geq k_0$. 
	\end{lemma}
\textit{Proof.} The claim follows from Proposition \ref{stablefiber} by arguments similar to those presented in \cite[Lemma 3.17]{N}. \qed
 \begin{prop} \label{ident} By associating to a quasimap $f$ the corresponding family $(F,\phi)$ twisted by an appropriate  line bundle $L'\in \Pic(E)$, we obtain an identification
 	\[Q^{0^+}_E(\Ms(\dsf),\wsf) \cong \Ms_{\usf}(X\times E,\dsf, \wsf), \quad f\mapsto (F\otimes p^*_EL',\phi).\]
 	\end{prop}
 \textit{Proof.} 
  Let us check the defining conditions of the moduli space $\Ms_{\usf}(X\times E,\dsf, \wsf)$. Firstly, by the discussion in Section \ref{slpre}, 
 \begin{align} \label{topinv}
 	\begin{split}
 	&\ch(F)_{\mathrm{f}}=\dsf \\
 	&\ch(F)_{\mathrm{m}}=0 \\
 	&\Delta(F)=2\wsf.
 	\end{split}
 	\end{align}
 Moreover, there exists a line bundle $L' \in \Pic(E)$,  such that for all sheaves $F$ associated to quasimaps in $Q^{0^+}_E(\Ms(\dsf),\wsf)$ the following holds
 \[0 \leq \ch_1(F\otimes p^*_CL')_{\mathrm{d}}< \rsf.\]
By Lemma \ref{degreechern}, its degree $\deg(L')$ is uniquely determined by $\dsf$ and $\wsf$. Moreover, by (\ref{detcond}),
 \[ \det(p_{E*}(p_{S}^{*}\usf \otimes F))\cong L'.\]
 Note that tensoring by $L'$ does not change any topological invariants in (\ref{topinv}). The vanishing of the trace $\tr(\phi)=0$ follows from the fact that it is zero at each fiber.

 Let us now deal with the stability. By Lemma  \ref{stability}, the  Higgs sheaf $(F,\phi)$ associated to a quasimap is $\mu$-stable with respect to a suitable polarisation. Conversely, any $\mu$-semistable Higgs sheaf $(F,\phi )$ (with the Chern character prescribed by Definition \ref{defnu} ) is semistable at a general fiber also by Lemma \ref{stability}, which in turn implies that it is stable at a general fiber, because $\dsf \neq0$. Therefore the association of a  Higgs sheaf to a quasimap is in fact surjective for all $B$-valued points. It is injective by the definition of the moduli stack $\FMs(\dsf)$. 
\qed
\\

We will now relate the moduli spaces $\Ms_{\usf}(X\times E,\dsf, \wsf)$ and $\Ms_{\Lsf}(X\times E,\dsf, \wsf)$. Let $\Jac(C)[\rsf] \subset \Jac(C)$ be the subgroup of $\rsf$-torsion points of the Jacobian $\Jac(C)$. The group $\Jac(C)[\rsf]$ acts on $\Ms_{\Lsf}(X\times E,\dsf, \wsf)$ by tensoring a sheaf with a $\rsf$-torsion line bundle.
\begin{lemma} \label{ul} Assume $\dsf\neq 0$.  Let $0\leq\asf<\rsf$ be defined by the following equation 
	\[\wsf=\dsf\cdot \asf 
	\mod \rsf. \]
	The natural map
	\[\Ms_{\Lsf}(X\times E,\dsf, \asf, \wsf) \rightarrow \Ms_{\usf}(X\times E,\dsf, \wsf),\]
	\[\quad (F,\phi) \mapsto (F\otimes  p^*_C\det(p_{C*}(p_{S}^{*}\usf \otimes F))^{-1},\phi)\]
	is a $E[\rsf]$-torsor with respect to the  action of $E[\rsf]$. 
	\end{lemma}
  \textit{Proof.} Since
  \[ \Delta(F)/2= \dsf \cdot \ch_1(F)_{\mathrm d}-\rsf \cdot \ch_2(F)_{\mathrm d}, \]
   $\asf=\ch_1(F)_{\mathrm d}$ can be obtained from $\wsf$ and $\dsf$ via the specified equation. For the rest of the proof we refer to  \cite[Section 3.5]{N}.  \qed 
\subsection{Moduli spaces of  Higgs azalgberas on $X\times E$}  \label{azal}
\begin{defn}
	A Higgs azalgebra on $X\times E$ is a pair $(\CA, \phi)$, where $\CA$ is an azalgbera on $X\times E$ and $\phi \in H^0(\CA\otimes  p^*_X\omega_X, X\times E)$.  
\end{defn}
Consider a Higgs azalgebra $(\CA,\phi)$ on $X\times E$. Recall the gerbe $p\colon \CX\mathcal \times \CC \rightarrow X\times E$ constructed in Definition \ref{brauergerbe}. Note that the notation is a bit unfortunate, as $\CX\mathcal \times \CE$ is not necessarily a product of gerbes from $X$ and $E$.  By construction of $\CX \times \CE$, there exists a sheaf $\overline{F}$ (do not confuse it with $\tilde{F}$) and $\overline{\phi} \in \Hom(\overline{F}, \overline{F}\otimes p^*\omega_S)$ on $\CX\times \CE$, such that $(R \CE nd(\overline F, \overline F), \overline \phi)$ descends to $(\CA,\phi)$. 

	\begin{defn} \label{stabaz}  Let $\CO_{X\times E}(1)$ be an ample line bundle on $X \times E$. A Higgs azalgebra $(\CA, \phi)$ on $X\times E$ is (Gieseker) semistable, if $\CA$ is torsion-free and  for all non-zero $\phi$-invariant subsheaves $F' \subset \overline{F}$ we have 
		\[p^g_{\CO_{X\times E}(1)}(F')(t)\leq p^g_{\CO_{X\times E}(1)}(\overline{F})(t),\]
		where $p^g_{\CO_{X\times E}(1)}(\dots)(t)$ is the reduced  geometric Hilbert polynomial taken with respect to the pullback  of $\CO_{X\times E}(1)$ to $\CX \times \CE$ (see \cite[Definition 2.2.7.2]{Lie3} for the definition). It is $\mu$-semistable, if
		\[ \mu(F')\leq \mu(\overline{F}).\]
		It is  $(\mu$-)stable, if the inequalities are strict. 
		
	\end{defn}
 \begin{prop} \label{nontrivialbrauer}
	A Higgs azalgebra $(\CA,\phi)$ on $X\times E$ with a non-trivial Brauer class  is stable  with respect to an ample line bundle $\CO_{X\times E}(1)$, if and only if it is torsion-free. In particular, its stability is independent of an ample line bundle $\CO_{X\times E}(1)$. 
\end{prop}	
\textit{Proof.} See \cite[Theorem 2.10]{JK}. \qed 

\begin{rmk} Proposition \ref{nontrivialbrauer} also shows that different constructions of moduli spaces of twisted Higgs sheaves and Higgs azalgebras agree (under the assumption that $\rsf$ is prime). For example, in \cite{JK}, moduli spaces of Higgs azalgebras were constructed via sheaves on projective bundles, while in \cite{Jiang} via twisted sheaves. Since stability in both cases depends just on torsion subsheaves and the associated categories are naturally isomorphic, the moduli spaces are the same. 
	\end{rmk}
By Proposition \ref{nontrivialbrauer} and Proposition \ref{stablefiber}, the following definitions are independent of a stability $\CO_X(1) \boxtimes \CO_E(k)$, as long as $k\gg0$. 
\begin{defn} We define
	$\Mp_{\msf}(X\times E,\dsf, \asf, \wsf)$ to be the moduli space of Higgs azalgebras on $X\times E$, which are subject to the following conditions
	\begin{itemize}
		\item[(a)] $\delta(\CA)_{\mathrm{f}}=\dsf$, $\delta(\CA)_{\mathrm{m}}=\msf$ and $\delta(\CA)_{\mathrm{d}}=\asf$;
		\item[(b)]$\mathrm{c}_2(\CA)=2\wsf$; 
		\item[(c)]$\tr(\phi)=0$;
		\item[(d)] semistable with respect to $\CO_X(1) \boxtimes \CO_E(k)$ for some $k\gg0$.
	\end{itemize}
	\end{defn}
  	
  	\begin{lemma} \label{torsionfree} The azalgbera $(\CA,\phi)$ associated to a quasimap $f\colon C \rightarrow \HAz(\dsf)$ is torsion-free. 
  	\end{lemma} 
  	\textit{Proof.} Consider the family $(\tilde F,\tilde \phi)$ associated to the quasimap $\tilde{f}$.  By the definition of the quasimap, the generic fiber of  $(\tilde F,\tilde \phi)$ is stable and therefore torsion-free. Hence by \cite[Lemma 3.18]{N}, the sheaf $\tilde F$ is torsion-free itself. Since $h^*\CA=R \CE nd(\tilde F,\tilde F)$, the azalgebra $\CA$ must be torsion-free too.
  	\qed
  	\\

  	
  
\begin{lemma}   \label{stability2}
 	 If $\dsf\neq 0$, then there exists $k_0 \in \BN$, such that the following holds. A Higgs azalgebra $(\CA,\phi)$ is $\mu$-semistable with respect to $\CO_X(1) \boxtimes \CO_E(k)$ for some $k\geq  k_0$, then a general fiber of $(\CA,\phi)$ is semistable. Moreover, if a general fiber of $(\CA,\phi)$ is stable, then $(\CA,\phi)$ is $\mu$-stable with respect  $\CO_X(1) \boxtimes \CO_E(k)$ for all $k\geq k_0$.  
 \end{lemma} 
  
\textit{Proof.}   The proof split into parts, depending on whether $\beta$ is algebraic or not (i.e.\ $\beta$ is of the image of algebraic class under the natural map $H^2(X\times E, \BZ) \rightarrow H^2(X\times E, \BZ_\rsf)$). 

Assume firstly that $\beta$ is algebraic, then by the exactness of the sequence 
\[H^1(X\times E, \BC^*) \rightarrow H^2(X\times E, \BZ_\rsf) \rightarrow H^2(X\times E, \BC^*)\]
we obtain that the Brauer class of $\CA$ is zero and therefore the gerbe $\CX \times \CC$ is trivial, which means that $\overline{F}$ descends to sheaf on $X\times E$ and the result then follows from Theorem \ref{stability}. 

Assume now that $\beta$ is not algebraic. By exactness of the sequence above, this is equivalent to the Brauer class being non-zero. Then we can use the argument from \cite[Theorem 2.10]{JK} concerning torsion-free sheaves, which relies on the result from \cite[Lemma 3.2]{Yo3}. More precisely, since the Brauer class is non-zero, there is no rank 1 sheaf on $\CX \times \CE$ (otherwise its double dual would be a line bundle and it will trivialise the gerbe $\CX \times \CE$, which contradicts the fact that Brauer class is non-zero). Sine $\rsf$ is prime and since there exists a sheaf of rank $\rsf$ on $\CX \times \CE$ (namely, the sheaf $\overline{F}$ itself), it must be that the minimal rank of a $\CX \times \CE$-twisted sheaf is $\rsf$ by \cite[Lemma 3.2]{Yo3}. 

Assume for the moment that $\overline{F}$ is locally free. Assume the general fiber of $(\CA,\phi)$ is unstable, then over an open subset $U \subset C$, there exists a destabilising family of $\phi$-invariant left ideals, $I_{|U} \subset \CA_{|U}$. We can extend it to a subsheaf $I' \subset \CA$ and then by taking the ideal generated by $I'$, we obtain an ideal $l\subset \CA$, such that the rank of $I$ is strictly smaller than the rank of $\CA$.  And ideal $I$ is a descend of $\CH om(F', \overline{F})$ for some subsheaf $F' \subset \overline F$. Indeed, locally  an ideal of $\CE nd(H,H)$ for a locally free sheaf $H$ is of the form $\CH om(H', H)$ for some subsheaf $H' \subset H$. Globally, such subsheaves glue to a $\CX \times \CC$-twisted subsheaf of  $\overline F$. We thereby produced a sheaf of rank strictly smaller $\rsf$, contradicting the minimality of $\rsf$. 

Assume now that $\overline{F}$ is not locally free. Then we can take the double dual of $\overline{F}$, which is locally free. The associated azalgebra is still unstable at the general fiber. We thereby reduced the situation to the locally-free case. 

The other direction follows from Lemma \ref{torsionfree}. 
\qed

  	 \begin{prop} \label{ident2} 
   By associating to a quasimap $f\colon E\rightarrow \HAz(\dsf) $ the corresponding  Higgs azalgebra $(\CA, \phi)$, we obtain the identification
  	\[Q^{0^+}_E(\hat{M}(\dsf), \asf, \wsf) \cong \coprod_{\msf} \Mp_{\msf}(X\times E,\dsf, \asf, \wsf), \quad f\mapsto (\CA,\phi).\]
  \end{prop}
\textit{Proof.} With Lemma \ref{topazal} and Lemma  \ref{stability2} the proof is the same as for Proposition \ref{ident2}. 
\qed
\\

There is also a relation between moduli spaces $\Mp_{\msf}(X\times E,\dsf, \asf, \wsf)$ and $\Ms_{\Lsf}(X\times E,\dsf, \asf, \wsf)$. Let $\Jac(X\times E)[\rsf]$ be a group of $\rsf$-torsion line bundles on $X\times E$, it acts on $\Ms_\Lsf(X\times E,\dsf, \asf, \wsf)$ by tensoring by line bundles.

\begin{lemma} \label{slpgl}  The natural map 
	\begin{align*} \Ms_{\Lsf}(X\times E,\dsf, \asf, \wsf) &\rightarrow  \Mp_{0}(X\times E,\dsf, \asf, \wsf) \\
(F,\phi) &\mapsto (R\CE nd(F,F),\phi)
\end{align*}
is a $\Jac(X\times E)[\rsf]$-gerbe with respect to the natural action of $\Jac(X\times E)[\rsf]$ on $\Ms_{\Lsf}(X\times E,\dsf, \asf, \wsf)$.
	\end{lemma}
\textit{Proof.} It is not difficult to see that the topological data matches on both sides. Moreover, since $\msf=0$, the Brauer class of azalgberas in  $\Ms(X\times E,\dsf, \asf, \wsf)$ is trivial and therefore all azalgebras are of the form $R\CE nd(F,F)$ for some sheaves $F$ on $X\times E$. This implies the map is surjective and the stabilities for Higgs sheaves and Higgs azalgebras are equivalent. Finally the fact that the map is indeed a gerbe follows from Lieblich's derived Skolem--Noether theorem \cite[Corollary 5.1.6]{Lie}. \qed
\subsection{Genus 1 Vafa--Witten invariants}
Throughout the section we will use  the moduli spaces $\Ms_{\usf}(X\times E,\dsf, \wsf)$ and $\Ms_{\Lsf}(X\times E,\dsf, \asf, \wsf )$  interchangeably. The former has a natural $E$-action, while the latter is defined for $\dsf = 0$. If $\dsf \neq 0$, the invariants that these moduli spaces produce are equivalent by Lemma \ref{ul}. The only case that requires a special treatment is when 
\[(\dsf, \msf, \asf) =(0,0,0),\]
because only in this case moduli spaces contain strictly semistable Higgs sheaves/Higgs azalgberas. Let us show that if $(\dsf, \msf, \asf) \neq(0,0,0)$, then there are no  strictly semistable Higgs sheaves/Higgs azalgberas. Assume firstly we are in the case of Higgs sheaves, then if $(\dsf, \msf, \asf) \neq (0,0,0)$, one can choose $\CO(1)\boxtimes \CO_E(k)$ with $k\gg 0$ in the way that  rank $\rsf$ and degree $\deg(F)$ are coprime, forcing $\mu$-stability and $\mu$-semistability to coincide. By Proposition \ref{stablefiber}, all other stabilities $\CO(1)\boxtimes \CO_E(k)$  with $k\gg 0$ have the same property. By Proposition \ref{nontrivialbrauer}, the same holds for Higgs azalgberas.  We will treat the case of $(\dsf, \msf, \asf) =(0,0,0)$ in Section \ref{noncoprime}. 
\\

There exists a $\BC_t^*$-action on  $\Ms_{\usf}(X\times E,\dsf, \wsf)$, $\Ms_{\Lsf}(X\times E,\dsf, \asf, \wsf)$ and $\Mp_{\msf}(X\times E,\dsf, \asf, \wsf)$ given by scaling the Higgs field $\phi$, such that the fixed locus is proper. This allows us to define invariants. 
\begin{defn}   Assuming that $(\dsf, \asf, \msf)\neq (0,0,0)$, we define Vafa-Witten invariants
	\begin{align*}
	&\VWs^{\asf}_{\dsf, \mathsf w}
:=\int_{[\Ms_{\Lsf}(X\times E,\dsf, \asf, \wsf)^{\BC_t^*}]^{\mathrm{vir}}}\frac{1}{e_{\BC_t^*}(N^{\mathrm{vir}})}\in \BQ \\
&\VWp^{\msf,\asf}_{\dsf, \mathsf w}
:=\int_{[\Mp_{\msf}(X\times E,\dsf, \asf, \wsf)^{\BC_t^*}]^{\mathrm{vir}}}\frac{1}{e_{\BC_t^*}(N^{\mathrm{vir}})}\in \BQ .
\end{align*}
\end{defn}
However, by the existence of surjective cosections constructed in the same way as in Proposition \ref{cosectionelliptic} and \ref{cosectionelliptic2}, these invariants vanish, unless $\wsf=0$. There are two ways to obtain non-trivial invariants in the case of $\wsf \neq0$. The first one is to quotient the moduli space by the action $E$, which corresponds to quasimaps up to translations of the domain curve. The second one is to take an insertion. Let us discuss both of them.

  \subsection{Reduced Vafa--Witten invariants} Assume $\dsf\neq 0$. 
  The elliptic curve $E$ acts on  $\check{M}_{\usf}(X\times E,\dsf, \wsf)$ and $\Mp_{\msf}(X\times E,\dsf, \asf, \wsf)$ by translation. If $\wsf\neq 0$, the action has only finite stabilisers. One can see that by invoking Proposition \ref{ident} and Proposition \ref{ident2}, the case of $\wsf\neq 0$ corresponds to non-constant quasimaps. Non-constant quasimaps cannot be fixed by entire $E$. 

 \begin{defn} Assuming $\dsf\neq 0$ and $\wsf\neq 0$, we define 
 \begin{align*}	
 &	\Ms_{\usf}(X\times E,\dsf,\wsf)^{\bullet }:= [\Ms_{\usf}(X\times E,\dsf,\wsf)/E]\\
 &	\Mp_{\msf}(X\times E,\dsf, \asf, \wsf)^{\bullet }:= [\Mp_{\msf}(X\times E,\dsf, \asf, \wsf)/E].
 	\end{align*}
Let $\asf \in \BZ_\rsf$. We define \textit{reduced} Vafa--Witten invariants 
 \begin{align*}
 			&\VWs^{\asf, \bullet}_{\dsf, \mathsf w}:=
 				\begin{cases}
 					\int_{[\Ms_{\usf}(X\times E,\dsf, \asf, \wsf)^{\bullet, \BC_t^*}]^{\mathrm{vir}}}\frac{1}{te_{\BC_t^*}(N^{\mathrm{vir}})} , & \text{if}\ \wsf=\dsf\cdot \asf \mod \rsf \\
 					0, &\text{otherwise} 
 				\end{cases} 
 			\\
 			&\VWp^{\msf, \asf, \bullet}_{\dsf, \mathsf w}
 			:=\int_{[\Mp_{\msf}(X\times E,\dsf, \asf, \wsf)^{\bullet, \BC_t^*}]^{\mathrm{vir}}}\frac{1}{te_{\BC_t^*}(N^{\mathrm{vir}})} .
 		\end{align*}
 	\begin{rmk} The class $\asf$ is included in the notation of $\VWs^{\asf, \bullet}_{\dsf, \mathsf w}$ for the same reason as in Definition \ref{qmsell}, also see Remark \ref{asfnot}.
 		\end{rmk}
 	
 	\begin{rmk}
 		  One can also define a (twisted) action of $E$ directly on $\Ms_{\Lsf}(X\times E,\dsf, \asf, \wsf)$, however, it is a slightly involved because  translations of $E$ do not fix line bundles of positive degree on $E$, so one has to twist the action by degree 0 line bundles. See \cite[Section 3]{Nqm} for more details. 
 		\end{rmk}
\end{defn}

\subsection{Vafa--Witten invariants with insertions} \label{insertions}
 Let $\mathsf{F}$ be the universal sheaf on  $X\times E\times  \Ms_{\usf}(X\times E,\dsf, \wsf)$. For a class $\beta \in H^*(X\times E,\BQ)$, we define 
\[\mu(\beta):=\pi_{\Ms*}\left( \Delta(\Fsf)/2\rsf\cdot \pi^*_{X\times E}\beta \right ) \in H_{\BC_t^*}^{4-*}(\Ms_{\Lsf}(X\times E,\dsf, \asf, \wsf),\BQ),\]
where 
\[\Delta(\Fsf):=\mathrm c_1^2(\mathsf F)-2\rsf\mathrm{ch}_2(\mathsf F).\]

Let $\mathsf A$ be the universal azalgbera on $X\times E\times \Mp_{\msf}(X\times E,\dsf, \asf, \wsf)$. For a class $\beta \in H^*(X\times E,\BQ)$, we define 
\[\mu(\beta):=\pi_{\Mp*} \left ( \mathrm{c}_2(\mathsf A)/2\rsf \cdot \pi^*_{X\times E}\beta \right ) \in H_{\BC_t^*}^{4-*}(\Mp_{\msf}(X\times E,\dsf, \asf, \wsf),\BQ).\]
Note that definition of insertions for  $\Ms_{\Lsf}(X\times E,\dsf, \asf, \wsf)$ and $ \Mp_{\msf}(X\times E,\dsf, \asf, \wsf)$ are in fact compatible, because 
\[\mathrm c_2(R \CE nd(\Fsf,\Fsf))=\Delta(\Fsf).\]
\begin{defn} Assuming  $(\dsf, \asf, \msf)\neq(0,0,0)$, we define Vafa-Witten invariants with $\mu$-\textit{insertions}
	\begin{align*}
			&\VWs^{\asf}_{\dsf, \mathsf w}(\beta):=\int_{[\Ms_{\Lsf}(X\times E,\dsf, \asf, \wsf)^{\BC_t^*}]^{\mathrm{vir}}}\frac{e^{\mu(\beta)}}{e_{\BC_t^*}(N^{\mathrm{vir}})}\in \BQ[t^{\pm}] \\
	&\VWp^{\msf, \asf}_{\dsf, \mathsf w}(\beta):=\int_{[ \Mp_{\msf}(X\times E,\dsf, \asf, \wsf)^{\BC_t^*}]^{\mathrm{vir}}}\frac{e^{\mu(\beta)}}{e_{\BC_t^*}(N^{\mathrm{vir}})}\in \BQ[t^{\pm}].
	\end{align*}
	\end{defn}
\subsection{Rigidification}
Following \cite[Section 3.3]{O} or \cite[Lemma 3.3]{MO1}, one can trade the action of $E$ on  $\Ms_{\usf}(X\times E,\dsf, \wsf)$ or $\Ms_{\msf}(X\times E,\dsf, \asf, \wsf)$ for a primary insertion. In fact the same argument also works for $\mu$-insertions. 
\begin{prop} \label{rigid}
Let 
\[B_{\mathsf w} := \frac{\rsf\mathbb 1 \boxtimes [\mathrm{pt}]}{\mathsf w} \in H^2(X\times E,\BQ).\]
Assuming $\dsf\neq 0$ and $\wsf \neq 0$, we have
\begin{align*}
	&[\VWs^{\asf}_{\dsf, \mathsf w}(B_\wsf )]_t = \rsf^2 \VWs_{\dsf, \mathsf w}^{\asf, \bullet}\\
	&[\VWp^{\msf,\asf}_{\dsf, \mathsf w}(B_{\wsf})]_t = \VWp_{\dsf, \mathsf w}^{\msf, \asf, \bullet},
	\end{align*}
where $[\dots ]_t$ is the coefficient of a Laurent series at $t$. 
	\end{prop}
\textit{Proof.}
To fix a notation, we present the argument for $\Ms_{\Lsf}(X\times E,\dsf, \asf, \wsf)$, the same argument applies to $\Mp_{\msf}(X\times E,\dsf, \asf, \wsf)$.  By Lemma \ref{ul}, 
\begin{multline*}
	\VWs^{\asf}_{\dsf, \mathsf w}(\beta):=\int_{[\Ms_{\Lsf}(X\times E,\dsf, \asf, \wsf)^{\BC_t^*}]^{\mathrm{vir}}}\frac{e^{\mu(\beta)}}{e_{\BC_t^*}(N^{\mathrm{vir}})}  \\
	= \int_{[\Ms_{\usf}(X\times E,\dsf, \wsf)^{\BC_t^*}]^{\mathrm{vir}}}\frac{\rsf^2e^{\mu(\beta)}}{e_{\BC_t^*}(N^{\mathrm{vir}})}.
	\end{multline*}
Hence we can work with invariants associated to $\Ms_{\usf}(X\times E,\dsf, \wsf)$. For simplicity, let  $\Ms:= \Ms_{\usf}(X\times E,\dsf, \asf, \wsf)$ and 
\[[\Ms]^{\mathrm{vir}}:= \frac{[\Ms^{\BC_t^*}]^{\mathrm{vir}}}{e_{\BC_t^*}(N^{\mathrm{vir}})}.\]
Consider now the following fiber diagram, 
\[
\begin{tikzcd}[row sep=scriptsize, column sep = scriptsize]
	X\times E & X\times E\times \Ms  \arrow[l, "\pi_{X\times E}"]\arrow[d,"\rho"] \arrow[r, "\pi_{\Ms}"] &  \Ms \arrow[d, "q"]\\
	   & (X\times E\times \Ms)/ E \arrow[r,"\pi_{\Ms/E}"] & \Ms/E 
\end{tikzcd}
\]
\noindent The obstruction theory of $M$ is pullbacked to the one of $M/E$, hence 
\[q^* [\Ms/E]^{\mathrm{vir}}=[\Ms]^{\mathrm{vir}}.\]
By definition, $[\VWs(B_{\wsf})]_t/\rsf^2$ is equal to the degree of the following class
\[\check{\mathbf{Z}}=\left( \left( \mathrm c_1^2(\mathsf F)/2\rsf-\mathrm{ch}_2(\mathsf F)\right )\cdot \pi^*_{X\times E}B_{\wsf}\right )\cap \pi_{\Ms}^*[\Ms]^{\mathrm{vir}}.\]
 We therefore have to show that 
\[(\pi_{M/E}\circ \rho)_*\check{\mathbf{Z}}=[\Ms/E]^{\mathrm{vir}}.\] 
Let $\bar{\mathsf{F}}$ be a descend of $\mathsf F$ to $(X\times E\times \Ms)/E$, define 
\[\check{\mathbf z}=\frac{\left( \mathrm c_1^2(\bar{\mathsf F})/2-\rsf\mathrm{ch}_2(\bar {\mathsf F})\right )}{\mathsf w  } \cap \pi_{\Ms/E}^*[\Ms/E]^{\mathrm{vir}}.\]
Now let 
\[ \iota \colon X\times 0_E \times \Ms \hookrightarrow X\times E\times \Ms\]
be the inclusion of a slice of the action of $E$, which mapped isomorphically onto $(X\times E\times \Ms) /E$ by $\rho$. Then
\[ \rho_*\check{\mathbf{Z}} = \rho_* (\pi^*_E(\mathrm{[pt]})\cap \rho^*\check{\mathbf{z}})=\rho_*(\iota_* \iota^*\rho ^*\check{\mathbf{z}})=(\rho\circ \iota)_*(\rho \circ \iota)^*\check{\mathbf{z}}=\check{\mathbf{z}}.\]
We conclude 
\[(\pi_{\Ms/E}\circ \rho)_*\check{\mathbf{Z}}=\pi_{\Ms/E*}\check{\mathbf{z}}=\frac{ \mathrm c_1^2(F)/2-\rsf\mathrm{ch}_2( F)}{\mathsf w }\cdot [\Ms/E]^{\mathrm{vir}}=[\Ms/E]^{\mathrm{vir}},\]
where we used that  
\[\mathsf w=  \mathrm c_1^2(F)/2-\rsf \mathrm{ch}_2( F)\]
for a sheaf in $F \in \Ms$, which follows from Lemma \ref{degreechern}.
\qed 
\subsection{Comparison between Vafa--Witten and quasimap invariants} \label{sectioncompar}
The results of the previous sections give us the following comparison statement between Vafa--Witten invariants and quasimap invariants. 
\begin{thm} \label{vwqm} Assume $\dsf\neq 0$. If $\wsf \neq 0$, then
\begin{itemize}
\item[(a)]$\check{\mathsf{QM}}^{\asf, \bullet}_{\dsf, \mathsf w}=	\VWs^{\asf, \bullet}_{\dsf, \mathsf w}=[\check{\mathsf{VW}}_{\dsf, \mathsf w}(B_\wsf ) ]_t/ 	\rsf^2$,\\
\item[(b)]$	\hat{\mathsf{QM}}^{\asf, \bullet}_{\dsf, \mathsf w}=\sum_{\msf} \VWp^{\msf, \asf, \bullet}_{\dsf, \mathsf w}=\sum_\msf [\hat{\mathsf{VW}}^{\msf,\asf}_{\dsf, \mathsf w}(B_\wsf )]_t $,\\
\item[(c)]	$\VWp^{0, \asf, \bullet}_{\dsf, \mathsf w}=	\VWs^{\asf, \bullet}_{\dsf, \mathsf w} /\rsf^{2g+2}$.		
\end{itemize}
If $\wsf=0$, then 
	\begin{itemize}
	\item[(a)]$\check{\mathsf{QM}}^{\asf}_{\dsf, \mathsf 0}=	\check{\mathsf{VW}}^{\asf}_{\dsf, \mathsf 0} / 	\rsf^2$\\
	\item[(b)]$\hat{\mathsf{QM}}^{\asf}_{\dsf, 0}=\sum_{\msf}\hat{\mathsf{VW}}^{\msf, \asf}_{\dsf, 0}$\\
	\item[(c)]$\VWp^{0,\asf}_{\dsf, \mathsf w}=	\VWs^{\asf}_{\dsf, \mathsf w} /\rsf^{2g+2}$.		
	\end{itemize}
\end{thm}
\textit{Proof.}
If $\wsf=0$, then (a) follows from Proposition \ref{ident} and Lemma \ref{ul}, (b) follows from Proposition \ref{ident2} and (c) from Lemma \ref{slpgl}.

If $\wsf\neq0$, it follows from (\ref{Eq}), (\ref{Eq2}), Proposition \ref{ident} and Proposition \ref{ident2}, that we have the following identifications
\begin{align*}
	&Q^{0^+}_E(\hat{M}(\dsf),\wsf)^{\bullet} \cong [\Ms_{\usf}(X\times E,\dsf,\wsf)/E]\\
	&Q^{0^+}_E(\hat{M}(\dsf),\asf,\wsf)^{\bullet}  \cong \coprod_{\msf} [\Mp_{\msf}(X\times E,\dsf, \asf, \wsf)/E],
\end{align*}
which imply the first equalities in (a) an (b). The second equalities follow from Proposition \ref{rigid}. Finally, (c) follows from (a), (b) and Lemma \ref{slpgl}.
\qed
\section{Quasimaps, $\dsf=0$} \label{noncoprime}
\subsection{Semistable Vafa--Witten invariants} We now deal with the case when $(\dsf,\msf,\asf)=(0,0,0)$.

Let $\CO(1):= \CO_X(1) \boxtimes \CO_E(k)$ for some $k\gg0$. Choose $m\gg0$  such that 
\[H^i(F(m))=0,\quad i>0, \]
 for all $(F, \phi)\in \Ms_{\CO}(X\times E, 0, 0, \wsf)$.
\begin{defn}
 A Joyce--Song pair $(F, \phi, s)$ consists of
 \begin{itemize}
 	\item a Higgs sheaf $(F,\phi) \in \Ms_\CO(X\times E, 0, 0, \wsf)$,
 	\item a nonzero section $s \in H^0(F(m))$. 
 	\end{itemize}
 A pair $(F, \phi, s)$ is \textit{stable}, if and only if 
 \begin{itemize}
 	\item if $\CF' \subset \CF$ is a proper $\phi$-invariant subsheaf which destabilises $F$,  then $s$ does not factor through $F'(m) \subset F(m)$. 
 	\end{itemize}
 Let $\Ms^{JS}(X\times E, 0, 0, \wsf)$ be the moduli space of stable Joyce--Song pairs associated to a moduli space $\Ms_\CO(X\times E, 0, 0, \wsf)$.
	\end{defn}

By \cite{JS} and \cite{TT}, the moduli space $\Ms^{JS}(X\times E, 0, 0, \wsf)$ carries a symmetric perfect obstruction theory. Moreover, there exits a $\BC_t^*$-action defined by scaling the Higgs field, such that the fixed locus is proper. There is no universal family over $\Ms^{JS}(X\times E, 0, 0, \wsf)$. Luckily, $\Delta(\mathsf F)$ can still be defined as $\mathrm c_2(R \CE nd(\mathsf F,\mathsf F)$), where $R \CE nd(\mathsf F,\mathsf F)$ is the descend from the moduli stack over which the universal family exists. The $\mu$-map is then constructed as in Section \ref{insertions}.  We proceed to define Vafa--Witten invariants with insertions via $\Ms^{JS}(X\times E, 0, 0, \wsf)$ as follows. 
\begin{defn} Let
	\begin{align*}
		&\VWs_{ 0,\mathsf w}^{JS,0}:=\int_{[\Ms^{JS,\BC_t^*}]^{\mathrm{vir}}}\frac{1}{e_{\BC_t^*}(N^{\mathrm{vir}})}\in \BQ\\
		&\VWs_{0, \mathsf w}^{JS,0}(\beta):=\int_{[\Ms^{JS,\BC_t^*}]^{\mathrm{vir}}}\frac{e^{\mu(\beta)}}{e_{\BC_t^*}(N^{\mathrm{vir}})}\in \BQ[t^{\pm}].
		\end{align*}
\end{defn}

\begin{conj}
	The invariants  $\VWs^{0}_{0, \mathsf w}$ and $\VWs_{0,\wsf}^{0}(\beta)$  are independent of the choice of the polarisation $\CO(m)$  and therefore are well-defined. 
\end{conj}
\begin{rmk} It was communicated to the author that the arguments which were used in \cite[Section 4.3.9]{Liu} to prove a similar statement in the K-theoretic set-up can be used to prove the conjecture above. 
	\end{rmk}
\begin{defn} \label{zeroinvariants}
	We define  Vafa--Witten invariants $\VWs^{0}_{0, \mathsf w}(\beta)$ by the following equation
	\begin{align*}
		&\VWs_{0, \wsf}^{JS,0}=(-1)^{\chi(\ch(F(n)))-1}\chi(\ch(F(m)))\cdot \VWs^{0}_{0, \mathsf w}\\
	&\VWs_{0, \wsf}^{JS,0}(\beta)=(-1)^{\chi(\ch(F(n)))-1}\chi(\ch(F(m)))\cdot \VWs^{0}_{0, \mathsf w}(\beta)\\ 
	&\VWp^{0,0}_{0, \mathsf w}=\VWs^{0}_{0, \mathsf w}/\rsf^{2g+2}\\
	&\VWp^{0,0}_{0, \mathsf w}(\beta)=\VWs^{0}_{0, \mathsf w}(\beta)/\rsf^{2g+2}.
	\end{align*}
\end{defn}

When there are no insertions, formulas in Definition \ref{zeroinvariants} are wall-crossings formulas between invariants associated to Joyce--Songs pairs and generalised (localised) Donaldson--Thomas invariants. Once the insertions are introduced, these formulas can be taken as a definition of generalised Donaldson--Thomas invariants with insertions, as proposed in \cite{J}. 

Notice that for Definition \ref{zeroinvariants} to be reasonable, one needs $\CO(1)$ to be generic in the sense of \cite[(2.4)]{TT2} (see also \cite[Definition 3.4]{JK}). And, indeed, a stability of the form $\CO_X(1)\boxtimes \CO_E(k)$ is generic in the sense of \cite[(2.4)]{TT2}, as long as the Jacobian of $X$ does not have $E$ as an isogenous factor and $k$ is a non-rational number. In this case, N\'eron--Severi group of $X\times E$ is equal to $H^2(X,\BZ) \oplus H^2(E,\BZ)$ (no  classes from $H^1(X,\BZ) \otimes H^1(E,\BZ)$), and since $k$ is a non-rational number, the genericity of  \cite[(2.4)]{TT2} is satisfied if $(\dsf,\msf,\asf)=(0,0,0)$. By deformation invariance of Vafa--Witten invariants, such setting is not restrictive,  because there always exists such curve $X$. 

\begin{rmk}A moduli space of stable sheaves $\check N_\CO(X\times E, 0, 0, \wsf)$ embeds into $\Ms_\CO(X\times E, 0, 0, \wsf)$ as Higgs sheaves with a zero Higgs field, and is referred to as the \textit{instanton branch} of $\Ms_\CO(X\times E, 0, 0, \wsf)$.  Over the instanton branch, the invariants from Definition \ref{zeroinvariants} are the invariants defined in \cite[Theorem 1.20 \textbf{(a)(iii)}]{J}. The definition of the invariants $\VWp^{0,0}_{0, \mathsf w}$ and $\VWp^{0,0}_{0, \mathsf w}(\beta)$ is justified by Lemma \ref{slpgl}. 
\end{rmk}
\subsection{Quasimap invariants for $\dsf=0$}
Using Theorem \ref{vwqm} as a justification, we make the following definitions.  
\begin{defn} \label{defnnoncoprime}
	We define quasimap invariants for $\Ms(0)$ and $\Mp(0)$ as follows. If $\wsf\neq0$, then 
	 \begin{align*}
	 	&\check{\mathsf{QM}}^{\asf, \bullet}_{0,  \wsf}:= [\check{\mathsf{VW}}^\asf_{0,\wsf}(\Asf_\wsf )]_t/\rsf^2\\
	 	&\hat{\mathsf{QM}}^{\asf, \bullet}_{0, \wsf}:=\sum _{\msf} [\hat{\mathsf{VW}}^{\msf, \asf}_{0, \wsf}(\Asf_\wsf )]_t.
	 \end{align*}
 If $\wsf=0$, then
	\begin{align*}
		&\check{\mathsf{QM}}^{\asf}_{0,  0}:= \check{\mathsf{VW}}^\asf_{0,0}/\rsf^2\\
		&\hat{\mathsf{QM}}^{\asf}_{0, 0}:=\sum_\msf \hat{\mathsf{VW}}^{\msf,\asf}_{0, 0}.
	\end{align*}
\end{defn} 
\begin{rmk}
	The Joyce--Song pairs on $X\times E$ are not $E$-equivariant, we therefore cannot define quasimap invariants in Definition \ref{defnnoncoprime} via quotient moduli spaces $[\Ms^{JS}(X\times E, 0, 0, \wsf)/E]$, as there is simply no $E$-action. One could define invariants via weighted  Euler characteristics  in the spirit of \cite{OS} and \cite{GPT}, which does not require Joyce--Song pairs. However, as it is explained in \cite{TT,TT2}, in general Vafa--Witten invariants defined via weighted Euler characteristics are not equal to Vafa--Witten invariants defined via virtual fundamental classes. Moreover, the latter matches predictions from mathematical physics and therefore are the right ones to be studied. However, the two types of Vafa--Witten invariants might agree in this case. 
	\end{rmk}
\section{Enumerative mirror symmetry} \label{sectionMS}

\subsection{Conjectures}
\begin{defn}  We define the following formal generating series in $\BQ[\![q]\!]$, 
	\begin{align*}
\QMs^{\asf}_{\dsf}(q)&:= \sum_{\wsf>0} \QMs^{\asf, \bullet}_{\dsf, \mathsf w}q^{\mathsf w}\\
\QMp^{\asf}_{\dsf}(q)&:= \sum_{\wsf>0} \QMp^{\asf, \bullet}_{\dsf, \mathsf w}q^{\mathsf w}.
\end{align*}
such that if $\dsf=0$, we use invariants from Definition \ref{defnnoncoprime}. 
\end{defn}
We will now present the conjectural expressions for the series above for a prime rank $\rsf$,  as well as expressions for $\wsf=0$ invariants.  In Section \ref{MM} we provide a derivation of these expressions based on physical computations from \cite{MM} for $\rsf=2$. We therefore separate $\rsf=2$ from arbitrary prime rank $\rsf$ for convenience.   Let us start with  $\rsf=2$. Let 
\[\tilde{\eta}(q)=\prod^{\infty}_{k=1}(1-q^k).\]
We then define
\begin{align*}
	\Usf_1(q)&=\log \tilde{\eta}(q^4)\\
	\Usf_2(q)&=  \log \tilde{\eta}(q)\\
	\Usf_3(q)&=  \log\tilde{\eta}(-q).
\end{align*}
\begin{rmk}
	Even though all of $\Usf_i(q)$ are related to each other by a substitution of variables, we want to distinguish them. The reason is that in the setting of Vafa--Witten invariants,  $\Usf_1(q)$ is responsible for the \textit{monopole branch}, while $\Usf_2(q)$  and $\Usf_3(q)$ are responsible for the \textit{instanton branch} (see \cite[(1.9),(1.10)]{TT} for the definitions of the instanton and monopole branches). In the language of quasimaps, $\Usf_2(q)$ and $\Usf_3(q)$ are invariants associated to quasimaps mapping to the main component of the nilpotent cone, while $\Usf_1(q)$ are the invariants associated to quasimaps mapping to other components. With respect to $S$-duality, these generating series are exchanged in simple way, see Conjecture \ref{MS1}.
	\end{rmk}
\begin{conjecture}[Expressions, $\wsf>0$] \label{wnonzero} If $g\geq 2$ and $\rsf=2$, then

	\begin{align*}
		\QMs^0_{\dsf}(q)&=(-1)^\dsf(2-2g)2^{4g-1}  \Usf_1(q) \\	
		\QMs^1_{\dsf}(q)&=(2-2g) 2^{2g-1 }(\Usf_2(q)+(-1)^\dsf\Usf_3(q)) \\
		\\
		\QMp^0_{\dsf}(q)&=(-1)^{\dsf}(2-2g)2^{2g-1}\Usf_1(q) \\
		\QMp^1_{\dsf}(q)&=(2-2g) (2^{4g-1}  \Usf_2(q) +(-1)^{\dsf} 2^{2g-1}  \Usf_3(q)).
	\end{align*}
\end{conjecture}

\begin{conjecture}[Expressions, $\wsf=0$] \label{wzero}If $g\geq 2$ and $\rsf=2$, then 
	\begin{align*}
	&	\QMs^0_{\dsf,0}=(-1)^\dsf2^{4g-3} & &\QMp^0_{\dsf,0}=(-1)^\dsf2^{2g-3}  \\
		&\QMs^{\asf\neq 0}_{0,0}=2^{2g-2}& &\QMp^{\asf\neq 0}_{\dsf,0}=2^{4g-3}+(-1)^\dsf2^{2g-3} \\
	    	&\QMs^{\asf \neq 0}_{\dsf\neq 0,0}=0 & 	&
	\end{align*}
If $g=1$ and $\rsf \geq 2$, then 
\begin{align*}
	&\QMs^0_{0,0}= \rsf^2+1 & &\QMp^0_{0,0}= \rsf+2-\rsf^{-1}\\
	&\QMs^0_{\dsf \neq0,0}=1  &  &\QMp^0_{\dsf\neq0,0}=\rsf+1-\rsf^{-1}  \\
	&\QMs^{\asf \neq0}_{0,0}=1 & &\QMp^{\asf\neq 0}_{0,0}=\rsf+1-\rsf^{-1} \\
	&\QMs^{\asf\neq0}_{\dsf\neq0,0}=0 & &\QMp^{\asf\neq0}_{\dsf\neq0,0}= \rsf-\rsf^{-1}.
\end{align*}
\end{conjecture}

Using the permutation action of $S_3$ on the linear span $ \BQ\langle \Usf_i(q) \rangle $ as in Section \ref{conjsec}, 
\[ \sigma \cdot \Usf_i(q):=\Usf_{\sigma(i)}(q), \quad \sigma \in S_3.\]
 we get the following relation from the formulas above.
\begin{conjecture}[Enumerative mirror symmetry, $\wsf>0$] \label{MS1}If $g\geq 2$ and $\rsf=2$, then 
\[(12)\cdot 2\QMs^{\asf}_{\dsf}(q)=\sum_{\dsf'} \sum_{\mathsf a'} (-1)^{  \mathsf d \cdot \mathsf a'+\dsf'\cdot \asf' } \hat{\mathsf{QM}}^{\asf'}_{\dsf'}(q).\]
where the sum is taken over all $\mathsf a' \in \BZ_2$ and all $\dsf' \in \BZ_2$. 
\end{conjecture}

\begin{conjecture}[Enumerative mirror symmetry, $\wsf=0$] \label{Enumerativemirror}
	 If $\rsf\geq 2$, then 
	 \begin{align*}
	 	& \rsf\check{\mathsf{QM}}^{\asf}_{\dsf,0}=\sum_{\dsf'} \sum_{\mathsf a'} \omega^{ ( \dsf'\cdot \asf + \mathsf d \cdot \mathsf a' )} \hat{\mathsf{QM}}^{\asf'}_{\dsf',0},
	\end{align*}
	where  $e^{\frac{2\pi \sqrt{-1}}{\rsf}}$; the sum is taken over all $\mathsf a' \in \BZ_\rsf$ and all $\dsf'\in \BZ_\rsf$.
	\end{conjecture}
\begin{rmk}
	If $g(X)=1$, a moduli space $\Ms(\dsf)$ is a point for $\dsf\neq 0$. Hence for $g(X)=1$, only the statement of Enumerative mirror symmetry for $\wsf=0$ makes sense. 
\end{rmk}
\subsection{Higher rank} \label{higherrank} Let us now consider invariants for higher prime rank $\rsf$. Manschot--Moore's calculations apply only to $\rsf=2$, but there is an obvious extension of  Conjecture \ref{wnonzero} and \ref{MS1}, which satisfy basic constraints and agree with \cite[Proposition 4.14]{Nqm} and \cite[Conjecture A]{Nqm}. Let us now explain this. We introduce the following generating series. 
\begin{align*}
	\Usf_1&= \log(\tilde\eta(q^{2\rsf})) \\
	\Usf_{2+l}&= \log( \tilde\eta(e^{\frac{2\pi \sqrt{-1} l}{\rsf}}q)), \quad l \in \{0	,\dots, \rsf-1 \}
\end{align*}
Let $\BC\langle \Usf_i \rangle$ be the linear span of $\Usf_i$. We define $n_l \in \{1	,\dots, \rsf-1 \}$ by the equation $l\cdot n_l=-1 \ \mathrm{mod} \ \rsf$. We then define a permutation $\sigma$ on $\{1	,\dots, \rsf+1 \}$ as follows
\begin{align*}
	&\sigma(1)=2, \\
&\sigma(2)=1, \\
	&\sigma(2+l)=2+n_l, \quad \ l \in \{1	,\dots, \rsf-1 \}.
\end{align*} 
The extension of Conjecture \ref{wnonzero} takes the following form.

\begin{conjecture} \label{wnonzeror} If $g\geq  2$ and $\rsf$ is prime, then
	\begin{align*}
		&\QMs^0_{0}(q)=(2-2g)(\rsf-1)\rsf^{4g-1} \Usf_1(q) \\
		&\QMs^0_{\dsf\neq 0}(q)=(2g-2)  \rsf^{4g-1}  \Usf_1(q) \\	
		&\QMs^{\asf \neq0}_{\dsf }(q)=(2-2g)  \rsf^{2g-1 } \left(\sum^{\rsf+1}_{i=2} \omega^{-\dsf \cdot \asf (i-2)}\Usf_i(q)\right) \\
		\\
		&\QMp^0_{0}(q)=(2-2g) (\rsf-1)\rsf^{2g-1}  \Usf_1(q) \\
		&\QMp^0_{\dsf\neq 0}(q)=(2g-2)  \rsf^{2g-1}  \Usf_1(q) \\
		&\QMp^{\asf \neq0}_{\dsf }(q)=(2-2g) \left( \rsf^{4g-1} \Usf_2(q) + \rsf^{2g-1} \sum^{\rsf+1}_{i=3} \omega^{-\dsf \cdot \asf (i-2)}\Usf_i(q)\right),
	\end{align*}
	where $\omega=e^{\frac{2\pi \sqrt{-1}}{\rsf}} $.
\end{conjecture}
Using the flux equation from \cite[Section 3.4]{JK}, one can verify that these formulas satisfy a basic expectation -  Enumerative mirror symmetry. 
\begin{conjecture}
	If  $g\geq 2$ and $\rsf$ is prime, then both $\QMs^{\asf}_{\dsf}(q)$ and $\hat{\mathsf{QM}}^{\asf}_{\dsf}(q)$ are in  $\BC\langle \Usf_i(q) \rangle $, such that 
	\[\sigma \cdot \rsf\QMs^{\asf}_{\dsf}(q)=\sum_{\dsf'} \sum_{\mathsf a'} \omega^{  \mathsf d \cdot \mathsf a'+\dsf'\cdot \asf' } \hat{\mathsf{QM}}^{\asf'}_{\dsf'}(q),\]
	where $\omega=e^{\frac{2\pi \sqrt{-1}}{\rsf}}$. 
\end{conjecture}

The expressions above are in fact much simpler than they appear and all of them lie in $\BQ[\![q]\!]$. For a power series $a(q)=\sum_n a_n q^n$ and a fixed $\rsf$, let us introduce the following notation 
\[[a(q)]_{k}:=\sum_{n=k \ \mathrm{mod} \ \rsf} a_n q^n.\]

 Then  by \cite[Lemma 3.6]{JK} we have the following identity
\begin{align*}
\sum^{\rsf+1}_{i=2} \omega^{-k (i-2)}\Usf_i(q) &= \rsf [\Usf_2]_k,
\end{align*}
which simplifies our expressions above

	\begin{align*}
	&\QMs^{\asf \neq0}_{\dsf }(q)=(2-2g)  \rsf^{2g }[\Usf_2]_{\dsf \cdot \asf },\\
	&\QMp^{\asf \neq0}_{\dsf }(q)=(2-2g) \left( (\rsf^{4g-1} -\rsf^{2g-1})\Usf_2(q) + \rsf^{2g}[\Usf_2]_{\dsf \cdot \asf }\right).
\end{align*}

In particular,  for $\Ms(\dsf\neq0)$ invariants of degree $(\wsf,\asf)$ vanish unless  $\wsf=\dsf \cdot \asf \ \mathrm{mod} \ \rsf$.  This is a basic constraint  for  these invariants, see Remark \ref{asfnot}.  Moreover, if we forget about $\asf$ (i.e.\ sum over it), 
\begin{align*} 
\QMs_{\dsf }(q)&:= \sum_\asf \QMs^\asf_{\dsf }(q) \\
\QMp_{\dsf }(q)&:= \sum_\asf \QMp^\asf_{\dsf }(q),
\end{align*}
 then we obtain that the resulting generating series are independent of  degree $\dsf$, which can rightfully be  called \textit{quantum $\chi$-independence}. 
\begin{conjecture}[\textit{Quantum $\chi$-independence}] \label{independence} If $\gcd(\rsf, \dsf)=1$ and $\gcd(\rsf, \dsf')=1$, then
\begin{align*} 
\QMs_{\dsf }(q)&= \QMs_{\dsf'}(q) \\
\QMp_{\dsf }(q)&= \QMp_{\dsf'}(q).
\end{align*}
\end{conjecture}

Differentiating or exponentiating generating series from Conjecture \ref{wnonzeror} (and after adding some terms to make them modular), one can obtain series that satisfy Enumerative mirror symmetry with respect to the change of variables $\tau \rightarrow -1/\tau$.  The form of permutation $\sigma$ can be justified by the following transformation properties of the Dedekind eta function, 
\begin{align*}
\eta(\rsf \tau)\mid_{-1/\tau}&=(\substack{\frac{\tau}{\rsf}})^{1/2}\eta(\substack{\frac{\tau}{\rsf}}),\\
\eta(\substack{\frac{\tau+l}{\rsf}})\mid_{-1/\tau}&=\tau^{1/2}\eta(\substack{\frac{\tau+n_l}{\rsf}}), \quad l \in \{1,\dots,\rsf -1\}.
\end{align*}

Using Proposition \ref{ggreat2} and \cite[Corollary 3.9]{Nqm}, Conjecture \ref{wzero} has an obvious extension to an arbitrary prime rank $\rsf$.
\begin{conjecture} \label{wzero2}If $g\geq 2$ and $\rsf$ is prime, then 
	\begin{align*}
		&	\QMs^0_{0,0}=(\rsf-1)\rsf^{4g-3} & &\QMp^0_{0,0}=(\rsf-1)\rsf^{2g-3}  \\
		&\QMs^0_{\dsf\neq 0, 0}= - \rsf^{4g-3} & 	&\QMp^0_{\dsf\neq 0,0}= -\rsf^{2g-3} \\
		&\QMs^{\asf\neq 0}_{0,0}=\rsf^{2g-2}& &\QMp^{\asf\neq 0}_{0,0}=\rsf^{4g-3}+(\rsf-1)\rsf^{2g-3} \\
		&\QMs^{\asf \neq 0}_{\dsf\neq 0,0}=0 & 	&\QMp^{\asf\neq0 }_{\dsf\neq0,0}=\rsf^{4g-3}-\rsf^{2g-3}.
	\end{align*}
\end{conjecture}
It can be checked that the formulas above also satisfy Enumerative mirror symmetry, Conjecture \ref{Enumerativemirror}.
\subsection{Gerbes and extended degree} \label{gerbesandclasses}
The class 
\[\check \alpha \in H^2_{\text{lis-\'et}}( \FMsr(\dsf),\BZ_\rsf)\]
 can be used to do two things. The first one is to define a $\BZ_\rsf$-gerbe $\FMs_{\mathrm{SL} }(\dsf)$ from Section \ref{gerbeclasses}. The second one is to define an extended degree of quasimap, i.e.\ the class $\asf$ in Definition \ref{qmsell},  see Remark \ref{asfnot}. We will now discuss the relation between the two. 
 
Firstly, using the moduli theoretic description of $\FMs_{\mathrm{SL} }(\dsf)$ from Section \ref{gerbeclasses}, we can identify a moduli space of quasimaps  $Q^{0^+}_E(\check{M}_{\alpha}(\dsf),\wsf)$ with a moduli space of sheaves in the spirit of Lemma \ref{ident}. 
\begin{prop} \label{ident3} Assume $\dsf\neq 0$. By associating to a quasimap $f\colon E \rightarrow \FMs_{\mathrm{SL} }(\dsf)$ the corresponding family $(F,\phi)$, we obtain an identification
	\[Q^{0^+}_E(\Ms(\dsf),\wsf) \cong \check{\CM}_{p_X^*L}(X\times E,\dsf, \wsf), \quad f\mapsto (F,\phi),\]
where $\check{\CM}_{p_X^*L}(X\times E,\dsf, \wsf)$ is a $\BZ_\rsf$-gerbe over $\Ms_{p_X^*L}(X\times E,\dsf, \wsf)$ obtained by adding the data of $\psi \colon  \det(F) \xrightarrow{\sim} p^*_XL$. 
\end{prop}

\textit{Proof.} By the description of the gerbe $\FMs_{\mathrm{SL} }(\dsf)$ from Section \ref{gerbeclasses}, a quasimap $f\colon C \rightarrow \FMs_{\mathrm{SL} }(\dsf)$ is given by a family of Higgs sheaves $(F,\phi)$ with  $\psi \colon  \det(F) \xrightarrow{\sim} p^*_XL$.  The rest is as in Lemma \ref{ident}.
 \qed
 \\

 By Lemma \ref{degreechern}, the degree of a quasimap $f \colon C  \rightarrow \FMsr(\dsf)$ has the following expression in terms of the associated family, 
 \[\deg(f)=\dsf \cdot \ch_1(F)_{\mathrm d}-\rsf \cdot  \ch_2(F)_{\mathrm d}.\]
 In particular, if
  \[\deg(f)=\mathrm{c}_1(F)_{\mathrm d}=0 \mod \rsf,\]
 then by  Lemma \ref{ul} there exists a natural map
 \begin{equation} \label{sltou}
 \Ms_{p^*_XL}(X\times E,\dsf, \asf, \wsf) \rightarrow \Ms_{\usf}(X\times E,\dsf, \wsf),
 \end{equation}
which is a $\Jac(C)[\rsf]$-torsor. Equivalently, one can obtain this map as follows. By 
composing a quasimap to $\check{\FM}_{\mathrm{SL} }(\dsf)$ with the natural projection to $\FMsr(\dsf)$, we obtain a map 
\[Q^{0^+}_C(\check{M}_{\mathrm{SL} }(\dsf),\wsf) \rightarrow  Q^{0^+}_C(\check{M}(\dsf),\wsf),\]
if $\wsf= 0\mod \rsf$. This is exactly the map (\ref{sltou}) after identifications of Lemma \ref{ident} and Lemma \ref{ident3} and rigidification of the $\BZ_\rsf$-gerbe structure. 

Let $\QMs^{\mathrm{SL} }_{\dsf, 0}$ and $\QMs^{\bullet,\mathrm{SL} }_{\dsf, \mathsf w}$ be the invariants associated to $Q^{0^+}_E(\check{M}_{\mathrm{SL} }(\dsf),\wsf)$ and defined in the same way as invariants from Definition \ref{qmsell}. By the preceding discussion, we obtain the following comparison result. 
\begin{lemma} \label{gerbelemma}If $\wsf=0 \mod \rsf$ and $\dsf\neq0$, then 
	\begin{align*}
	\QMs^{\bullet,\mathrm{SL} }_{\dsf, \mathsf w}=\rsf \cdot \QMs^{0, \bullet}_{\dsf, \mathsf w} \\
	\QMs^{\mathrm{SL} }_{\dsf, \mathsf w}=\rsf \cdot \QMs^{0}_{\dsf, \mathsf w}.
	\end{align*}
	\end{lemma}

\section{Derivation of conjectures from Manschot--Moore} \label{MM}
Let 
\[B:=\mathbb 1\boxtimes [\mathrm{pt}] \in H^2(X\times E,\BZ).\] We define 
\begin{align*}
&\VWs^{\asf}_{\dsf}(B,q):= \sum_{\wsf>0}\VWs^{\asf}_{\dsf, \wsf}(B)q^{\wsf} \\
&\VWp^{\msf,\asf}_{\dsf}(B,q):= \sum_{\wsf>0}\VWp^{\msf, \asf}_{\dsf, \wsf}(B)q^{\wsf} \\
&\VWp^{\asf}_{\dsf}(B,q):= \sum_{\msf}\VWp^{\msf,\asf}_{\dsf}(B,q).
\end{align*}
These are formal powers series in the ring $\BQ[\![q]\!][t^{\pm}]$. Proposition \ref{rigid} then implies the following lemma. 
\begin{lemma} \label{difeq}
\begin{align*}
	&q\frac{ d\QMs^{\asf}_{\dsf}(q)}{dq}= [\VWs^{\asf}_{\dsf}(B,q)]_t/\rsf \\ 
	&q \frac{ d\QMp^{\asf}_{\dsf}(q)}{dq}= \rsf  [\VWp^{\asf}_{\dsf}(B,q)]_t,
	\end{align*}
where  $[\dots ]_t$ is the coefficient of an element in $\BQ[\![q]\!][t^{\pm}]$ at $t$.

\begin{rmk} \label{diveq}

	The operation $q\frac{d}{dq}$ is also natural from the point of view of Gromov--Witten theory, as it can be used to express the divisor equation
	$\langle \gamma_{1}, \dots, \gamma_n, D \rangle^{\infty}_{g,n+1,\beta}= \int_\beta D \cdot \langle \gamma_{1}, \dots, \gamma_n \rangle^{\infty}_{g,n,\beta},$
	which suggests that perhaps a different kind of invariants is more suitable for modularity properties discussed in Section \ref{sduality}. 
\end{rmk}

\end{lemma}
From now on, we restrict to the case $\rsf=2$. Using \cite[Section 7]{MM}, we will derive (conjectural) expressions for $\VWs^{\mathsf a}_{\dsf}(B,q)$ and $\VWp^{\mathsf a}_{\dsf}(B,q)$, as well as for  $\VWs^{\mathsf a}_{\dsf,0}$ and $\VWp^{\mathsf a}_{\dsf,0}$, if $\rsf=2$.  We therefore also obtain conjectural expressions for quasimap invariants, using Theorem \ref{vwqm}. We refer  to  \cite{GK1,GK,GKL, GKW}) for a more mathematical treatment of related invariants, which however does not cover $\mu$-insertions for the monopole branch. 

\vspace{0.3cm}
\noindent \textbf{Important.} We take the liberty to extend formulas of \cite[Section 7]{MM} from the case of surfaces with $b_1=0$ and $p_g(S)>0$ to the case of surfaces with $b_1\neq0$ and $p_g(S)>0$ in the way it is done for invariants without insertions in \cite[Section 6]{DPS}.
\vspace{0.3cm}

To derive formulas from \cite{MM}, we need to define certain modular forms, compute  the Seiberg--Witten invariants and calculate the basic topological invariants of $X\times E$.  Formulas for $g=1$ can be obtained from \cite[(5.25)]{LL} without any effort, they are verified in Section \ref{calc}.
\subsection{Modular forms} Recall the following fundamental modular forms. 
\begin{align*}
	\theta_2(q)&=\sum_{n \in \BZ+\frac{1}{2}} q^{\frac{n^2}{2}}  &  E_2(q)&=1-24\sum_{n=1}^{\infty}\sigma_{1}(n)q^n \\
	\theta_3(q)&=\sum_{n \in \BZ} q^{\frac{n^2}{2}} &  G_2(q)&=\frac{1}{24}E_2(q). \\
	\theta_4(q)&=\sum_{n \in \BZ} (-1)^nq^{\frac{n^2}{2}} & 	
\end{align*}
The series $\theta_i(q)^4$ are modular forms for $\Gamma(2)$ of weight 2, the space of which is 2-dimensional (see \cite[Proposition 3]{Z} and \cite[Section 3]{Z}). It is therefore enough to check the first two coefficients to verify the following identities,

\begin{align} \label{genser}
	\begin{split}
	&\CU_1(q):=\frac{\theta_3^4+\theta_4^4+2E_2}{12}=-8G_2(q^2)\\
&\CU_2(q):=\frac{-\theta_2^4-\theta_3^4+2E_2}{12}= -2G_2(q^{\frac{1}{2}})\\
&\CU_3(q):=\frac{\theta_2^4-\theta_4^4+2E_2}{12}= -2G_2(-q^{\frac{1}{2}}).
\end{split}
\end{align}
\subsection{Seiberg--Witten invariants} \label{SWch}By \cite[Proposition 4.4]{FM} or by \cite[Proposition 39]{Brus}, Seiberg--Witten basic classes of $X\times E$ are 
 \[a_k:=(k,0,0) \in H^2(X\times E, \BZ), \quad 0\leq  k\leq 2g-2,\] 
such that the associated Seiberg--Witten invariants are
\[\SW(a_k)=(-1)^k \binom{2g-2}{k}.\]
Notice the following identities, which we will need later. 

\begin{align} \label{d0}
\begin{split}
&\sum^{2g-2}_{i=0}(-1)^k \binom{2g-2}{k}=0 \\
& \sum^{2g-2}_{i=0}(-1)^k \binom{2g-2}{k}2k=0\\
& \hspace{-0.4cm}\sum_{\mathrm{odd/even}\ k} \binom{2g-2}{k}(2g-2-2k)=0\\
&\sum_{\mathrm{odd} \ k}\binom{2g-2}{k}= 2^{2g-3}. \\
\end{split}
\end{align}

\subsection{Topological invariants}
Finally, we list the basic topological invariants of $X\times E$ needed for the formulas in \cite{MM}.
\begin{align} \label{topdata} 
	&\chi(X\times E)=0 & &a_k^2=0 \nonumber \\
	&\sigma(X\times E)=0 &&B^2=0\\ 
	& \chi(\CO_{X\times E})=0 &&\mathrm{c}_1(X\times E) \cdot B=2g-2 \nonumber \\
	&\mathrm{c}_1(X\times E)=0 &&(\mathrm{c}_1(X\times E)-2a_k)\cdot B=2g-2-2k.\nonumber
\end{align}	

\subsection{Formulas for $\wsf=0$} We now can state the formulas, starting with $\wsf=0$ case. We make the following notational shortcuts for this section, 
\[H^2:=H^2(X\times E,\BZ) ,\quad w:=(\dsf, \msf,\asf) \in H^2(X\times E,\BZ_2).\]
We will consider the standard  intersection pairing on  $H^2(X\times E,\BZ)$,
\begin{align*}
H^2(X\times E,\BZ)&\times H^2(X\times E,\BZ) \rightarrow \BZ \\
(v,v')&\mapsto v\cdot v'.
\end{align*}
Since  $H^2(X\times E,\BZ)$  is an even lattice,  $\frac{v\cdot v}{2}$ takes values in integers. Hence for an element $w\in H^2(X\times E,\BZ_2)$, the expression $\frac{w\cdot w}{2}$  is defined via a lift of $w$ to $H^2(X\times E,\BZ)$.

Plugging (\ref{topdata}) into the formulas from \cite[(7.2),(7.23),(7.32),(7.40)]{MM} and using the dictionary of \cite[Table 7]{MM}, we obtain
\begin{equation} \label{VWZ}
	 \VWp^{\msf,\asf}_{\dsf,0}=\sum^{j=3}_{j=1} Z^{\msf,\asf}_{j,\dsf} \hspace{4cm} \quad (7.2) \text{ in \cite{MM},}
	\end{equation}
where
\begin{align*} 
	&Z^{\msf,\asf}_{1,\dsf}=\sum_{a_k=w \mod 2H^2} \SW(a_k) &(7.23) \text{ in \cite{MM},}\\
	&Z^{\msf,\asf}_{2,\dsf}=2^{-2g-1}\sum_{a_k} (-1)^{ a_k\cdot w }\SW(a_k) &(7.32) \text{ in \cite{MM},}\\
	&Z^{\msf,\asf}_{3,\dsf}= 2^{-2g-1} (-1)^{\frac{w\cdot w }{2}} \sum_{a_k} (-1)^{a_k\cdot w}\SW(a_k) &(7.40) \text{ in \cite{MM}.}
\end{align*}
Using Section \ref{SWch}, we can readily obtain expressions for $Z^{\msf,\asf}_{j,\dsf}$, which are displayed in the table below.
\begin{table}[h!]
	\begin{tabular}{|c|c|c|c| } 
		\hline
	$w$	& $Z^{\msf,\asf}_{1,\dsf}$ & $Z^{\msf,\asf}_{2,\dsf}$ & $Z^{\msf,\asf}_{3,\dsf}$ \\ 
		\hline
	$(0,0,0)$	& $2^{2g-3}$ & 0 & 0  \\ 
	$(1,0,0)$	& $-2^{2g-3}$ & 0 &0 \\ 
	$(0, \msf \neq 0 ,0)$	& 0& 0 &0  \\
	$(1, \msf \neq 0, 0)$	& 0& 0 &0   \\
		$(0, \msf ,1)$	& 0& $2^{-3} $&  $(-1)^{\frac{ \msf\cdot \msf}{2}}2^{-3}$  \\
		$(1, \msf ,1)$	& 0& $2^{-3} $&  $-(-1)^{\frac{ \msf \cdot \msf}{2}}2^{-3}$   \\
		\hline
	\end{tabular}
\caption{}
\end{table}

Using Theorem \ref{vwqm} and Definition \ref{defnnoncoprime}, we therefore obtain the following expressions for quasimap invariants 
\begin{align*}
&\QMp^0_{0,0}=\VWp^0_{0,0}=\sum _{\msf } \VWp^{\msf,0}_{0,0}=2^{2g-3} \\
&\QMp^0_{1,0}=\VWp^0_{1,0}=\sum _{\msf } \VWp^{\msf,0}_{1,0}=-2^{2g-3}\\
&\QMp^1_{0,0}=\VWp^1_{0,0}=\sum _{\msf } \VWp^{\msf,1}_{0,0}=2^{-2}N_0 \\
&\QMp^1_{1,0}=\VWp^1_{1,0}=\sum _{\msf } \VWp^{\msf,1}_{1,0}=2^{-2}N_1 \\
\\
&\QMs^0_{0,0}= 2^{2g}\cdot\VWp^{0,0}_{0, 0}=2^{4g-3} \\
&\QMs^0_{1,0}=  2^{2g}\cdot\VWp^{0,0}_{1, 0}=- 2^{4g-3} \\	
&\QMs^1_{0,0}=2^{2g}\cdot\VWp^{0,1}_{0, 0}=2^{2g-2}\\
&\QMs^1_{1,0}=2^{2g}\cdot\VWp^{0,1}_{1, 0}=0,
\end{align*}
where 
\begin{align*}
&N_0:= |\{ \msf \in H^1(X, \BZ_2 )\otimes H^1(E, \BZ_2 ) \mid \frac{\msf \cdot \msf }{2}=0\}| \\
 &N_1:= |\{ \msf\in H^1(X, \BZ_2 )\otimes H^1(E, \BZ_2 ) \mid \frac{\msf \cdot \msf}{2}=1\}|. 
 \end{align*}
Let us now determine $N_0$ and $N_1$. By picking a symplectic basis of $H^1(X,\BZ)$ and $H^1(E,\BZ)$, one can readily establish the isomorphism of lattices
\[H^1(X, \BZ )\otimes H^1(E, \BZ )\cong U^{2g},\]
where $U$ is a hyperbolic plane. By induction on $g$, one can then  show that 
 \begin{align} \label{npm}
 	\begin{split}
& N_0= 2^{4g-1}+2^{2g-1} \\
& N_1= 2^{4g-1}-2^{2g-1}.
 \end{split}
 \end{align}
 The above formulas therefore give rise to Conjecture \ref{wzero} for $g\geq 2$ and $\rsf=2$.  
\subsection{Formulas for $\wsf>0$}	\label{formulas}
Consider the ring $\BQ[\![q,t]\!]$, we define the additive group
\[\overline {\BQ[\![q,t]\!]}:=\BQ[\![q,t]\!]/\BQ,\] 
i.e.\ $\overline {\BQ[\![q,t]\!]}$ is the group of formal power series up to a constant term. Given two elements $f,g\in \BQ[\![q,t]\!]$, we denote 
\[f \equiv g \iff \overline f = \overline g, \]
where $\overline f$ and $ \overline g$ are the images of $f$ and $g$ in $\overline {\BQ[\![q,t]\!]}$, i.e.\ $f \equiv g$, if $f$ and $g$ are equal up to a constant term.

Substituting the topological data from (\ref{topdata}) into the formulas from \cite[(7.47),(7.48),(7.49)]{MM}, 
 we obtain that 
\begin{equation} \label{vwz}
\VWp^{\msf,\asf}_{\dsf}(B,q)\equiv \sum_j Z^{\msf,\asf}_{j,\dsf}(B,q^2),
\end{equation}
where 
\begin{align*}
	&Z^{\msf,\asf}_{1,\dsf}(B,q)=\sum_{a_k=w \mod 2H^2} \SW(a_k)  \cdot e^{(-((\mathrm{c}_1-2a_k)\cdot B )t \CS_1+(\mathrm{c}_1\cdot B)t  \CU_1) } \\
	&Z^{\msf,\asf}_{2,\dsf}(B,q)= 2^{-2g-1} \sum_{a_k} (-1)^{a_k\cdot w }\SW(a_k)  \cdot e^{( ((\mathrm c_1-2a_k)\cdot B) t \CS_2+ (\mathrm c_1\cdot B) t  \CU_2)} \\
	&Z^{\msf,\asf}_{3,\dsf}(B,q)= 2^{-2g-1} (-1)^{\frac{w\cdot w }{2}}  \sum_{a_k} (-1)^{a_k\cdot w}\SW(a_k) \cdot e^{( ((c_1-2a_k)\cdot B)t \CS_3+(\mathrm c_1\cdot B)t  \CU_3)},
	\end{align*}
which are the formulas (7.47),(7.48) and (7.49) in \cite{MM}, respectively. 
We refer to \cite[Table 6]{MM} for the definition of $\CS_j$, here we do not need it, because it cancels out (see Table \ref{table}).  Note that we make a substitution of variables 
\[q \rightarrow q^2,\]
because \cite{MM}  sums over the \textit{instanton number}, while we sum over $\Delta(F)/2$, which is related to the instanton number by a factor of two (see \cite[Table 7]{MM} and \cite[Appendix C]{GK}). Moreover, \cite{MM} has constant terms in their expression, which are not present in our case. For that reason, the equality of (\ref{vwz}) is only up to a constant term. 

Using the identities from (\ref{d0}) and the topological data from (\ref{topdata}), we obtain the following expressions for the coefficient of $Z^{\msf,\asf}_{j,\dsf}(B,q)$ at $t$, which we denote by $[Z^{\msf,\asf}_{j,\dsf}(B,q)]_t$. 
\begin{table}[h!] 
	\begin{tabular}{|c|c|c|c|} 
		\hline
		$w$	& $[Z^{\msf,\asf}_{1,\dsf}(B,q)]_t$ & $[Z^{\msf,\asf}_{2,\dsf}(B,q)]_t$ & $[Z^{\msf,\asf}_{3,\dsf}(B,q)]_t$ \\ 
		\hline
		$(0,0,0)$	& $(2g-2)2^{2g-3} \CU_1$& 0 & 0 \\ 
		$(1,0,0)$	& $(2-2g)2^{2g-3}\CU_2$ & 0 &0 \\ 
		$(0, \msf \neq 0 ,0)$	& 0& 0 &0  \\
		$(1, \msf \neq 0, 0)$	& 0& 0 &0  \\
		$(0, \msf ,1)$	& 0& $(2g-2)2^{-3}\CU_2$&  $(-1)^{\frac{\msf\cdot \msf }{2}}(2g-2)2^{-3}\CU_3$  \\
		$(1, \alpha ,1)$	& 0& $(2g-2)2^{-3}\CU_2$ & $(-1)^{\frac{\msf\cdot \msf }{2}}(2-2g)2^{-3}\CU_3$  \\
		\hline
	\end{tabular}
	\caption{}\label{table}
\end{table}

\noindent The coefficient of $t$ in $\VWp^{\asf}_{\dsf}(B,q)$ can now be obtained from Table \ref{table}, 
\begin{align} \label{exps}
	\begin{split}
	&[\VWp^0_{0}(B,q)]_t=\sum _{\msf} [\VWp^{\msf,0}_{0}(B,q)]_t\equiv (2g-2)2^{2g-3}\CU_1(q^2)\\
	&[\VWp^0_{1}(B,q)]_t=\sum _{\msf} [\VWp^{\msf,0}_{1}(B,q)]_t\equiv(2-2g)2^{2g-3} \CU_1(q^2)\\
	&[\VWp^1_{0}(B,q)]_t=\sum _{\msf} [\VWp^{\msf,1}_{0}(B,q)]_t\equiv(2g-2) (2^{4g-3}\CU_2(q^2)+2^{2g-3}\CU_3(q^2))\\
	&[\VWp^1_{1}(B,q)]_t=\sum _{\msf} [\VWp^{\msf,1}_{1}(B,q)]_t\equiv(2g-2) (2^{4g-3}\CU_2(q^2)-2^{2g-3}\CU_3(q^2)).
	\end{split}
\end{align}
For the last expressions we also used (\ref{npm}).

Let us now determine quasimap invariants, using Lemma \ref{difeq}. We define
\begin{align*}
	&\eta(q)=q^{\frac{1}{24}}\prod^{\infty}_{k=1}(1-q^k). 
	\end{align*}
 Recall that 
\[G_2(q)=q \frac{d \log(\eta(q))}{dq},\]
hence 
\[ 4G_2(q^4)=q \frac{d \log(\eta(q^4))}{dq},\quad G_2(-q)=q \frac{d \log(\eta(-q))}{dq}. \]
Using the expressions from (\ref{genser}), we therefore obtain the antiderivatives of $\CU_j(q^2)$, 
\begin{align*}
	&q \frac{d(-2 \log(\eta(q^4)))}{dq} =\CU_1(q^2)\\
	&q \frac{d(-2 \log(\eta(q)))}{dq}=\CU_2(q^2) \\
	&q \frac{d(-2 \log(\eta(-q)))}{dq}=\CU_3(q^2).
	\end{align*}
Finally, notice that the factor $q^{\frac{1}{24}}$ in $\eta$ is responsible for the constant coefficient of $G_2(q)$, 
\[q \frac{d\log(\tilde\eta(q))}{dq}=G_2(q)-\frac{1}{24}.\]
By Lemma \ref{difeq}, Theorem \ref{vwqm} and Definition \ref{defnnoncoprime}, we therefore obtain 
	\begin{align*}
	\QMs^0_{0}(q)&=(2-2g) 2^{4g-1} \Usf_1(q) \\
	\QMs^0_{1}(q)&=(2g-2) 2^{4g-1}\Usf_1(q) \\	
	\QMs^1_{0}(q)&=(2-2g)2^{2g-1}\left (\Usf_2(q)+\Usf_3(q)\right )\\
	\QMs^1_{1}(q)&=(2-2g) 2^{2g-1} (\Usf_2(q)-\Usf_3(q))) \\
	\\
	\QMp^0_{0}(q)&=(2-2g)2^{2g-1}\Usf_1(q) \\
	\QMp^0_{1}(q)&=(2g-2)  2^{2g-1} \Usf_1(q) \\
	\QMp^1_{0}(q)&=(2-2g) (2^{4g-1} \Usf_2(q)+2^{2g-1} \Usf_3(q))\\
	\QMp^1_{1}(q)&=(2-2g) (2^{4g-1}\Usf_2(q) - 2^{2g-1}  \Usf_3(q)).	
\end{align*}
The above formulas give rise to Conjecture \ref{wnonzero}. 
\subsection{S-duality} \label{sduality}
Recall the transformation law of the theta series and second Eisenstein series, 
\begin{align*}
	\frac{1}{\tau^2}E_2(-1/\tau)&=E_2+\frac{6}{\pi i\tau}\\
	\theta_2(-1/\tau)^4&=-\tau^2\theta_4(\tau)^4\\
	\theta_3(-1/\tau)^4&=-\tau^2\theta_3(\tau)^4 \\
	\theta_4(-1/\tau)^4&=-\tau^2\theta_2(\tau)^4
\end{align*}
from which we obtain that 
\begin{align} \label{transform}
	\begin{split}
	&\frac{1}{\tau^2}\CU_1(-1/\tau)+\frac{1}{\pi i\tau}=\CU_2(\tau)\\\
	&\frac{1}{\tau^2}\CU_2(-1/\tau)+\frac{1}{\pi i\tau }=\CU_1(\tau)\\
	&\frac{1}{\tau^2}\CU_3(-1/\tau)+\frac{1}{\pi i\tau}=\CU_3(\tau).
	\end{split}
\end{align}
We define 
\[ \tilde{\VWp}^{\msf,\asf}_{\dsf}(B,q)_t:=\sum_j Z^{\msf,\asf}_{j,\dsf}(B,q)_t, \]
the generating series $\tilde{\VWp}^{\msf,\asf}_{\dsf}(B,q)$ is a minor modification of $\VWp^{\msf,\asf}_{\dsf}(B,q)$, the two are related as follows
\[ \tilde{\VWp}^{\msf,\asf}_{\dsf}(B,q)=\VWp^{\msf,\asf}_{\dsf}(B,q^{\frac{1}{2}})+c,\] 
where $c$ is some constant, which can be determined from (\ref{exps}). However, $\tilde{\VWp}^{\msf,\asf}_{\dsf}(B,q)$ has better modularity properties. Using transformation laws above and expressions in Table \ref{table}, one can then readily verify the following statement (see \cite[Section 7.4]{MM}). 
\begin{conj}[S-duality with $\mu$-insertions] \label{conj1} Let $q=e^{2\pi i \tau}$. If $\rsf=2$, then
	\begin{multline*}
	\frac{1}{ \tau^2}[\tilde{\VWp}^{\msf,\asf}_{\dsf}(B, -1/\tau)]_t+\frac{\mathrm{const}}{\pi i \tau }
	\\ = 2^{-2g-1}  \sum_{(\dsf',\msf',\asf')}(-1)^{(\dsf,\msf,\asf)\cdot (\dsf',\msf',\asf') }[\tilde{\VWp}^{\msf',\asf'}_{\dsf'}(B,\tau)]_t. 
	\end{multline*}
	\end{conj}
Using the transformation relation (\ref{transform}), we deduce that under the change of variables $\tau \mapsto -1/\tau$, the generating series are exchanged in the following way (up to the factors and terms involving $\tau$), 
\begin{align}
	\begin{split} \label{transform2}	
	\CU_1(q) &\mapsto \CU_2(q)\\
\CU_2(q) &\mapsto \CU_1(q) \\
\CU_3(q) &\mapsto \CU_3(q).
\end{split}
\end{align}
 Using the transformations laws  (\ref{transform2}) and Conjecture \ref{conj1}, we obtain the following relation, which is true both for $\VWp^{\msf,\asf}_{\dsf}(B,q)$ and $\tilde{\VWp}^{\msf,\asf}_{\dsf}(B, q)$. For the statements we use the permutation action from Section \ref{conjsec},
 \[ \sigma\cdot \CU_i(q):= \CU_{\sigma(i)}(q), \quad \sigma \in S_3.\]
\begin{conj}[Geometric S-duality with $\mu$-insertions] \label{sdualityq} If $\rsf=2$, then 
	\[ (12)\cdot [\VWp^{\msf,\asf}_{\dsf}(B,q)]_t= 2^{-2g+1} \sum_{(\dsf',\msf',\asf')}(-1)^{(\dsf,\msf,\asf)\cdot (\dsf',\msf',\asf')}[\VWp^{\asf'}_{\dsf'}(B,q)]_t.\]
	\end{conj}

Using Theorem \ref{vwqm}, Definition \ref{defnnoncoprime} and Lemma \ref{difeq}, we thereby obtain Conjecture \ref{MS1}. 

The invariants for $\wsf=0$ similarly satisfy the $S$-duality transformation, which was also conjectured by \cite{JK}.

\begin{conj}[S-duality without $\mu$-insertions] If $\rsf\geq 2$, then
		\[\VWp^{\msf,\asf}_{\dsf,0}= \rsf^{-2g-1}  \sum_{(\dsf',\msf',\asf')}\omega^{ (\dsf,\msf,\asf)\cdot (\dsf',\msf',\asf') }\VWp^{\msf',\asf'}_{\dsf',0},\]
	where $\omega=e^{\frac{2\pi i}{\rsf}}$.
	\end{conj} 

Again using Theorem \ref{vwqm}, Definition \ref{defnnoncoprime} and Lemma \ref{difeq}, we thereby obtain Conjecture \ref{Enumerativemirror}.

\section{Some calculations for $\wsf=0$} \label{calc}
\subsection{$g=1,  \wsf=0$}
\begin{prop} Conjecture \ref{wzero} is true for $g=1$, i.e.\ if $g=1$, then 
	\begin{align*}
		&\QMs^0_{0,0}= \rsf^2+1 & &\QMp^0_{0,0}= \rsf+2-\rsf^{-1}\\
		&\QMs^0_{\dsf\neq0,0}=1  &  &\QMp^0_{\dsf\neq0,0}=\rsf+1-\rsf^{-1}  \\
		&\QMs^{\asf\neq0}_{0,0}=1 & &\QMp^{\asf\neq0}_{0,0}=\rsf+1-\rsf^{-1} \\
		&\QMs^{\asf\neq0}_{\dsf \neq 0,0}=0 & &\QMp^{\asf\neq0}_{\dsf\neq0,0}= \rsf-\rsf^{-1}.
	\end{align*}
	\end{prop}
\textit{Proof.} 
In fact, we will compute Vafa--Witten invariants $\VWp^{\msf,\asf}_{\dsf,0}$, which by Theorem \ref{vwqm} and Definition \ref{defnnoncoprime} determine quasimap invariants. Firstly, for every class $(\dsf,\msf,\asf) \in H^2(X\times E, \BZ)$, there exists an abelian variety $A$, which is deformation equivalent to $X\times E$ (in a projective family), such that with respect to the induced identification of cohomologies
\begin{equation} \label{cohiden}
H^2(X\times E,\BZ) \cong H^2(A, \BZ), 
\end{equation}
the class $(\dsf, \msf, \asf)$ is algebraic on $A$. We will denote $\mathrm{SU}_\rsf$ and $\mathrm{SU}_\rsf/\BZ_\rsf$ Vafa--Witten invariants on $A$ associated to the class $(\dsf, \msf, \asf) \in H^2(X\times E,\BZ)  \cong H^2(A, \BZ)$ and discriminant $\Delta=0$ by 
\[\VWs^{\msf,\asf}_{\dsf,0}(A) \quad \text{and} \quad \VWp^{\msf,\asf}_{\dsf,0}(A),\]
respectively. In particular, by deformation invariance,  $\mathrm{SU}_\rsf/\BZ_\rsf$ invariants  $\VWp^{\msf,\asf}_{\dsf,0}$ can be computed as $\mathrm{SU}_\rsf/\BZ_\rsf$ invariants  $\VWp^{\msf,\asf}_{\dsf,0}(A)$ on $A$, which by algebraicity of $(\dsf, \msf, \asf)$ on $A$, can be computed via $\VWs^{\msf,\asf}_{\dsf,0}(A)$, 
\[\VWp^{\msf,\asf}_{\dsf,0}=\VWp^{\msf,\asf}_{\dsf,0}(A)=\VWs^{\msf,\asf}_{\dsf,0}(A)/\rsf^4.\]
 We therefore have to compute  $\mathrm{SU}_\rsf$ Vafa--Witten invariants  on an abelian surface.  Firstly assume that $(\dsf, \msf, \asf) \neq (0,0,0)$, i.e.\ the corresponding moduli spaces do not have strictly semistable sheaves. By \cite[Proposition 7.4]{TT}, the stable $\mathrm{SU}_\rsf$ Vafa--Witten invariants on $A$ are just Euler characteristics of moduli spaces of sheaves with fixed determinants  in the class $(\dsf, \msf, \asf)$ and discriminant $\Delta=0$ (the moduli spaces are smooth, because there are no strictly semistable sheaves). The dimension of these moduli spaces is equal to 
\[ \sum (-1)^i \mathrm{ext}^i(F,F)=\mathrm c_2(R\CE nd(F,F))=\Delta(F)=0,\] 
 where the first equality follows from Hirzebruch--Riemann--Roch theorem and the fact the tangent bundle of $A$ is trivial. Hence there moduli spaces 
 are just finitely many points, if they are non-empty. By \cite[Lemma 4.19]{Y4}, there are $\rsf^2$ of these points. Hence 
\[\VWs^{\msf,\asf}_{\dsf,0}(A)=\rsf^2,\]
if $(\dsf, \msf, \asf) \neq (0,0,0)$ and the associated moduli space is non-empty. Since the moduli spaces of semistable sheaves on abelian surfaces are non-empty, if and only if the expected dimension is non-negative, the only constraint for our moduli spaces to be non-empty is compatibility of $\Delta=0$  with a choice of $(\dsf, \msf, \asf)$.  By the definition of discriminant $\Delta$, 
\[\Delta/2 =\mathrm c_1^2/2-\rsf \mathrm{ch}_2=\dsf  \asf +\msf ^2/2 -\rsf \mathrm{ch}_2,\]
hence the constraint that we obtain is the following one 
\begin{equation} \label{constraint}
 \dsf  \asf +\msf ^2/2 =0 \mod \rsf,
 \end{equation}
and this is the only constraint (because we are free to choose an arbitrary $\ch_2$). Hence if $(\dsf, \msf, \asf)$ satisfies (\ref{constraint}), then the moduli space of sheaves in the class $(\dsf, \msf, \asf)$ and $\Delta=0$ is non-empty. 
\\

We now deal with the case  $(\dsf, \msf, \asf) =(0,0,0)$. For that we have consider another correspondence of invariants. Let $E'$ be a smooth connected genus 1 curve. Let $P$ be a nodal stable genus 1 curve. The arguments of \cite{MT} can be easily adopted to the case of abelian surfaces. Hence by degenerating $A\times E'$ to $A\times P$ and localising with respect to $\BC^*$-action on $P$, one can show that Vafa--Witten invariants on $A$ (i.e.\ local Donaldson--Thomas invariants on $A\times \BC$ ) are equal to reduced 2-dimensional Donaldson-Thomas invariants on the abelian threefold $A\times E'$. 
Let $\CM$ be the moduli space of semistable sheaves on $A\times E'$  in following class
\[ \iota_*(\rsf,(0,0,0),0)\in H^{\mathrm{ev}}(A\times E',\BZ),\]
where $\iota \colon A \hookrightarrow A\times E'$ is an embedding of a slice of $A$ and $(\rsf,(0,0,0),0)\in H^\mathrm{ev}(A,\BZ)$ is the class of rank $\rsf$ sheaves with $(\dsf, \msf, \asf) =(0,0,0)$ and $\Delta/2=0$.  Let 
\[\DT^{0,0}_{0,0}(A\times E') \in \BQ\] 
be the generalised Donaldson--Thomas invariant associated to the moduli space $[\CM_\CL/E']$, where $\CM_\CL \subset \CM$ is the locus sheaves with a fixed determinant $\CL$. 
 Then by the arguments of \cite{MT}, we have
\[\DT^{0,0}_{0,0}(A\times E')=\VWs^{0,0}_{0,0}(A).\] 
We therefore have to compute the invariants $\DT^{0}_{0,0}(A\times E')$. In \cite{GPT}, the generalised Donaldson--Thomas invariants were computed for different moduli spaces. More precisely, for $[\CM/\Jac(A)\times E']$, i.e.\ instead of fixing the determinant, the quotient  is taken with respect to $\Jac(A)\times E'$. Let $\tilde{\DT}^{0}_{0,0}(A\times E') \in \BQ$ be generalised Donaldson--Thomas invariants associated to $[\CM/\Jac(A)\times E']$. In \cite[Lemma 4.13]{GPT}, these invariants are calculated, 
\[\tilde{\DT}^{0,0}_{0,0}(A\times E')=1+\rsf^{-2}.\] 
The moduli spaces $[\CM_\CL/E']$ and $[M_\vsf/\Jac(A)\times E']$ are related by the action of $\Jac(A)[\rsf]$, 
\[ [[\CM_\CL/E']/\Jac A[\rsf]]\cong [M_\vsf/\Jac(A)\times E'].\]
 Since this relation also holds in the equivariant Hall algebra,  the associated generalised Donaldson--Thomas invariants are related as follows
\[\DT^{0,0}_{0,0}(A\times E')=\rsf^4\tilde{\DT}^{0,0}_{0,0}(A\times E'),\]
we therefore obtain that 
\[ \VWs^{0,0}_{0,0}(A)=\DT^{0,0}_{0,0}(A\times E')=\rsf^4+\rsf^2.\]

Let us now go back to quasimap invariants. On the $\mathrm{SL}_\rsf$ side, by Theorem \ref{vwqm} and Definition \ref{defnnoncoprime}, we just have to divide $\mathrm{SU}_\rsf$ Vafa--Witten invariants by $\rsf^2$. We therefore obtain that 
 \begin{align*}
 &	\QMs^0_{0,0}= \VWs^{0,0}_{0,0}(A)/\rsf^2=\rsf^2+1 \\
 & \QMs^0_{\dsf\neq0,0}= \VWs^{0,0}_{\dsf\neq 0,0}(A)/\rsf^2=1.
 	\end{align*}
 On the $\mathrm{PGL}_\rsf$ side, again by Theorem \ref{vwqm} and Definition \ref{defnnoncoprime}, we have to divide $\rsf^4$ and sum over $\msf$. For that we have to count the number of $\msf$, which satisfy the associated constraint (\ref{constraint}) for the chosen values of $\dsf$ and $\asf$. As in (\ref{npm}), 
 let 
 \begin{align*}
 &N_{\mathsf k}:= |\{ \msf \in H^1(X, \BZ_\rsf  )\otimes H^1(E, \BZ_\rsf) \mid \frac{\msf \cdot \msf }{2}=\mathsf k\}|. 
\end{align*}
After some simple calculations, one can show that
\begin{align*}
&N_{0}= \rsf^3+ \rsf^2-\rsf^{-1}\\
&N_{\mathsf k\neq 0}= \rsf^3 -\rsf. 
\end{align*}
Then using the constraint (\ref{constraint}), we obtain
 \begin{align*}
 & \QMp^0_{0,0}=\sum_\msf \VWp^{\msf,0}_{0,0}(A)=1+N_0 / \rsf^{2}=1+\rsf+1+\rsf^{-1} \\
  &\QMp^0_{\dsf\neq 0,0}=\sum_\msf \VWp^{\msf,0}_{\dsf,0}(A)=N_0 / \rsf^{2}= \rsf+1+\rsf^{-1} \\
  &\QMp^{\asf\neq 0}_{0,0}=\sum_\msf \VWp^{\msf,\asf }_{0,0}(A)=N_0 / \rsf^{2}=\rsf+1+\rsf^{-1}\\
  &\QMp^{\asf\neq 0}_{\dsf\neq 0,0}=\sum_\msf \VWp^{\msf,\asf}_{\dsf,0}(A)=N_{\asf\cdot\dsf} / \rsf^{2}=\rsf-\rsf^{-1}.
 \end{align*}
\qed
\begin{rmk} Note that if $\Delta\neq 0$ and $(\dsf,\msf, \asf)\neq 0$, then the associated moduli space is positive-dimensional and therefore has an extra smooth map to an abelian variety, which is given by the dual determinant. This makes the corresponding invariants zero. 
	\end{rmk}
\begin{rmk} In fact, the invariants $\QMs_{\dsf,0}$ and $ \QMp_{\dsf,0}$ can be easily computed, using the point of view of maps, because $\Ms(\dsf)$ is a point, and $\Mp(\dsf)\cong B\BZ_\rsf^2$. To compute the invariants for refined degrees in terms of maps, one has to acquire a better understanding of the extended degree in terms of the geometry of $\Mp(\dsf)$. In particular, it must be related to the decomposition of the inertia stack $\CI \Mp(\dsf)=\coprod_{\gamma\in \Gamma_X}[\Ms(\dsf)_\gamma/\Gamma_X]$.  As the proof of Proposition $\ref{ggreat2}$ implicitly shows, if $\rsf=2$, the invariants $\QMp^1_{1,0}$ compute the equivariant Euler characteristics of the complement $\CI\Mp(1) /\Mp(1)$. However, we do not how know to establish a statement like that in general. 
	\end{rmk}
\subsection{$g\geq 2,\wsf=0$}
\begin{prop} \label{ggreat2}If $g\geq 2$, then 
	\begin{align*}
		&\QMs^0_{\dsf\neq 0,0}=-\rsf^{4g-3} \\
		&\QMp_{\dsf \neq 0,0} = \sum_\asf \QMp^\asf_{\dsf\neq 0,0}= - \rsf^{2g-3} +
		(\rsf-1)(\rsf^{4g-3}-\rsf^{2g-3}).
		\end{align*}
	In particular, it agrees with Conjecture \ref{wzero2}. 
	\end{prop}
\textit{Proof.} To compute these invariants, we use the fact that  genus 1 degree 0 Gromov--Witten invariants are given by (orbifold) Euler characteristics of the (orbifold) target space. In the case of $\Ms(\dsf)$, this can be seen directly. We can describe the moduli space $Q^{0^+}_{E,0}(\hat{M}(\dsf),0)$ as follows 
\[Q^{0^+}_{E,0}(\Ms(\dsf),0) \cong   \Ms(\dsf),\]
such that under this identification the obstruction theory is given by the complex
\[T_{ \Ms(\dsf)}\oplus T_{\Ms(\dsf)}[-1],\]
therefore 
\[\QMs_{\dsf,0}= \int_{\Ms(\dsf)^{\BC_t^*}}e(T_{\Ms(\dsf)^{\BC_t^*}}) =e(\Ms(\dsf)^{\BC_t^*})=e(\Ms(\dsf)).\]
The Euler characteristics of $\Ms(\dsf)$ is computed in \cite[Corollary 3.5]{Mer}, 
 \[e(\Ms(\dsf))=-\rsf^{4g-3}.\]
 
 The case of $\Mp(\dsf)$ is more complicated and we will not even sketch an argument. For that we refer to \cite[Theorem 1.1]{Ts}, where it is shown\footnote{More precisely, in \cite{Ts} the result is stated for $\langle\tau_1(\mathbb 1)\rangle_{1,1,0} $, which is equal to 24$\langle [(E,0)]\rangle_{1,1,0} $.} that 
 \[ \QMp_{\dsf,0}=\int_{\CI\CI\Mp(\dsf)^{\BC_t^*}}e(T_{\CI\CI\Mp(\dsf)^{\BC_t^*}}),\]
 where $e(T_{\CI\CI\Mp(\dsf)})$ is the top Chern class of the tangent bundle of the double inertia stack. Then 
\begin{equation*}
 \int_{\CI\CI\Mp(\dsf)^{\BC_t^*}}e(T_{\CI\CI\Mp(\dsf)^{\BC_t^*}})= e_{\mathrm{eq}}(\CI\Mp(\dsf))=e(I\Mp(\dsf))=e_{\mathrm{orb}}(\Mp(\dsf)), 
 \end{equation*}
where $I\Mp(\dsf)$ is the coarse space of the inertia stack  $\CI\Mp(\dsf)$. 
We therefore have to compute the orbifold Euler characteristics of $\Mp(\dsf)$.  A detailed analysis of the inertia stack was conducted in \cite[Section 7]{HT}, 
\[ I\Mp(\dsf)= \Ms(\dsf)/\Gamma_X\cup \coprod_{\gamma\neq \mathrm{id} } \Ms(\dsf)_\gamma/\Gamma_X,\]
such that the fixed locus $\Ms(\dsf)_\gamma/\Gamma_X$ of an element $\gamma\neq \mathrm{id} \in \Gamma_X$ is a certain Prym variety, i.e.\ an abelian variety. By \cite[Section 8]{HT}, the action of $\Gamma_X$ on  $\Ms(\dsf)_\gamma$ is also explicit. In particular, the Euler characteristics of the quotients  $\Ms(\dsf)_\gamma/\Gamma_X$ can be computed from the formulas in the end of \cite[Proposition 8.2]{HT} by taking away the twist by the character $\rho$ (which is due to the fact that authors consider the twisted orbfiold cohomology) and substituting $u=1,v=1$. One then obtains 
\[e(\Ms(\dsf)_\gamma/\Gamma_X)=(\rsf-1)\rsf^{2g-3},\] 
and the result is independent of $\gamma\neq \mathrm{id}$.
On the other hand, by \cite[Corollary 1.1.1]{HRV}, we have
\[e(\Ms(\dsf)/\Gamma_X)=-\rsf^{2g-3}.\]
Since $|\Gamma_X|=\rsf^{2g}$, we  obtain that 
\begin{multline*}
\QMp_{\dsf,0}=e(I\Mp(\dsf))=e(\Ms(\dsf)/\Gamma_X)+(\rsf^{2g}-1)e(\Ms(\dsf)_\gamma/\Gamma_X)\\
=- \rsf^{2g-3} +
(\rsf-1)(\rsf^{4g-3}-\rsf^{2g-3}).
\end{multline*}
\qed

\section{Gromov--Witten invariants and Quot schemes} \label{sectionws}
\subsection{Graph space for $\mathrm{SL}_\rsf$-quasimaps}
\begin{defn}  We define
	\[GQ(\Ms(\dsf),\wsf)\]
	to be the moduli space of prestable quasimaps  $f \colon \p^1 \rightarrow \FMs(\dsf)$ of degree $\wsf$  which map to  $\Ms(\dsf)$ at $\infty \in \p^1$ (without identifications by automorphisms of $\p^1$).
\end{defn}
\begin{rmk} The prinicple difference between the spaces $GQ(\Ms(\dsf),\wsf)$ and $Q^{0^+}_{\p^1}(\Ms(\dsf),\wsf)$ from Definition \ref{fixedC} is that we include \textit{all} prestable quasimaps, even those that might be fixed by infinitely many automorphisms of $\p^1$. 
	\end{rmk}

Throughout this section, by a Higgs sheaf on $X\times \p^1$ will mean a pair $(F,\phi)$, where\footnote{Which is a slight abuse of the terminology, because $p_{X}^*\omega_X\ncong \omega_{X\times \p^1}$.} $\phi \in \Hom(F,F\otimes p_{X}^*\omega_X)$. By Lemma \ref{ul} and the fact that $\Jac(\p^1)=0$, there is no need to distinguish between the moduli spaces $\Ms_\usf(X\times \p^1,\dsf, \wsf)$ and $\Ms_\Lsf(X\times \p^1,\dsf, \wsf)$, hence we will denote them simply by $\Ms(X\times \p^1,\dsf, \wsf)$. By arguments from Lemma \ref{ident}, we have an open embedding,
\[GQ(\Ms(\dsf),\wsf) \hookrightarrow \Ms(X\times \p^1,\dsf, \wsf).\]
There exists an additional $\BC^*$-action on $GQ(\Ms(\dsf),\wsf)$ given by scaling the $\p^1$-factor,
\[\lambda [x,y]=[\lambda x,y], \quad \lambda \in \BC^{*},\]
with weight 1 at $0\in \p^1$.
To distinguished it from the $\BC^*$-action given by scaling the field $\phi$, we incorporate the equivariant parameters into the notation,  
\[\BC_z^*- \text{torus scaling of $\p^1$} , \quad \BC_t^* - \text{torus scaling of $\phi$.} \]

\begin{defn}
	Let 
	\[\Ws_{\dsf,\wsf} :=GQ(\Ms(\dsf),\wsf)^{\BC_z^*} \]
the $\BC_z^*$-fixed locus of $GQ(\Ms(\dsf),\wsf)$. Equivalently, we can define
\[\Ws_{\dsf,\wsf} \subset M(X\times \p^1,\dsf,\wsf)^{\BC_z^*} \]
as a locus of sheaves which are characterized by the following condition
\[ (F,\phi) \in \Ws_{\dsf,\wsf} \iff (F,\phi)_{\infty} \text{ is stable,}\]
where $(F,\phi)_{\infty}$ is the fiber of $(F,\phi)$ over $\infty \in \p^1$.
\end{defn}
  By the definition of $\Ws_{\dsf,\wsf}$, there exists an  evaluation map 
\[\evsf\colon W_{\dsf,\wsf} \rightarrow \Ms(\dsf), \quad (F,\phi) \mapsto (F,\phi)_{\infty}. \]
To pushforward classes without localisation with respect to the $\BC^*_t$-action, we need the following lemma. 
\begin{lemma} \label{evproper}
	The evaluation map $\evsf$ is proper. 
\end{lemma}
\textit{Proof.} The moduli space $\Ms(X\times \p^1,\dsf,\wsf)$ has a natural compactification, denoted by $\overline M(X\times \p^1,\dsf,\wsf)$, which is given by the moduli space of 2-dimensional sheaves on $\p(\omega_X\oplus \CO_X) \times \p^1$.   This also provides a compactification of  $\Ws_{\dsf,\wsf}$, denoted by $\overline W_{\dsf,\wsf}$. The space $\overline W_{\dsf,\wsf}$ admits an evaluation map to the moduli space of 1-dimension sheaves on $\p(\omega_X\oplus \CO_X)$,  
\[\overline \evsf\colon \overline W_{\dsf,\wsf} \rightarrow \overline M(\dsf),\]
whose restriction over $\Ms(\dsf)\subset \overline M(\dsf)$ is $\evsf$. The map  $\overline \evsf$ is proper by properness of the moduli spaces, hence so is $\evsf$. 
\qed
\\

The moduli space $GQ(\Ws(\dsf),\wsf)$ carries a natural perfect obstruction theory, defined in the same way as the one in Section \ref{obstructionsection}. By \cite{GP}, the fixed part of the obstruction theory defines a $\BC^*_{t}$-equivariant virtual fundamental class
\[[\Ws_{\dsf,\wsf}]^{\mathrm{vir}} \in H^{\BC^*_{t}}_*(\Ws_{\dsf,\wsf},\BQ),\]
while the moving part defines the virtual normal bundle $N^{\mathrm{vir}}_{\Ws_{\dsf,\wsf}}$ and therefore the $\BC_t^*\times \BC_z^*$-equivariant Euler class 
\[e_{\BC_t^*\times \BC_z^*}(N_{\Ws_{\dsf,\wsf}}^{\mathrm{vir}})\in H_{\BC_t^*}^*(\Ws_{\dsf,\wsf},\BQ)[z^{\pm}].\]
We now in the position to make the definition of invariants that will contribute to the wall-crossing formulas. 
\begin{defn} \label{Ifunction} We define \textit{I-function}  
	\[I(q,z)=1+\sum_{\wsf>0}q^\wsf\evsf_*\left( \frac{[ \Ws_{\dsf,\wsf}]^{\mathrm{vir}}}{e_{\BC_t^*\times \BC_z^*}(N_{\Ws_{\dsf,\wsf}}^{\mathrm{vir}})}\right)\in H_{\BC_t^*}^{*}(\Ms(\dsf),\BQ)[z^\pm]\otimes\BQ[\![q]\!].\]
	We also define
	\[\mu(z):= [zI(q,z)-z]_{+}\in H_{\BC_t^*}^{*}(\Ms(\dsf),\BQ)[z]\otimes\BQ[\![q]\!],\] 
	where $[\dots ]_{+}$ is the truncation given by taking only non-negative powers of $z$. Let 
	\[\mu_{\wsf}(z)\in H_{\BC_t^*}^{*}(\Ms(\dsf),\BQ)[z]\] 
	be the coefficients of $q^{\wsf}$ in $\mu(z)$. 
\end{defn}
As in \cite{N}, $\mu(z)$ will appear in a wall-crossing formula between quasimap invariants for different values of $\epsilon$. 
\subsection{Flags of shaves and $\mathrm{SL}_\rsf$-quasimaps} \label{sectionflags}
\subsubsection{Graph space and flags of sheaves} 
We will now express $\Ws_{\dsf,\wsf}$ in terms of flags of sheaves.  The following will be largely a repetition of \cite{Ob}. For this section, we will treat Higgs sheaves in $M(X\times \p^1,\dsf,\wsf)$ as 2-dimensional sheaves on $\rT^*X\times \p^1$. 

We write $\p^1$ in terms the standard affines charts around $\infty$ and $0$,
\[\p^1=\BA^1_\infty \cup \BA^1. \]
Let $\CF  \in \Ws_{\dsf,\wsf}$, the sheaf $\CF$ is fixed by the $\BC_z$-action and is stable over $\infty$. This implies that
\[\CF_{|\rT^*X\times \BA^1_{\infty}}=\pi_{\rT^*X}^*\CG\]
where $\CG=\CF_{\infty}$ is fiber of $\CF$ over $\infty$. 
Since $\CF$ is fixed by the $\BC_z$-action, $\CF$ is naturally a $\BC_z$-equivariant sheaf. Hence over 0, the equivariant sheaf $\CF_{|\rT^*X\times \BA^1}$ can be identified with a graded $\CO_{\rT^*X}[x]$-module on $\rT^*X$, 
\[\CF_{|\rT^*X\times \BA^1}=\bigoplus_{i\geq i_0} \CF_i\mathbf z^i,\]
for some coherent sheaves $\CF_i$ on $\rT^*X$ and $i_0 \in \BZ$.  Here,  
\[\mathbf z:= \text{standard representation of } \BC_z^* \text{ on } \BC,\]
Since $\CF$ is a coherent sheaf, there exists $r$, such tat $\CF_r=\CF_{r+1}=\dots =\CG$.  Hence we have the following expression for our sheaf, 
\[\CF_{|\rT^*X\times \BA^1}=\CF_{i_0}\mathbf z^{i_0} \oplus \CF_{i_0+1}\mathbf z^{i_0+1} \oplus \dots \oplus  \CF_{r-1} \mathbf z^{r-1} \oplus \CG \mathbf z^{r} \oplus \dots  \]
Since $\CF$ is torsion free, the multiplication by $x$ gives injective morphisms,
\[f_{i}\colon \CF_i \hookrightarrow \CF_{i+1}.\]  
We thereby obtain a flag of sheaves associated to $\CG$, 
\[\CF_{\bullet}=(\CF_i \subseteq\CF_{i+1} \subseteq \dots\subseteq  \CF_{r-1}\subseteq \CF_r= \CG),\]
such that consecutive terms are allowed to be equal.  The space $W_{\dsf,\wsf}$ can therefore be decomposed in terms of the length of the associated flags.  Moreover, as it is shown in \cite[Proposition 4.1]{Ob}, the deformation theories of the flag $\CF_{\bullet}$ and the sheaf $\CF$ are equivalent. Most importantly, by the second-cosection argument of \cite{KKV}  and its modification \cite[Proposition 4.8]{Ob} to our set-up, only 1-step flags contribute to $\mu(z)$, as we will now explain.   

\subsubsection{Second cosection} By using Langton's inductive argument (see \cite[Theorem 2.B.1]{HL}), one can construct an injection\footnote{More precisely, one constructs an injection $x^k\CF \hookrightarrow   \pi^*_{\rT^*X}\CG$, where $x$ is the variable on $\BA^1$ and  $k\gg0$. We get the desired injection, since  $x^k\CF\cong\CF$.}
\[\CF \hookrightarrow   \pi^*_{\rT^*X}\CG.\]
Let $Q$ be its cokernel, then we can treat $\CF$ as a complex,
\[\CF=[\pi^*_{\rT^*X}\CG \rightarrow Q],\]
where  $Q$ is supported on $\rT^*X\times \Spec(\BC[x]/(x^m))$ and can be identified with 
\[Q=\CG/\CF_0 \oplus \CG/\CF_1 \mathbf z^1 \oplus \dots \oplus \CG/\CF_{m-1}\mathbf z^{m-1}.\]
Applying the construction of \cite[Section 5.4]{KKV}, one obtains two $\BC^*_{z}$-equivariant infinitesimal deformations 
\[\BC v_1 \oplus \BC v_2 \subseteq \mathbf z^{-1}\Ext_{\rT^*X\times \BP^1}^1 (\CF, \CF)_0,\]
such that $v_1$ corresponds to the deformation associated to the translation on $\BA^1$, while  $v_2$ corresponds to the deformation of  $\Spec(\BC[x]/(x^m))$ given by 
\[\{x^{m-1}(x-\epsilon)=0\} \subset \BA^1 \times \BA^1_\epsilon.\]
 These deformations can also be seen as deformations of $\CF$ viewed as a sheaf on $\rT^*X\times \p^1$ which is framed as $\infty \in \p^1$, i.e.\ we remember the identification $\CF_{\infty}\cong  \CG$, which is unique, because $\CG$ is stable. By \cite[Section 4]{HuyLe}, the deformation space and the obstruction space of framed sheaves are 
\[ \Ext_{\rT^*X\times \p^1}^1 (\CF, \CF(-\infty))_0, \quad  \Ext_{\rT^*X\times \p^1}^2 (\CF, \CF(-\infty))_0.\]
Now, with respect to the equivariant identification 
\[\omega_{\rT^*X \times \p^1} \cong (\mathbf{tz})^{-1} \CO_{\p^1}(-2\infty),\]
Serre's duality gives two $\BC^*_{t}$-equivariant cosections associated to these deformations 
\[  \Ext_{\rT^*X\times \BA^1}^2(\CF,\CF(-\infty))_0 \cong \mathbf t\Ext_{\rT^*X\times \BA^1}^1(\CF,\CF(-\infty))_0^\vee \rightarrow \mathbf t (\BC v_1 \oplus \BC v_2)^{\vee}.\]
These maps globalise to equivariant cosections of the obstruction theory
\begin{equation*}
	\sigma=(\sigma_1, \sigma_2) \colon h^1(\BE^{\mathrm{vir}}_{W_{\dsf,\wsf}}) \rightarrow \mathbf t \CO^{\oplus 2}_{W_{\dsf,\wsf}}, 
\end{equation*}
then \cite[Proposition 12]{KKV} translates into the following statement. 
\begin{prop} \label{vanishing}
	Let $\Ws_{\dsf,\wsf,s}\subset \Ws_{\dsf,\wsf}$ be a component, parametrising flags  $\CF_{\bullet}$ with $s$ non-trivial steps, i.e. $s=|\{i: \CF_i\neq \CF_{i+1}\}|$. If $s\geq 2$, then $\sigma_{|\Ws_{\dsf,\wsf,s}}$ is surjective, in particular,
	\[[\Ws_{\dsf,\wsf,s}]^{\mathrm{vir}}= t^2[\Ws_{\dsf,\wsf,s}]^{\mathrm{red}}, \]
	where 	$[\Ws_{\dsf,\wsf,s}]^{\mathrm{red}}$ is the \textit{reduced} virtual fundamental class associated to the cone of $\sigma$. 
	
	If $s=1$, only  $\sigma_{1| \Ws_{\dsf,\wsf,s}}$ is surjective, in particular, 
	\[[\Ws_{\dsf,\wsf,1}]^{\mathrm{vir}}= t[\Ws_{\dsf,\wsf,1}]^{\mathrm{red}}. \] 
\end{prop}
\textit{Proof.} We refer to \cite[Proposition 12]{KKV}, see also \cite[Proposition 4.8]{Ob}. 
\qed
\subsection{Quot schemes and $\mathrm{SL}_\rsf$-quasimaps} \label{quotschemes}

 The virtual dimension of the graph space $GQ(\Ms(\dsf),\wsf)$ is equal to 
\[ \dim(\Ms(\dsf)),\]
by the virtual localisation, it is also equal to the homological degree of the class 
\[\frac{[ \Ws_{\dsf,\wsf}]^{\mathrm{vir}}}{e_{\BC_t^*\times \BC_z^*}(N_{\Ws_{\dsf,\wsf}}^{\mathrm{vir}})}.\]
By Proposition \ref{vanishing}, we have
\[\frac{[ \Ws_{\dsf,\wsf}]^{\mathrm{vir}}}{e_{\BC_t^*\times \BC_z^*}(N_{\Ws_{\dsf,\wsf}}^{\mathrm{vir}})}= t\frac{[ \Ws_{\dsf,\wsf,1}]^{\mathrm{red}}}{e_{\BC_t^*\times \BC_z^*}(N_{\Ws_{\dsf,\wsf,1}}^{\mathrm{vir}})} + t^2\sum_{s\geq 2} \frac{[ \Ws_{\dsf,\wsf,s}]^{\mathrm{red}}}{e_{\BC_t^*\times \BC_z^*}(N_{\Ws_{\dsf,\wsf,s}}^{\mathrm{vir}})},\]
such that 
\[\frac{[ \Ws_{\dsf,\wsf,1}]^{\mathrm{red}}}{e_{\BC_t^*\times \BC_z^*}(N_{\Ws_{\dsf,\wsf,1}}^{\mathrm{vir}})} \quad  \text{and}  \quad \frac{[ \Ws_{\dsf,\wsf,s\geq 2}]^{\mathrm{red}}}{e_{\BC_t^*\times \BC_z^*}(N_{\Ws_{\dsf,\wsf,s\geq 2}}^{\mathrm{vir}})}\]
are of homological degree equal to $ \dim(\hat{M}(\dsf))+1$ and $\dim(\hat{M}(\dsf))+2$, respectively. Moreover, there are no poles in the variable $t$ by properness of the evaluation map Lemma \ref{evproper}. By degree constraints, we therefore obtain that
\begin{multline}  \label{theclass}
	\evsf_*\left( \frac{[ \Ws_{\dsf,\wsf}]^{\mathrm{vir}}}{e_{\BC_t^*\times \BC_z^*}(N_{\Ws_{\dsf,\wsf}}^{\mathrm{vir}})}\right)= \evsf_* \left(\frac{t[ \Ws_{\dsf,\wsf,1}]^{\mathrm{red}}}{e_{\BC_t^*\times \BC_z^*}(N_{\Ws_{\dsf,\wsf,1}}^{\mathrm{vir}})} \right)+ O(1/z^2)\\
	=\frac{\Qsf_{\dsf, \wsf}t  \mathbbm{1}}{z} + O(1/z^2)\in H_{\BC_t^*}^{*}(\hat{M}(\dsf),\BQ)[z^\pm], 
\end{multline}
where $\Qsf_{\dsf, \wsf}  \in \BQ$.  Hence for the purposes of determining $\mu(z)$, we only need to determine contributions of Quot schemes. 

\subsubsection{}
Following \cite[Section 4.5]{Ob}, we will now express $\Qsf_{\dsf,\wsf}$ in terms of virtual Euler characteristics of Quot schemes. Let 
\[\Ws^{m}_{\dsf,\wsf}\subset \Ws_{\dsf,\wsf}\]
be the component parametrising sheaves of the form
\[\CF_{|\rT^*X\times \BA^1}=\CK \mathbf z^{i_0} \oplus \CK\mathbf z^{i_0+1} \oplus \dots \oplus  \CK \mathbf z^{i_0+m-1} \oplus \CG \mathbf z^{i_0+m} \oplus \CG \mathbf z^{i_0+m+1}\oplus \dots   \]
Let $h_m$ be the unique integer such that 
\[h_m\rsf -\frac{\rk(\check{\wsf})}{m}\in [0, \rsf-1].\]
We define 
\[u_m:=h_m\vsf-\frac{\check{\wsf}}{m},\]
where $\check{\wsf}$ is defined in Lemma \ref{eqchern}. 
As it is explained in \cite[Section 4.4]{Ob}, the  flag $\CF_{\bullet}=(\CK\subset \CG)$ associated to $\CF$ is an element of a Quot scheme of $\CG$ with $\ch(\CG/\CK)=u_m$.  Let $\BG$ be the universal sheaf of $\Ms(\dsf)$ on $\Ms(\dsf)\times \rT^*X$. Then 
\[\Ws^{m}_{\dsf,\wsf} \cong \mathrm{Quot}(\BG/\Ms(\dsf),u_m).\]
To compute $\Qsf_{\dsf,\wsf}$, it is enough to evaluate the class at a $\BC^*$-fixed point in $\check{M}(\dsf)$.  Moreover, since we have the freedom to choose a $\BC^*$-fixed sheaf $\CG \in \check{M}(\dsf)$, we may just take a sheaf supported on the zero section of $\rT^*X$, i.e.\ a sheaf associated to a Higgs sheaf with zero Higgs field $(G,0)$. This allows us to express the class 
(\ref{theclass}) in terms of Quot schemes  $\mathrm{Quot}(G,u_m)$ on the curve $X$. 

However, as explained in detail in \cite{Nqm}, we have to endow $\mathrm{Quot}(G,u_m)$ with the obstruction theory of sheaves on $\rT^*X$ (not of sheaves on $X$), which nevertheless has a simple form.  Let $\CK$ and $\CQ$ be the universal kernel and the universal quotient on $\mathrm{Quot}(G,u_m)\times X$ and $\pi \colon \mathrm{Quot}(G,u_m)\times X \rightarrow X$ be the natural projection. 

\begin{lemma} \label{virtang}The virtual tangent bundle $\BT^{\mathrm{vir}}_{\mathrm{Quot}}$ of $\mathrm{Quot}(G,u_m)$ viewed as a Quot scheme on $\rT^*X$ is 
	\[\BT^{\mathrm{vir}}_{\mathrm{Quot}} \cong \CH om_\pi(\CK,\CQ)\oplus \CH om_\pi(\CK \boxtimes \omega^\vee_X,\CQ)\mathbf t[-1].\]
\end{lemma}
\textit{Proof.} See \cite[Lemma 5.1]{Nqm}. 
\qed 
\\

\noindent We define the virtual Euler characteristics of $\mathrm{Quot}(G,u_m)$  with respect to the $\BT^{\mathrm{vir}}_{\mathrm{Quot}}$ as 
\begin{multline*}
e^{\mathrm{vir}}(\mathrm{Quot}(G,u_m) )=[\mathrm{Quot}(G,u_m)]^{\mathrm{vir}}\cap \mathrm{c}_{\mathrm{vdim}}(\BT^{\mathrm{vir}}_{\mathrm{Quot}}) \\
=e_{\BC^*_t}( \CH om_\pi(\CK \boxtimes \omega^\vee_X,\CQ)\mathbf t) \cdot \left[\frac{\mathrm{c}(\mathrm{\CH om_\pi(\CK,\CQ)})}{\mathrm{c}_{\BC^*_t}({\CH om_\pi(\CK \boxtimes \omega^\vee_X,\CQ)\mathbf t})}\right]_{\mathrm{vdim}}. 
\end{multline*}
Concluding the preceding discussion, we obtain the following expression of $\Qsf_{\dsf,\wsf}$. 
\begin{thm} \label{wallcrossinv}
	Given a $\BC^*_t$-fixed Higgs sheaf $(G,0) \in \hat{M}(\dsf)$, then
	
	\[\frac{\Qsf_{\dsf,\wsf}t}{z}=\int_{[\Ws_{\dsf,\wsf}]^{\mathrm{vir}}}\frac{\evsf^*[(G,0)]}{e_{\BC_t^*\times \BC_z^*}(N_{\Ws_{\dsf,\wsf}}^{\mathrm{vir}})} = \sum_{m| \check{\wsf}}\frac{(-1)^{\chi(\check{\wsf}\cdot \vsf) }}{mz} e^{\mathrm{vir}}(\mathrm{Quot}(G,u_m)).\]
\end{thm}

\textit{Proof.} See \cite[Theorem 4.1]{Ob}.
\qed
\\

\noindent An immediate corollary of Lemma \ref{virtang}  is the  following vanishing result, which will be extremely useful for determining $\Qsf_{\dsf,\wsf}$ in some cases. 

\begin{cor} \label{vansh} Assume $g\geq 2$. If $u_m$ is not a class of a skyscraper sheaf, then 
	 \[ e^{\mathrm{vir}}(\mathrm{Quot}(G,u_m) )=0.\]
	\end{cor}
\textit{Proof.} 
Consider a quotient $Q$ of $G$ and its kernel $K$. By stability of $G$, we obtain that
\[\mathrm{ext}^1(K\otimes \omega^\vee_C, Q)=0.\]
Hence if $u_m$ is the class of a skyscraper sheaf, then 
\[\mathrm{hom}(K,Q)<\mathrm{hom}(K\otimes \omega^\vee_C,Q), \]
which implies that 
\[ \dim(\mathrm{Quot}(G,u_m))<\mathrm{vdim}(\mathrm{Quot}(G,u_m)),\]
therefore the Euler class vanishes, $e( \CH om_\pi(\CK \boxtimes \omega^\vee_X,\CQ)) =0$. 

Alternatively, this can be shown by realizing that $\omega_C$ is ample, therefore we have  $\mathrm{vdim}=\ch(K^\vee)\cdot \ch(Q)<0$ on $\rT^*X$. 
\qed
\subsection{Wall-crossing for $\mathrm{SL}_\rsf$-quasimaps}

The space $\BR_{>0} \cup \{0^+,\infty\}$ of $\epsilon$-stabilities is divided into chambers, inside of which the moduli space $Q^{\epsilon}_{g,n}(\Ms(\dsf),\wsf)$ stays the same, and as $\epsilon$ crosses  a wall between chambers, the moduli space changes discontinuously. Let $\epsilon_{0}=1/\wsf_{0}$ be a wall and $\epsilon^-$, $\epsilon^+$ be some values that are close to $\epsilon_{0}$ from left and right of the wall respectively. 

\begin{thm} \label{wallcrossing}Assuming $2g-2+n+\epsilon_0\wsf>0$, we have 
	\begin{multline*}
		\langle  \psi^{ m_1}\gamma_{1}, \dots, \psi^{m_{n}}\gamma_{n} \rangle^{\epsilon^{-}}_{g,n,\wsf}-\langle  \psi^{m_{1}}\gamma_{1}, \dots, \psi^{m_{n}}\gamma_{n} \rangle^{\epsilon^{+}}_{g,n,\wsf}\\
		=\sum_{k\geq 1} \langle \psi^{m_{1}}\gamma_{1}, \dots, \psi^{m_{n}}\gamma_{n},  \mu_{\wsf_0}, \dots,  \mu_{\wsf_0} \rangle^{\epsilon^{+}}_{g,n+k,\wsf-k\wsf_0}\\
		=\sum_{k\geq 1}\frac{(\Qsf_{\dsf, \wsf_0}t)^k}{k!}\langle \psi^{m_{1}}\gamma_{1}, \dots, \psi^{m_{n}}\gamma_{n},  \mathbbm{1}, \dots,  \mathbbm{1} \rangle^{\epsilon^{+}}_{g,n+k,\wsf-k\wsf_0}. 
	\end{multline*}
	
\end{thm}

\textit{Proof.} The wall-crossing formula is proven in \cite[Section 6]{YZ} in the GIT set-up. The proof does not require any modifications for invariants associated to quasimaps to moduli spaces of sheaves or equivariant invariants. See also \cite[Theorem 6.3]{N}. The last equality follows from (\ref{theclass}). 
\qed 
\\

We are interested in comparing quasimap invariants for $\epsilon=0^+$ and $\epsilon=\infty$.  Applying Theorem \ref{wallcrossing} inductively,  the wall-crossing difference between $\epsilon=0^+$ and $\epsilon=\infty$ is a sum of invariants of the form
\[ \frac{(\Qsf_{\dsf, \wsf_0}t)^k}{k!} \langle \psi^{m_{1}}\gamma_{1}, \dots, \psi^{m_{n}}\gamma_{n},  \mathbbm{1}, \dots,  \mathbbm{1} \rangle^{\infty }_{g,n+k,\wsf'},\]
which are trivial in many cases. For example, if we consider just primary insertions (no $\psi$-classes or  tautological classes coming from $\Mbar_{g,n}$ ), then by the string equation, 
\[ \langle \gamma_{1}, \dots, \gamma_{n},  \mathbbm{1}, \dots,  \mathbbm{1} \rangle^{\infty}_{g,n+k,\wsf'}=0,\]
unless $(\wsf', g, n)=(0,0,3)$ or $(0,1,1)$. Another source of vanishing results is the power of $t$ in the wall-crossing formula - if we know the degree of invariants with respect to $t$, then the wall-crossing contributions whose degree of $t$ is higher must therefore vanish. 


\begin{thm} \label{wsE} If $\wsf\neq 0$, then
	\begin{align*} 
		\QMs_{\dsf, \wsf}^{\bullet}&=\GWs_{\dsf,\wsf}^{\bullet}+\Qsf_{\dsf, \wsf}\cdot \QMs_{\dsf,0}\\
	&=\GWs_{\dsf,\wsf}^{\bullet}+\Qsf_{\dsf, \wsf}\cdot e(\Ms(\dsf)),
	\end{align*}
	where $e(\Ms(\dsf))$ is the topological Euler characteristics of $\Ms(\dsf)$.
\end{thm}
\textit{Proof.} The arguments of \cite[Section 6]{YZ} are applicable the invariants associated to the moduli space $Q^{\epsilon}_{E}(\Ms(\dsf),\wsf)^{\bullet, \BC^*}$,
	\[
\langle \emptyset\rangle^{\epsilon}_{E,0,\wsf} :=\int_{[Q^{\epsilon}_{E}(\Ms(\dsf),\wsf)^{\bullet, \BC^*}]^{\mathrm{vir}}}\frac{1}{e_{\BC^*}(N^{\mathrm{vir}})} \in\BQ t.\]
 see also \cite{NK3}. The resulting wall-crossing formula is 
\[\langle \emptyset\rangle^{\epsilon^-}_{E,0,\wsf} =\langle \emptyset\rangle^{\epsilon^+}_{E,0,\wsf}+\sum_{k\geq 1}\frac{(\Qsf_{\dsf, \wsf_0}t)^k}{k!}\langle  \mathbbm{1}, \dots,  \mathbbm{1} \rangle^{\epsilon^+}_{E,k,\wsf'},\]
where $\langle  \dots \rangle^{\epsilon}_{E,n,\wsf}$ are equivariant invariants associated to the moduli spaces of $\epsilon$-stable quasimaps from curves, whose smoothing is $E$ (i.e.\ $E$ with rational tails). 
Applying the wall-crossing formula above multiple times to go from $0^+$ to $\infty$ and dividing by $t$, the difference between $\QMs_{\dsf, \wsf}^{\bullet}$ and $\GWs_{\dsf, \wsf}^{\bullet}$ will be a sum of the terms of the form
\[\frac{\prod^{i=k}_{i=1}\Qsf_{\dsf, \wsf_i}t^{k-1}}{k!}\langle  \mathbbm{1}, \dots,  \mathbbm{1} \rangle^{\infty}_{E,k,\wsf-\sum \wsf_i}.\] 
Hence the only non-vanishing contributions in the wall-crossing formula above is therefore given by
\[\Qsf_{\dsf, \wsf} \langle \mathbb{1}\rangle_{E,1,0}^{\infty}\in \BQ , \]
because other summands vanish by the string equation.  The claim then follows from the identification of invariants,
\[\langle \mathbb{1}\rangle_{E,1,0}^{\infty}=\langle [E,0]\rangle_{1,1,0}^{\infty}=\QMs_{\dsf,0}.\]
\qed

\begin{cor}\label{gweqqm} Assume $\gcd(\rsf,\dsf)=1$ and $\gcd(\rsf,\wsf)=1$, then 
	\[\QMs_{\dsf, \wsf}^{\bullet}=\GWs_{\dsf,\wsf}^{\bullet}.\] 
	If $\gcd(\rsf,\wsf)\neq 1$, then only Quot schemes  parametrising $0$-dimensional quotients contribute to the wall-crossing.
	\end{cor}
\textit{Proof.} By Theorem \ref{wallcrossinv}, the wall-crossing invariants are expressed via virtual Euler characteristics of Quot schemes.  We will show that in this case, $u_m$ is never the class of a skyscraper sheaf. Recall that the definition of $u_m$, 
\[u_m:=h_m\vsf-\frac{\check{\wsf}}{m},\]
where $h_m$ is the unique integer such that 
 \[h_m\rsf -\frac{\rk(\check{\wsf})}{m}\in [0, \rsf-1].\]
Hence $u_m$  is a skyscraper, if and only if
\[ \rk(\check \wsf) = 0 \textrm{ mod } \rsf.\]
Now $\check \wsf$ is define by the system of equation of Lemma \ref{lemmaeqchern}. In particular, 
\[ \dsf \cdot \rk(\check \wsf) -\rsf \cdot\mathrm{c}_1(\check\wsf) =\wsf,\]
from which we obtain that 
\[  \dsf \cdot \rk(\check \wsf) =\wsf \textrm{ mod } \rsf.\]
By assumption this means that 
\[ \rk(\check \wsf) \neq 0 \textrm{ mod } \rsf,\]
hence $u_m$ is not the class of a skyscraper sheaf. The claim then follows from Theorem \ref{wsE} and the vanishing of Corollary \ref{vansh}. 
\qed 
\\

Let us say a few words about the punctorial Quot schemes. Consider the simplest case, $u_m=(0,1)$, then the Quot scheme is just the projectivization of $G$ 
\[\mathrm{Quot}(G,(0,1))=\p G.\]
Let $p \colon \p G \rightarrow X$ be the natural projection, then
\[\CH om_\pi(\CK \boxtimes \omega^\vee_X,\CQ)\cong T_{\p G}\otimes p^*\omega_X.\]
Hence 
\begin{align*}e^{\mathrm{vir}}(\p G)&=\left(\mathrm{c}_1(\p G)-\mathrm{c}_1(X)\right)\cdot \mathrm{c}_1(X) =e(\p G),
	\end{align*}
where the second equality follows from the sequence 
\[0\rightarrow T_p \rightarrow T_{\mathrm{Quot}} \rightarrow p^*\omega_X^\vee \rightarrow 0.\] 
By \cite[Section 1.7]{POp} and the fact that Euler characteristics does not depend on $\dsf$ (see Lemma \ref{dindependence}) , we  obtain that $e^{\mathrm{vir}}(\mathrm{Quot}(G(0,1)))=\rsf(2-2g)$. In particular, the wall-crossing invariants do not vanish in this case. We do not know how to compute these invariants for an arbitrary $u_m$, however, there is a simple $\dsf$-independence property that these invariants satisfy. 

\begin{lemma} \label{dindependence} If $\gcd(\rsf,\dsf)=1$ and $\gcd(\rsf,\dsf')=1$, then 
		\[ \Qsf_{\dsf, \wsf}=\Qsf_{\dsf', \wsf}. \]
	\end{lemma}

\textit{Proof.} If $\gcd(\rsf, \wsf)=1$, invariants vanish, hence assume $\gcd(\rsf, \wsf)\neq 1$. We use Lemma \ref{lemmaeqchern} and definition of $u_m$ to deduce that $u_m$ is independent of $\dsf$ for a fixed $\wsf$, if $\wsf$ is a multiple of $\rsf$. 

Next, we use the expression of these invariants in terms of Quot schemes from Theorem \ref{wallcrossinv}. One can deform any vector bundle on a curve to a direct sum of line bundles. Quot schemes of direct sums of lines bundles admit torus-actions, whose fixed loci are product of symmetric powers of the curve. In particular, the fixed loci  are independent of $\dsf$, see \cite[Section 3]{MR}. This can be done for sheaves scheme-theoretically supported on the zero section of $\mathrm{T}^*X$, such that the virtual Euler characteristics also localises (in this case, moduli spaces are the same, but the obstruction theories are given by Lemma \ref{virtang}). Since these fixed loci are independent of $\dsf$, we obtain the claim. \qed
\subsection{Graph space for $\mathrm{PGL}_\rsf$-quasimaps}

Since $\Mp(\dsf)$ is an orbifold, we need to consider graph spaces with $\p^1$ twisted at $\infty \in \p^1$. 
\begin{defn} Let $\p^1(m,1)$ be an orbifold $\p^1$ with an isotropy group $\BZ_m$ at $\infty$. In other words, around $\infty$, the curve $\p^1(m,1)$ is of the form $[\BA^1/\BZ_m]$, where  $\BZ_m$ acts on $\BA^1$ by multiplication by $m$-roots of unity. We define
	\[GQ(\Mp(\dsf),\wsf)\]
	to be the moduli space of prestable quasimaps  $f \colon \p^1(m,1) \rightarrow \HAz(\dsf)$ of degree $\wsf$ ($m$ is arbitrary) which map to  $\Mp(\dsf)$ at $\infty$.
\end{defn}
Recall that there is a natural decomposition of the inertia stack
\begin{equation} \label{decomp}
	\CI \Mp(\dsf)=\coprod_{\gamma \in \Gamma_X} [\Mp(\dsf)_\gamma/\Gamma_X],
	\end{equation}
where $\hat{M}(\dsf)_\gamma$ is the locus fixed by $\gamma \in \Gamma_X$. There exists a natural involution 
\[ \iota \colon \CI \Mp(\dsf)  \xrightarrow{\sim} \CI  \Mp(\dsf),\]
which with respect to the decomposition (\ref{decomp}) identifies $\Mp(\dsf)_\gamma$ with $\Mp(\dsf)_{\gamma^{-1}}$.

The graph space $GQ(\Mp(\dsf),\wsf)$ admits an evaluation map, which we define as a composition
\[\evsf \colon GQ(\Mp(\dsf),\wsf) \rightarrow \CI  \Mp(\dsf) \xrightarrow{\iota} \CI \Mp(\dsf).\]
\noindent The graph space decomposes according to (\ref{decomp}), 
\begin{equation} \label{decomposition2}
	GQ(\Mp(\dsf),\wsf)= \coprod_{\gamma \in \Gamma_X}  GQ^{\gamma}(\Mp(\dsf),\wsf), 
\end{equation}
where $GQ_{\gamma }(\Mp(\dsf),\wsf):= \evsf^{-1}([\Mp(\dsf)_\gamma/\Gamma_X])$.  

The graph space $GQ(\Mp(\dsf),\wsf)$ also admits a $\BC^*_t\times \BC^*_z$-action given by scaling of the target $\Mp(\dsf)$ and of the source curve $\p^1(m,1)$.  
\begin{defn}Let 
\[\Wp_{\dsf,\wsf}:= GQ (\Mp(\dsf),\wsf)^{\BC^*_z} \]
be the $\BC_z^*$-fixed locus of $GQ (\Mp(\dsf),\wsf)$.
\end{defn}
The space is endowed with the virtual fundamental class $[\Wp_{\dsf,\wsf}]^{\mathrm{vir}}$ and the virtual normal bundle $N_{\Wp_{\dsf,\wsf}}^{\mathrm{vir}}$ defined by the fixed and the moving parts of the obstruction theory of $GQ(\Mp(\dsf),\wsf)$, which is constructed in the same way as the one given in Proposition \ref{obst2}.

\begin{defn} \label{Ifunction2} We define \textit{I-function}  
	\[I(q,z)=1+\sum_{\wsf>0}q^\wsf\evsf_*\left( \frac{[ \Wp_{\dsf,\wsf}]^{\mathrm{vir}}}{e_{\BC_t^*\times \BC_z^*}(N_{W_{\dsf,\wsf}}^{\mathrm{vir}})}\right)\in H_{\mathrm{orb}, \BC_t^*}^{*}(\Mp(\dsf),\BQ)[z^\pm]\otimes\BQ[\![q]\!].\]
	We also define
	\[\mu(z):= [zI(q,z)-z]_{+}\in H_{orb, \BC_t^*}^{*}(\check{M}(\dsf),\BQ)[z]\otimes\BQ[\![q]\!],\] 
	where $[\dots ]_{+}$ is the truncation given by taking only non-negative powers of $z$. Let 
	\[\mu_{\wsf}(z)\in H_{\mathrm{orb}, \BC_t^*}^{*}(\check{M}(\dsf),\BQ)[z]\] 
	be the coefficients of $q^{\wsf}$ in $\mu(z)$. 
\end{defn}

The space $\Wp_{\dsf,\wsf}$ inherits the decomposition from (\ref{decomposition2}),
\[\Wp_{\dsf,\wsf}= \coprod_{\gamma \in \Gamma_X} \Wp_{\dsf,\wsf,\gamma }.\]
 We will now show that only  $\Wp_{\dsf,\wsf,\mathrm{id}}$ contribute to the $\mu$-function.
\begin{prop} \label{vanishingL}
	If $g \neq \mathrm{id}$, then 
	\[\evsf_*\left( \frac{[\Wp_{\dsf,\wsf,\gamma}]^{\mathrm{vir}}}{e_{\BC_t^*\times \BC_z^*}(N_{\Wp_{\dsf,\wsf,\gamma}}^{\mathrm{vir}})}\right)_{z^{\geq -1}}=0.\]
\end{prop}
\textit{Proof.} The claim follows from a virtual-dimension count. Firstly,  by \cite[Remark 1.6]{FG},
\begin{equation} \label{ineq1} \mathrm{age}(\gamma)+\mathrm{age}(\gamma^{-1})=\mathrm{dim}(\hat{M}(\dsf))-\mathrm{dim}(\hat{M}(\dsf)_\gamma).
\end{equation}
Secondly, by the discussion in \cite[Section 7.5]{Pic} the virtual dimension of a $GQ_{\gamma}(\check{M}(\dsf),\wsf)$  is equal to 
\begin{equation} \label{ineq2}
	\dim(\hat{M}(\dsf)) -\mathrm{age}(\gamma).
\end{equation}
If $\gamma \neq \mathrm{id}$, then $\mathrm{age}(\gamma)>0$, therefore (\ref{ineq1}) and (\ref{ineq2}) imply that the homological degree of the class
\begin{equation} \label{twistedw}
	\evsf_*\left( \frac{[\Wp_{\dsf,\wsf,\gamma}]^{\mathrm{vir}}}{e_{\BC_t^*\times \BC_z^*}(N_{\Wp_{\dsf,\wsf,\gamma}}^{\mathrm{vir}})}\right)\in H^{\BC_t^*}_*(\Mp(\dsf)_\gamma,\BQ)[z^{\pm}]
\end{equation}
exceeds the top degree of $H^{\BC_t^*}_*(\hat{M}(\dsf)_\gamma,\BQ)$. Moreover, the obstruction theory of  $\Wp_{\dsf,\wsf,\gamma} \subset GQ_{\gamma}(\Mp(\dsf),\wsf)^{\BC^*_z}$ relatively to $[\hat{M}(\dsf)_\gamma/\Gamma_X]$ is given by the obstruction theory of flags. Indeed, once a quasimap $f \in \Wp_{\dsf,\wsf,\gamma}$ is restricted to $\BA^1 \subset \p^1$ (the complement of $\infty \in \p^1$), the analysis of Section \ref{sectionflags} applies applies, because the Brauer group of $X\times \BA^1$ is zero, hence a family $(\CA,\phi)_{|\BA^1}$ associated to $f_{|\BA^1}$ is given by $(R\CE nd(F,F),\phi)$ for some $\BC^*_z$-fixed sheaf $F$. The deformation theories of $(R\CE nd(F,F),\phi)$ and $(F,\phi)$ are equivalent. Hence by the arguments of Section \ref{sectionflags}, the class (\ref{twistedw}) is a multiple of the equivariant parameter $t$. Therefore the truncation of this class by powers of $z$ of degree at least $-1$ is zero, 
\[\evsf_*\left( \frac{[\Wp_{\dsf,\wsf,\gamma}]^{\mathrm{vir}}}{e_{\BC_t^*\times \BC_z^*}(N_{\Wp_{\dsf,\wsf,\gamma}}^{\mathrm{vir}})}\right)_{z^{\geq -1}}=0.\]
\qed

\subsection{Quot schemes and $\mathrm{PGL}_\rsf$-quasimaps} Let us analyse the space $\Wp_{\dsf,\wsf,\mathrm{id}}$.   By the arguments of Proposition \ref{ident2}, we obtain an open embedding 
\begin{equation} \label{embedding}
 GQ_{\mathrm{id}}(\Mp(\dsf),\wsf) \hookrightarrow \Mp_0(X\times \p^1, \dsf, \wsf), 
 \end{equation}
note that $H^1(\p^1,\BQ)=0$, therefore the middle class is always zero, $\msf=0$.  We will denote  $\Mp_0(X\times \p^1, \dsf, \wsf)$ simply by $\Mp(X\times \p^1, \dsf, \wsf)$.  By Lemma \ref{ul}, we obtain that 
\[[\Mp(X\times \p^1, \dsf, \wsf)\cong [\Ms(X\times \p^1, \dsf, \wsf)/\Jac(X)[r]], \]
which via the embedding (\ref{embedding}) implies that
\[\Wp_{\dsf,\wsf,\mathrm{id}}\cong [\Ws_{\dsf,\wsf}/\Gamma_X ]. \]
Moreover, under the identification above, the evaluation map 
$ \evsf \colon \Wp_{\dsf,\wsf,\mathrm{id}} \rightarrow \Mp(\dsf)$ is just a quotient of $\evsf \colon \Ws_{\dsf,\wsf} \rightarrow \Ms(\dsf)$ by $\Gamma_X$. Combining Proposition \ref{vanishing} with the analysis above, we obtain the following description of the $\mu$-function. 
\begin{cor}
	\[	\evsf_*\left( \frac{[ \Wp_{\dsf,\wsf}]^{\mathrm{vir}}}{e_{\BC_t^*\times \BC_z^*}(N_{\check W_{\dsf,\wsf}}^{\mathrm{vir}})}\right)
	=\frac{\Qsf_{\dsf, \wsf}t  \mathbbm{1}}{z} + O(1/z^2)\in H_{\BC_t^*}^{*}(\Mp(\dsf),\BQ)[z^\pm], \]
	where $\Qsf_{\dsf, \wsf}$ is as in Section \ref{quotschemes}.
\end{cor}
\subsection{Wall-crossing for $\mathrm{PGL}_\rsf$-quasimaps} 
As before, the space $\BR_{>0} \cup \{0^+,\infty\}$ of $\epsilon$-stabilities is divided into chambers, inside of which the moduli space $Q^{\epsilon}_{g,n}(\Mp(\dsf),\wsf)$ stays the same, and as $\epsilon$ crosses  a wall between chambers, the moduli space changes discontinuously. Let $\epsilon_{0}=1/\wsf_{0}$ be a wall and $\epsilon^-$, $\epsilon^+$ be some values that are close to $\epsilon_{0}$ from left and right of the wall respectively. 

\begin{thm} \label{wallcrossing2}Assuming $2g-2+n+\epsilon_0\wsf>0$, we have 
	\begin{multline*}
		\langle  \psi^{ m_1}\gamma_{1}, \dots, \psi^{m_{n}}\gamma_{n} \rangle^{\vee,\epsilon^{-}}_{g,n,\wsf}-\langle  \psi^{m_{1}}\gamma_{1}, \dots, \psi^{m_{n}}\gamma_{n} \rangle^{\vee,\epsilon^{+}}_{g,n,\wsf}\\
			=\sum_{k\geq 1} \langle \psi^{m_{1}}\gamma_{1}, \dots, \psi^{m_{n}}\gamma_{n},  \mu_{\wsf_0}, \dots,  \mu_{\wsf_0} \rangle^{\vee,\epsilon^{+}}_{g,n+k,\wsf-k\wsf_0}\\
		=\sum_{k\geq 1}\frac{(\Qsf_{\dsf, \wsf_0}t)^k}{k!}\langle \psi^{m_{1}}\gamma_{1}, \dots, \psi^{m_{n}}\gamma_{n},  \mathbbm{1}, \dots,  \mathbbm{1} \rangle^{\vee,\epsilon^{+}}_{g,n+k,\wsf-k\wsf_0}. 
	\end{multline*}

\end{thm}
\textit{Proof.} Similar to Theorem \ref{wallcrossing}. \qed 
\begin{thm} \label{wsE2} \
	\begin{align*} 
		\QMp_{\dsf, \wsf}^{\bullet}&=\GWp_{\dsf,\wsf}^{\bullet}+\Qsf_{\dsf, \wsf}\cdot \QMp_{\dsf,0} \\
		&=\GWp_{\dsf,\wsf}^{\bullet}+\Qsf_{\dsf, \wsf}\cdot e_{\mathrm{orb}}(\Mp(\dsf)),
		\end{align*}
where $e_{\mathrm{orb}}(\Mp(\dsf))$ is the orbifold topological Euler characteristics of $\Mp(\dsf)$.
\end{thm}
\textit{Proof.} By the arguments identical to those presented in Theorem \ref{wsE}, we obtain that the only non-zero contribution to the wall-crossing is given by
\[\Qsf_{\dsf, \wsf} \langle \mathbb{1}\rangle_{E,1,0}^{\vee,\infty}\in \BQ.\] 
The claim then follows from the simple identification of invariants, 
\[\langle \mathbb{1}\rangle_{E,1,0}^{\vee, \infty}=\langle [E,0]\rangle_{1,1,0}^{\vee,0^+}=\QMp_{\dsf,0},\]
where the second equality follows from \cite{Ts}, see also the proof of Proposition \ref{ggreat2}. 
\qed
\begin{cor}\label{gweqqm2} Assume $\gcd(\rsf,\dsf)=1$ and $\gcd(\rsf,\wsf)=1$, then 
	\[\QMp_{\dsf, \wsf}^{\bullet}=\GWp_{\dsf,\wsf}^{\bullet}.\] 
	If $\gcd(\rsf,\wsf)\neq 1$, then only Quot schemes parametrising $0$-dimensional quotients contribute to the wall-crossing.
\end{cor}
\textit{Proof.} See Corollary \ref{gweqqm}. \qed
\\

Moreover, combining wall-crossings with Lemma \ref{dindependence} and Conjecture \ref{independence}, we obtain \textit{quantum $\chi$-independence} for Gromov--Witten invariants. 
\begin{conjecture}\label{independence2} Assume $\gcd(\rsf,\dsf)=1$ and $\gcd(\rsf,\dsf')=1$, then 
	\begin{align*}
	\GWs_{\dsf,\wsf}^{\bullet}&=\GWs_{\dsf',\wsf}^{\bullet} \\
	\GWp_{\dsf,\wsf}^{\bullet}&=\GWp_{\dsf',\wsf}^{\bullet}.
	\end{align*}
	\end{conjecture}
\bibliographystyle{amsalpha}
\bibliography{Mirr}

\providecommand{\bysame}{\leavevmode\hbox to3em{\hrulefill}\thinspace}
\providecommand{\MR}{\relax\ifhmode\unskip\space\fi MR }
\providecommand{\MRhref}[2]{%
  \href{http://www.ams.org/mathscinet-getitem?mr=#1}{#2}
}
\providecommand{\href}[2]{#2}
\begin{thebibliography}{HMMS22}

\bibitem[ACV03]{ACV}
D.~Abramovich, A.~Corti, and A.~Vistoli, \emph{Twisted bundles and admissible
  covers}, Comm. Algebra \textbf{31} (2003), no.~8, 3547--3618.

\bibitem[AO17]{AO}
M.~Aganagic and A.~Okounkov, \emph{Quasimap counts and {Bethe} eigenfunctions},
  Mosc. Math. J. \textbf{17} (2017), no.~4, 565--600.

\bibitem[AP22]{ArP}
D.~Aranha and P.~Pstragowski, \emph{{The Intrinsic Normal Cone For Artin
  Stacks}}, arXiv:1909.07478 (2022).

\bibitem[BLS98]{BLS}
A.~Beauville, Y.~Laszlo, and C.~Sorger, \emph{The {Picard} group of the moduli
  of {{\(G\)}}-bundles on a curve}, Compos. Math. \textbf{112} (1998), no.~2,
  183--216.

\bibitem[BNR89]{BNR}
A.~Beauville, M.~S. Narasimhan, and S.~Ramanan, \emph{Spectral curves and the
  generalized theta divisor}, J. Reine Angew. Math. \textbf{398} (1989),
  169--179.

\bibitem[Bru96]{Brus}
R.~Brussee, \emph{The {Canonical} {Class} and the {{\(C^ \infty\)}}
  {Properties} of {K{\"a}hler} {Surfaces}}, New York J. Math. \textbf{2}
  (1996), 103--146.

\bibitem[CKM14]{CFKM}
I.~{Ciocan-Fontanine}, B.~{Kim}, and D.~{Maulik}, \emph{{Stable quasimaps to
  GIT quotients}}, {J. Geom. Phys.} \textbf{75} (2014), 17--47.

\bibitem[CP18]{SP}
S.~Chakraborty and A.~Paul, \emph{Picard group and fundamental group of the
  moduli of {Higgs} bundles on curves}, Complex Manifolds \textbf{5} (2018),
  146--149.

\bibitem[DLI22]{DI}
A.~Di~Lorenzo and G.~Inchiostro, \emph{{Degenerations of twisted maps to
  algebraic stacks}}, arXiv:2210.03806 (2022).

\bibitem[DN89]{DN}
J.~M. Dr{\'e}zet and M.~Narasimhan, \emph{Groupe de picard des vari{\'e}t{\'e}s
  de modules de fibr{\'e}s semi-stables sur les courbes alg{\'e}briques},
  Inventiones mathematicae \textbf{97} (1989), 53--94.

\bibitem[DP12]{DP}
R.~Donagi and T.~Pantev, \emph{Langlands duality for {Hitchin} systems},
  Invent. Math. \textbf{189} (2012), no.~3, 653--735.

\bibitem[DPS98]{DPS}
R.~Dijkgraaf, J.-S. Park, and B.~Schroers, \emph{{N=4 Supersymmetric Yang-Mills
  Theory on a K\"ahler Surface}}, arXiv:hep-th/9801066 (1998).

\bibitem[FG03]{FG}
B.~Fantechi and L.~G\"{o}ttsche, \emph{Orbifold cohomology for global
  quotients}, Duke Math. J. \textbf{117} (2003), no.~2, 197--227.

\bibitem[FM99]{FM}
R.~Friedman and J.~W. Morgan, \emph{Obstruction bundles, semiregularity, and
  {Seiberg}-{Witten} invariants}, Commun. Anal. Geom. \textbf{7} (1999), no.~3,
  451--495.

\bibitem[GK18]{GK1}
L.~G{\"o}ttsche and M.~Kool, \emph{Refined $\text{SU}(3)$ {Vafa}-{Witten}
  invariants and modularity}, Pure Appl. Math. Q. \textbf{14} (2018), no.~3-4,
  467--513.

\bibitem[GK20a]{GK}
\bysame, \emph{Virtual refinements of the {Vafa}-{Witten} formula}, Commun.
  Math. Phys. \textbf{376} (2020), no.~1, 1--49.

\bibitem[GK20b]{GK2}
\bysame, \emph{Virtual refinements of the {Vafa}-{Witten} formula}, Commun.
  Math. Phys. \textbf{376} (2020), no.~1, 1--49.

\bibitem[GKL21]{GKL}
L.~Göttsche, M.~Kool, and T.~Laarakker, \emph{{$\mathrm{SU}(r)$ Vafa-Witten
  invariants, Ramanujan's continued fractions, and cosmic strings}},
  arXiv:2108.13413 (2021).

\bibitem[GKW21]{GKW}
L.~G{\"o}ttsche, M.~Kool, and R.~A. Williams, \emph{Verlinde formulae on
  complex surfaces: $\text{K}$-theoretic invariants}, Forum Math. Sigma
  \textbf{9} (2021), 31.

\bibitem[GP99]{GP}
T.~Graber and R.~Pandharipande, \emph{Localization of virtual classes}, Invent.
  Math. \textbf{135} (1999), no.~2, 487--518.

\bibitem[GWZ20]{GWZ}
M.~Groechenig, D.~Wyss, and P.~Ziegler, \emph{Mirror symmetry for moduli spaces
  of {Higgs} bundles via {{\(p\)}}-adic integration}, Invent. Math.
  \textbf{221} (2020), no.~2, 505--596.

\bibitem[HL95]{HuyLe}
D.~Huybrechts and M.~Lehn, \emph{Framed modules and their moduli}, Int. J.
  Math. \textbf{6} (1995), no.~2, 297--324.

\bibitem[HL97]{HL}
D.~{Huybrechts} and M.~{Lehn}, \emph{{The geometry of moduli spaces of
  sheaves}}, vol. E31, Braunschweig: Vieweg, 1997.

\bibitem[HMMS22]{HMMS}
T.~Hausel, A.~Mellit, A.~Minets, and O.~Schiffmann, \emph{{$P=W$ via
  $\mathcal{H}_2$}}, arXiv:2209.05429 (2022).

\bibitem[HRV08]{HRV}
T.~Hausel and F.~Rodriguez-Villegas, \emph{Mixed {Hodge} polynomials of
  character varieties. {With} an appendix by {Nicholas} {M}. {Katz}.}, Invent.
  Math. \textbf{174} (2008), no.~3, 555--624.

\bibitem[HT03]{HT}
T.~Hausel and M.~Thaddeus, \emph{Mirror symmetry, {Langlands} duality, and the
  {Hitchin} system}, Invent. Math. \textbf{153} (2003), no.~1, 197--229.

\bibitem[Jia22]{Jiang}
Y.~Jiang, \emph{Counting twisted sheaves and {S}-duality}, Adv. Math.
  \textbf{400} (2022), 75.

\bibitem[JK22]{JK}
Y.~Jiang and M.~Kool, \emph{Twisted sheaves and $\mathrm{SU}(r) / \mathbb{Z}_r$
  {Vafa-Witten theory}}, Math. Ann. \textbf{382} (2022), no.~1-2, 719--743.

\bibitem[Joy21]{J}
D.~Joyce, \emph{{Enumerative invariants and wall-crossing formulae in abelian
  categories}}, arXiv:2111.04694 (2021).

\bibitem[JS12]{JS}
D.~Joyce and Y.~Song, \emph{A theory of generalized {Donaldson}-{Thomas}
  invariants}, Mem. Am. Math. Soc., vol. 1020, Providence, RI: American
  Mathematical Society (AMS), 2012.

\bibitem[Kim22]{Kim}
Y.-J. Kim, \emph{{The dual Lagrangian fibration of known hyper-Kähler
  manifolds}}, arXiv:2109.03987 (2022.

\bibitem[KN97]{KN}
S.~Kumar and M.~S. Narasimhan, \emph{Picard group of the moduli spaces of
  {{\(G\)}}-bundles}, Math. Ann. \textbf{308} (1997), no.~1, 155--173.

\bibitem[KW07]{KW}
A.~Kapustin and E.~Witten, \emph{Electric-magnetic duality and the geometric
  {Langlands} program}, Commun. Number Theory Phys. \textbf{1} (2007), no.~1,
  1--236.

\bibitem[Lan75]{Lang}
S.~G. Langton, \emph{Valuative criteria for families of vector bundles on
  algebraic varieties}, Ann. Math. (2) \textbf{101} (1975), 88--110.

\bibitem[Lan15]{La}
A.~Langer, \emph{Bogomolov's inequality for {Higgs} sheaves in positive
  characteristic}, Invent. Math. \textbf{199} (2015), no.~3, 889--920.

\bibitem[Lie06]{Lie2}
M.~Lieblich, \emph{Moduli of complexes on a proper morphism}, J. Algebr. Geom.
  \textbf{15} (2006), no.~1, 175--206.

\bibitem[Lie07]{Lie3}
\bysame, \emph{Moduli of twisted sheaves}, Duke Math. J. \textbf{138} (2007),
  no.~1, 23--118.

\bibitem[Lie09]{Lie}
\bysame, \emph{Compactified moduli of projective bundles}, Algebra Number
  Theory \textbf{3} (2009), no.~6, 653--695.

\bibitem[Liu22]{Liu}
H.~Liu, \emph{{Equivariant K-theoretic enumerative invariants and wall-crossing
  formulae in abelian categories}}, arXiv:2207.13546 (2022).

\bibitem[LL99]{LL}
J.~M.~F. Labastida and C.~Lozano, \emph{The {Vafa}-{Witten} theory for gauge
  group {SU}{{\((N)\)}}}, Adv. Theor. Math. Phys. \textbf{3} (1999), no.~5,
  1201--1225.

\bibitem[Lur12]{Lur}
J.~Lurie, \emph{{Derived Algebraic Geometry XIV: Representability Theorems}},
  \url{https://www.math.ias.edu/~lurie/papers/DAG-XIV.pdf}.

\bibitem[Mer15]{Mer}
M.~Mereb, \emph{On the {\(e\)}-polynomials of a family of
  $\mathrm{SL}_n$-character varieties}, Math. Ann. \textbf{363} (2015),
  no.~3-4, 857--892.

\bibitem[MM21]{MM}
J.~Manschot and G.~W. Moore, \emph{{Topological correlators of $SU(2)$,
  $\mathcal{N}=2^*$ SYM on four-manifolds}}, arXiv:2104.06492 (2021).

\bibitem[MO09]{MO1}
D.~Maulik and A.~Oblomkov, \emph{Donaldson-{Thomas} theory of {{\(\mathcal
  A_n\times \mathbb{P}^{1}\)}}}, Compos. Math. \textbf{145} (2009), no.~5,
  1249--1276.

\bibitem[MO19]{MO}
D.~Maulik and A.~Okounkov, \emph{Quantum groups and quantum cohomology},
  Ast{\'e}risque, vol. 408, Paris: Soci{\'e}t{\'e} Math{\'e}matique de France
  (SMF), 2019.

\bibitem[MR22]{MR}
S.~Monavari and A.~T. Ricolfi, \emph{On the motive of the nested {Quot} scheme
  of points on a curve}, J. Algebra \textbf{610} (2022), 99--118.

\bibitem[MT18]{MT}
D.~Maulik and R.~P. Thomas, \emph{Sheaf counting on local {{\(K3\)}} surfaces},
  Pure Appl. Math. Q. \textbf{14} (2018), no.~3-4, 419--441.

\bibitem[Nes21a]{N}
D.~Nesterov, \emph{Quasimaps to moduli spaces of sheaves}, arXiv:2111.11417
  (2021).

\bibitem[Nes21b]{NK3}
\bysame, \emph{{Quasimaps to moduli spaces of sheaves on a K3 surface}},
  arXiv:2111.11425 (2021).

\bibitem[Nes23]{Nqm}
\bysame, \emph{{On quasimap invariants of moduli spaces of Higgs bundles}},
  2023.

\bibitem[NO21]{NO}
D.~Nesterov and G.~Oberdieck, \emph{Elliptic curves in hyper-{K{\"a}hler}
  varieties}, Int. Math. Res. Not. \textbf{2021} (2021), no.~4, 2962--2990.

\bibitem[Obe18]{O}
G.~Oberdieck, \emph{On reduced stable pair invariants}, Math. Z. \textbf{289}
  (2018), no.~1-2, 323--353.

\bibitem[Obe19]{OP1}
\bysame, \emph{{Gromov-Witten theory of \(K3 \times \mathbb{P}^1\) and
  quasi-Jacobi forms}}, {Int. Math. Res. Not.} \textbf{2019} (2019), no.~16,
  4966--5011.

\bibitem[Obe21]{Ob}
\bysame, \emph{{Multiple cover formulas for K3 geometries, wallcrossing, and
  Quot schemes}}, arXiv:2111.11239 (2021).

\bibitem[OP10]{OP10b}
A.~Okounkov and R.~Pandharipande, \emph{Quantum cohomology of the {Hilbert}
  scheme of points in the plane}, Invent. Math. \textbf{179} (2010), no.~3,
  523--557.

\bibitem[OP16]{OPa}
G.~Oberdieck and R.~Pandharipande, \emph{Curve counting on {$K3\times E$}, the
  {I}gusa cusp form {$\chi_{10}$}, and descendent integration}, K3 surfaces and
  their moduli, Progr. Math., vol. 315, Birkh\"{a}user/Springer, 2016,
  pp.~245--278.

\bibitem[OP21]{POp}
D.~Oprea and R.~Pandharipande, \emph{Quot schemes of curves and surfaces:
  virtual classes, integrals, {Euler} characteristics}, Geom. Topol.
  \textbf{25} (2021), no.~7, 3425--3505.

\bibitem[OPT22]{GPT}
G.~Oberdieck, D.~Piyaratne, and Y.~Toda, \emph{Donaldson-{Thomas} invariants of
  abelian threefolds and {Bridgeland} stability conditions}, J. Algebr. Geom.
  \textbf{31} (2022), no.~1, 13--73.

\bibitem[OS19]{OS}
G.~Oberdieck and J.~{Shen}, \emph{{Reduced Donaldson-Thomas invariants and the
  ring of dual numbers}}, {Proc. Lond. Math. Soc. (3)} \textbf{118} (2019),
  no.~1, 191--220.

\bibitem[OT22]{OT}
Z.-C. Ong and M.-C. Tan, \emph{{Vafa-Witten Theory: Invariants, Floer
  Homologies, Higgs Bundles, a Geometric Langlands Correspondence, and
  Categorification}}, arXiv:2203.17115 (2022).

\bibitem[Pic21]{Pic}
R.~Picciotto, \emph{{Moduli of stable maps with fields}}, arXiv:2009.04385
  (2021).

\bibitem[PT16]{KKV}
R.~Pandharipande and R.~P. Thomas, \emph{The {Katz}-{Klemm}-{Vafa} conjecture
  for {{\(K3\)}} surfaces}, Forum Math. Pi \textbf{4} (2016), 111.

\bibitem[Rua03]{Ru}
Y.~Ruan, \emph{Discrete torsion and twisted orbifold cohomology}, J. Symplectic
  Geom. \textbf{2} (2003), no.~1, 1--24.

\bibitem[ST15]{STV}
T.~Sch{\"u}rg and G.~To{\"e}n, B.and~Vezzosi, \emph{Derived algebraic geometry,
  determinants of perfect complexes, and applications to obstruction theories
  for maps and complexes}, J. Reine Angew. Math. \textbf{702} (2015), 1--40.

\bibitem[To{\"e}12]{To}
B.~To{\"e}n, \emph{Derived {Azumaya} algebras and generators for twisted
  derived categories}, Invent. Math. \textbf{189} (2012), no.~3, 581--652.

\bibitem[Tse11]{Ts}
H.-H. Tseng, \emph{On degree-0 elliptic orbifold {Gromov}-{Witten} invariants},
  Int. Math. Res. Not. \textbf{2011} (2011), no.~11, 2444--2468.

\bibitem[TT18]{TT2}
Y.~Tanaka and R.~P. Thomas, \emph{Vafa-{Witten} invariants for projective
  surfaces. {II}: {Semistable} case}, Pure Appl. Math. Q. \textbf{13} (2018),
  no.~3, 517--562.

\bibitem[TT20]{TT}
\bysame, \emph{Vafa-{Witten} invariants for projective surfaces. {I}: {Stable}
  case}, J. Algebr. Geom. \textbf{29} (2020), no.~4, 603--668.

\bibitem[TV07]{TV}
B.~To{\"e}n and M.~Vaqui{\'e}, \emph{Moduli of objects in dg-categories}, Ann.
  Sci. {\'E}c. Norm. Sup{\'e}r. (4) \textbf{40} (2007), no.~3, 387--444.

\bibitem[Wal85]{Wa}
F.~Waldhausen, \emph{Algebraic {K}-theory of spaces}, Algebraic and geometric
  topology, {Proc}. {Conf}., {New} {Brunswick}/{USA} 1983, {Lect}. {Notes}
  {Math}. 1126, 318-419 (1985)., 1985.

\bibitem[Wit10]{Wit}
E.~Witten, \emph{Mirror symmetry, {Hitchin}'s equations, and {Langlands}
  duality}, The many facets of geometry. A tribute to Nigel Hitchin, Oxford:
  Oxford University Press, 2010, pp.~113--128.

\bibitem[Yos96]{Yo}
K.~Yoshioka, \emph{Chamber structure of polarization and the moduli of stable
  sheaves on a ruled surface}, Int. J. Math. \textbf{7} (1996), no.~3,
  411--431.

\bibitem[Yos99]{Yo2}
\bysame, \emph{Some notes on the moduli of stable sheaves on elliptic
  surfaces}, Nagoya Math. J. \textbf{154} (1999), 73--102.

\bibitem[Yos01]{Y4}
\bysame, \emph{Moduli spaces of stable sheaves on abelian surfaces.}, Math.
  Ann. \textbf{321} (2001), no.~4, 817--884.

\bibitem[Yos06]{Yo3}
\bysame, \emph{Moduli spaces of twisted sheaves on a projective variety},
  Moduli spaces and arithmetic geometry. Papers of the 13th International
  Research Institute of the Mathematical Society of Japan, Kyoto, Japan,
  September 8--15, 2004, Tokyo: Mathematical Society of Japan, 2006, pp.~1--30.

\bibitem[Yun11]{Yun}
Z.~Yun, \emph{Global {Springer} theory}, Adv. Math. \textbf{228} (2011), no.~1,
  266--328.

\bibitem[Zag08]{Z}
D.~Zagier, \emph{Elliptic modular forms and their applications}, The 1-2-3 of
  modular forms. Lectures at a summer school in Nordfjordeid, Norway, June
  2004, Berlin: Springer, 2008, pp.~1--103.

\bibitem[{Zho}22]{YZ}
Y.~{Zhou}, \emph{{Quasimap wall-crossing for GIT quotients}}, {Invent. Math.}
  \textbf{227} (2022), no.~2, 581--660.

\end{thebibliography}
	\end{document}